\documentclass[10pt,a4paper]{article}
\usepackage[utf8]{inputenc}
\usepackage[T1]{fontenc}
\usepackage[english]{babel}
\usepackage[bottom=3cm, top=3cm, left=2cm, right=2cm]{geometry}
\usepackage{amsmath}
\usepackage{amsfonts}
\usepackage{amssymb}
\usepackage{amsthm}
\usepackage{graphicx}
\usepackage{color}
\usepackage{mathrsfs}
\usepackage{lmodern}
\usepackage[most]{tcolorbox}
\usepackage{framed}
\usepackage{fancybox}
\usepackage{tikz}
\usepackage{bbm}
\usepackage{dsfont}
\usepackage{wrapfig}
\usepackage{hyperref}
\usepackage[pagewise]{lineno}
\allowdisplaybreaks

\newtcolorbox{boxrd}[2][]{enhanced,colback=white,width={\textwidth},
attach boxed title to top left={yshift={-0.5\baselineskip},xshift=1cm}, 
title={#2},
boxrule=0.5pt,
coltitle=black,
boxed title style={
  borderline={-0.5mm}{black}
  colframe=white,
  colback=gray!50,
  colupper={black},
},
}

\newtcolorbox{boxsp}[2][]{%
  enhanced,colback=white,colframe=black,coltitle=black,
  boxrule=0.4pt,
  fonttitle=\itshape,
  attach boxed title to top left={yshift=-0.5\baselineskip-0.3pt,xshift=2mm},
  boxed title style={tile,size=minimal,left=0.5mm,right=0.5mm,
    colback=white,before upper=\strut},
  title=#2,#1
}

\newcommand{\N}{\mathbb{N}}

\newcommand{\R}{\mathbb{R}}
\newcommand{\C}{\mathbb{C}}

\newcommand{\HH}{\mathbf{H}}

\newcommand{\Ree}{\, \text{Re}}
\newcommand{\Imm}{\, \text{Im}}

\newcommand{\sgn}{\text{sgn}}

\newcommand{\midspace}{\, \, \, \, \, \, \, \, \, }

\usepackage{fancyhdr}
\usepackage{changepage}

\title{\textbf{Multi-solitary waves for the one-dimensional Zakharov system}}
\author{Guillaume Rialland}
\date{\footnotesize{Université de Paris-Saclay, UVSQ, CNRS, Laboratoire de Mathématiques de Versailles, 78000 Versailles \\ \texttt{guillaume.rialland@uvsq.fr}}}

\begin{document}
\maketitle

\begin{adjustwidth}{80pt}{80pt}
\small{\textsc{Abstract.} Given different speeds $c_1$ , ... , $c_K$, in the present paper we establish the existence of a solution to the Zakharov system in dimension $1$ that behaves asymptotically like a $K$-solitary wave, each wave travelling with speed $c_k$. The proof is adapted from previous results for the NLS and gKdV equations.}
\end{adjustwidth}

\textcolor{white}{a} \\ \\
\\ In this paper we study the Zakharov system
\begin{equation}
\left \{ \begin{array}{l} \partial_t u = i \partial_x^2 u - inu \\ \partial_t^2 n = \partial_x^2 n + \partial_x^2 ( |u|^2 ). \end{array} \right.
\label{Z}
\end{equation}
for $(t \, , x) \in \R \times \R$. The function $u$ is complex-valued while the function $n$ is real-valued. See \cite{Z} for the first introduction of this system by V. E. Zakharov to describe the propagation of Langmuir turbulence in a plasma ($u$ denotes the envelope of the electric field while $n$ denotes the deviation of the ion density from its equilibrium). One can also find in \cite{Su} a physical derivation of this system. We also refer to \cite{Kiv} for the physical interest of this system and a description of its solitary waves. \\
\\  The initial data for such a system is $(u(0) \, , n(0) \, , \partial_t n (0))$. The associated Cauchy problem is well-posed for $(u \, , n \, , \partial_t n)$ in $H^1 ( \R \, , \C ) \times L^2 ( \R \, , \R ) \times H^{-1} ( \R \, , \R )$. See \cite{Gin}, \cite{Bou}, \cite{Su} or \cite{Sa} for reviews of the wellposedness theory for the Zakharov system. Writing $\partial_t n \in H^{-1} ( \R )$ as $\partial_t n = - \partial_x v$ with $v \in L^2 ( \R )$ leads to the following equivalent system:
\begin{equation}
\left \{ \begin{array}{l} \partial_t u = i \partial_x^2 u - inu \\
\partial_t n = - \partial_x v \\
\partial_t v = - \partial_x n - \partial_x (|u|^2 ). \end{array} \right.
\label{ZL}
\end{equation}
See \cite{Gib} for the derivation of the system \eqref{ZL} from a Lagrangian formalism. As said above, the Cauchy problem associated with \eqref{Z} is well-posed for $(u \, , n \, , \partial_t n ) \in H^1 \times L^2 \times H^{-1}$; equivalently, the Cauchy problem associated with \eqref{ZL} is well-posed for $(u \, , n \, , v ) \in H^1  \times L^2 \times L^2$. \\
\\ The Zakharov system \eqref{Z} is the combination of a nonlinear Schrödinger equation (the first line) and a wave equation (the second line). If $n$ has small variation through time, $\partial_t^2 n \simeq 0$ thus $n \simeq - |u|^2$ and, consequently, the first line becomes the well-known focusing cubic nonlinear Schrödinger equation: $\partial_t u = i \partial_x^2 u + i |u|^2 u$. Unlike the latter equation, the Zakharov equation is not known to be integrable. \\
\\ The one-dimensional Zakharov system is $L^2$-subcritical. However, while both the nonlinear Schrödinger equation (NLS) and the wave equation have scale invariances, these two relevant scale invariances are incompatible and fail to provide a global scale invariance to the Zakharov system. We still have phase and translation invariances: if $(u \, , n )$ is a solution to \eqref{Z}, then $(u(t \, , x- \sigma ) e^{i \gamma} \, , n(t \, , x- \sigma ))$ is too, for any $\sigma$,$\gamma \in \R$. \\
\\ The following quantities are preserved through time by a solution $(u \, , v \, , n)$ to \eqref{ZL}:
\begin{itemize}
	\item the mass $M(u ) = \displaystyle{ \int_{\R} |u|^2}$;
	\item the energy $E(u \, , n \, , v) = \displaystyle{\int_{\R} \left ( | \partial_x u |^2 + n |u|^2 + \frac{n^2}{2} + \frac{v^2}{2} \right )}$;
	\item the momentum $P(u \, , n \, , v) = \displaystyle{\text{Im} \left (  \int_{\R} \overline{u} \partial_x u \right ) + \int_{\R} nv}$.
\end{itemize}
The presence of the function $v = - \partial_x^{-1} \partial_t n$ in the expressions of the energy and the momentum is one of the reasons why we prefer to work on the system \eqref{ZL} rather than on the system \eqref{Z}. \\
\\ In the present paper we study a family of solitary waves of \eqref{ZL}. Let $\phi_\omega (x) = \frac{\sqrt{2 \omega}}{\text{cosh} ( \sqrt{\omega} \, x)}$ be the solitary wave of pulsation $\omega > 0$ for the nonlinear Schrödinger equation: namely, 
\[ \partial_x^2 \phi_\omega = \omega \phi_\omega - \phi_\omega^3 . \]
Then $(u \, , n \, , v) = (e^{i \omega t} \phi_\omega (x) \, , - \phi_\omega^2 (x) \, , 0)$ provides a solitary wave of pulsation $\omega$ for the Zakharov system \eqref{ZL}. More generally, although no Galilean transform holds for the Zakharov equation (to our knowledge), \eqref{ZL} has a family of travelling waves, given by
\[ \left \{ \begin{array}{rcl} u_{\omega , c} (t \, , x) &=& \sqrt{1-c^2} \, \phi_\omega (x-ct) e^{i \left ( \frac{c   x}{2} - \frac{c^2 t}{4} + \omega t \right )} \\ 
n_{\omega , c} (t \, , x) &=& - \phi_\omega^2 (x-ct) \\
v_{\omega , c} (t \, , x) &=& - c \phi_\omega^2 (x-ct) \end{array} \right. \]
for $\omega > 0$ and $c \in ( -1 \, , 1)$. The speed of this travelling wave $(u_{\omega , c} \, , n_{\omega , c} \, , v_{\omega , c})$ is $c$. As proved in \cite{Wu} by spectral analysis or \cite{Oh} by variational arguments, these solitary waves for the Zakharov system (in dimension $1$) are stable, in the following sense: for any $\varepsilon > 0$ there exists $\delta > 0$ such that, if $(u_0 \, , n_0 \, , v_0) \in H^1 \times L^2 \times L^2$ satisfies
\[ \left \| (u_0 \, , n_0 \, , v_0) - \left ( u_{\omega , c} (0) \, , n_{\omega , c} (0) \, , v_{\omega , c} (0) \right ) \right \|_{H^1 \times L^2 \times L^2} \leqslant \delta , \]
then the solution $(u(t) \, , n(t) \, , v(t))$ of \eqref{ZL} with initial data $(u(0) \, , n(0) \, , v(0)) = (u_0 \, , n_0 \, , v_0 )$ satisfies
\[ \inf\limits_{\gamma , \sigma \in \R } \left \| (u(t) \, , n(t) \, , v(t)) - \left ( e^{i \gamma} u_{\omega , c} (t - \sigma ) \, , n_{\omega , c} (t - \sigma ) \, , v_{\omega , c} (t - \sigma ) \right ) \right \|_{H^1 \times L^2 \times L^2} \leqslant \varepsilon \]
for all $t \geqslant 0$. \\
\\ The main result of this paper is the existence of multi-solitary waves for the equation \eqref{ZL}. It is stated as Theorem 1 below. The existence of multi-solitons has been proved for several equations: see \cite{Ma1} for the generalized Korteweg-de Vries equation, \cite{Ma3}, \cite{Me2} for the NLS equation, \cite{Cot2} for the nonlinear Klein-Gordon equation (first occurrence of such a theorem for a wave-type equation) or \cite{Cot3}, \cite{Gus}, \cite{LeC1}, \cite{LeC2} for similar results for NLS-like equations. See \cite{Cot4} for the existence of multi-solitons for the supercritical gKdV and NLS equations.

\begin{leftbar}
\noindent \textbf{Theorem 1.} Let $K \in \N^*$. For all $k \in \{ 1 \, , ... \, , K \}$, take $\omega_k^0 > 0$, $c_k \in (-1 \, , 1)$, $\sigma_k^0 \in \R$ and $\gamma_k^0 \in \R$. Assume that $c_k \neq c_\ell$ for any $1 \leqslant k \neq \ell \leqslant K$. For all $k \in \{ 1 \, , ... \, , K \}$, define
\begin{align*}
& R_k^u = \sqrt{1-c_k^2} \, \phi_{\omega_k^0} (x-c_kt - \sigma_k^0) e^{i \left ( \frac{c_kx}{2} - \frac{c_k^2t}{4} + \omega_k^0 t + \gamma_k^0 \right )} , \\
& R_k^n = - \phi_{\omega_k^0}^2 (x-c_kt- \sigma_k^0 ) \\
\text{and} \, \, \, & R_k^v = - c_k \phi_{\omega_k^0}^2 (x-c_kt- \sigma_k^0 ). 
\end{align*}
Set $R^u = \sum\limits_{k=1}^K R_k^u$, $R^n = \sum\limits_{k=1}^K R_k^n$ and $R^v = \sum\limits_{k=1}^K R_k^v$. \\
\\ Then there exists $C>0$, $\theta_0 >0$ and a solution $(u(t) \, , n(t) \, , v(t)) \in \mathscr{C}^0 \left ( [0 \, , + \infty ) \, , H^2 \times H^1 \times H^1 \right )$ of \eqref{ZL} such that
\[ \left \| u(t) - R^u(t) \right \|_{H^2} + \left \| n(t) - R^n(t) \right \|_{H^1} + \left \| v(t) - R^v (t) \right \|_{H^1} \leqslant C e^{- \theta_0 t} \]
for all $t \geqslant 0$.
\end{leftbar}

\noindent \textbf{Remarks.}
\begin{itemize}
	\item This theorem ensures the existence of a particular solution $(u \, , n \, , v)(t)$ behaving like the multi-soliton $(R^u \, , R^n \, , R^v)$ as $t \to + \infty$. The error is exponentially decreasing (in norms $H^2 \times H^1 \times H^1$ here). Note that, since the Zakharov system is reversible (if $(u \, , n \, , v)$ is solution then so is $(\overline{u} (-t \, , x) \, , n(-t \, , x) \, , v(-t \, , x))$), Theorem 1 can be stated in a similar way for $t \to - \infty$. However, if one focuses on the solution behaving like a multi-soliton for $t \to + \infty$, its behavior for $t \to - \infty$ is unknown; there is no appearing reason for the solution to resemble a multi-soliton when $t \to - \infty$ as well. 
	\item Such solutions correspond to an exceptional behavior. Indeed, by strong $H^1$ convergence, one has
	\begin{align*}
	& M(u(t)) = \int_{\R} |u(t)|^2 = \sum\limits_{k=1}^K \int_{\R} |R_k^u(t)|^2 \\
	\text{and} \, \, \, & E(u(t) \, , n(t) \, , v(t)) = \int_{\R} \left ( | \partial_x u(t)|^2 + n(t) |u(t)|^2 + \frac{n(t)^2 + v(t)^2}{2} \right ) = \sum\limits_{k=1}^K E(R_k^u (t) \, , R_k^n (t) \, , R_k^v(t)). 
	\end{align*}
	This means that all the mass and all the energy of the solution are due to the solitary waves. This is an exceptional effect, since generally some of the mass is spread out by the dispersive effect of the equation.
	\item The use of higher-order modified energies would certainly enable to prove a better regularity on the solution ($\mathscr{C}^\infty$); such a purpose is not pursued here. See \cite{Cot1} for an analogous smoothness result for NLS multi-solitons and \cite{Ma1} for the gKdV equation. The modified energies estimated in section 4 are new (for the Zakharov system) and may have applications in a different framework, for example to prove global wellposedness in higher Sobolev spaces.
	\item The question of the uniqueness of this solution is an open problem, as well as the stability of multi-solitons. Similar issues are still open for the nonlinear Schrödinger equation (see for instance \cite{Ma4} for the stability of NLS multi-solitons in the case of particular nonlinearities). For the gKdV equation, the solution is known to be unique, see \cite{Ma1}. For the NLS equation, the solution is known to be unique in some sense, namely if the difference between the solution and the multi-soliton is asked to be \textit{sufficiently} decreasing (exponentially, or with a power sufficiently large). See \cite{Cot1} for this result. It is reasonable to think that a similar result holds for the Zakharov system; such a purpose is not pursued in this paper. 
	\item In dimension $d \geqslant 2$, the Zakharov system takes the form
	\begin{equation}
	\left \{ \begin{array}{l} \partial_t u = i \Delta u - inu \\
	\partial_t n = - \nabla \cdot v \\
	\partial_t v = - \nabla n - \nabla (|u|^2 ). \end{array} \right.
	\label{Z2}
	\end{equation}
	For example this system is well posed in $H^1 ( \R^2 \, , \C ) \times L^2 ( \R^2 \, , \R ) \times L^2 ( \R^2 \, , \R^2 )$. Denoting $\phi_\omega$ the solution of $\Delta \phi_\omega = \omega \phi_\omega - \phi_\omega^3$, then $(u \, , n \, , v ) = \left ( e^{i \omega t} \phi_\omega (x) \, , - \phi_\omega^2 (x) \, , 0 \right )$ provides a standing wave of pulsation $\omega$ for \eqref{Z2}. However, this solitary wave does not generate travelling waves. Actually, it seems that no expression is known for travelling waves in dimension higher than one. While $\phi_\omega$ is a radial function (and this is crucial in the proof), determining travelling waves in dimension higher than one is a non-radial problem. Note that no Galilean invariance holds for the Zakharov system (even in one dimension), contrary to the NLS equation.
\end{itemize}

\noindent \textbf{Outline of the proof.} The general outline of the proof presented in this paper is inspired from \cite{Ma1} and \cite{Ma3}, where analogous results are established for the gKdV equation or the nonlinear Schrödinger equation. The idea is to construct the solution \textit{backwards}: we consider a sequence $(u_p \, , n_p \, , v_p)$ of solutions to \eqref{ZL} with the convenient initial data $(u_p \, , n_p \, , v_p) (T_p) = (R^u \, , R^n \, , R^v)(T_p)$, where $T_p \to + \infty$. As $p \to + \infty$ (and $T_p \to + \infty$), we use compactness arguments to extract a subsequence $(u_{\psi (p)}\, , n_{\psi (p)} \, , v_{\psi (p)}) \underset{p \to + \infty}{\longrightarrow} \, (u \, , n \, , v)$ and we prove that $(u \, , n \, , v)$ is the desired solution behaving like the multi-soliton $(R^u \, , R^n \, , R^v)$. \\
\\ In section 1 we modulate the pulsations, translations and phases of the multi-soliton in order to gain orthogonality relations that will come in helpful for spectral properties. In section 2 we establish the pseudo-conservation (in some sense) of local energies and local momentums. In section 3 we analyse the so-called \textit{Weinstein functional}, namely a Lyapunov functional constructed in order to be quadratic in the difference $(u - R^u \, , n - R^n \, , v-R^v )$. Controlling this coercive functional gives good estimates on the norm $||(u-R^u \, , n-R^n \, , v-R^v)||^2$. In section 4, we extend these estimates to spaces of higher regularity via the use of modified energies. Such higher regularity is indeed helpful for the sake of the final compactness arguments and enables to rely on the Cauchy theory in the natural energy space $H^1 \times L^2 \times L^2$. An alternative approach (see \cite{Ma3} for NLS) is to use a more sophisticated Cauchy theory in a lower regularity space. Here, such wellposedness results exist (see \cite{Bou}, \cite{Gin}, \cite{Sa}) but require more delicate estimates in suitable spaces. We choose instead to obtain higher regularity estimates and use the natural Cauchy theory; as an additional consequence, we prove that the solution is $H^2 \times H^1 \times H^1$. Such an approach can be found in \cite{Ma1} for gKdV. \\
\\ \textbf{Notation.} We denote $\HH = H^1 ( \R \, , \C ) \times L^2 ( \R \, , \R ) \times L^2 ( \R \, , \R )$ with norm
\[ ||(u \, , n \, , v)||_{\HH} = ||u||_{H^1} + ||n||_{L^2} + ||v||_{L^2}. \]
We denote by $\langle f \, , g \rangle = \text{Re} \int_{\R} f \overline{g}$ the scalar product in $L^2$. In $\R^{3K}$ (the space of the set of parameters), we consider the norm
\[ |(a_1 \, , ... \, , a_K \, , b_1 \, , ... \, , b_K \, , c_1 \, , ... \, , c_K)| := \sum\limits_{k=1}^K \left ( |a_k| + |b_k| + |c_k| \right ) \]
and we denote by $\overline{B} ( \Pi \, , r )$ (respectively $\overset{o}{B} (\Pi \, , r)$) the closed (respectively open) ball centered at $\Pi$ and of radius~$r$. The notation $| \cdot |$ will also denote any matricial norm that we assume to be submultiplicative. We will often write $x_k^t$ instead of $(x-c_k t - \sigma_k (t))$. The letter $C$ will denote various positive constants whose expression change from one to another. The concerned constants do not depend on $t$, $x$, $p$, $A_0$ or $T_0$. Some information on the dependence of these constants $C$ can be found in some lemmas below. The notation $\mathcal{O} (g(t))$ will denote any function $f(t)$ such that $|f(t)| \leqslant C |g(t)|$ on the appropriate interval, where $C>0$ is a constant that verifies the remark above. \\
\\ \textbf{Acknowledgments.} This paper is the result of many discussions with Yvan Martel. May he be warmly thanked for it here.

\section{Modulation}
\noindent Let $K \in \N^*$. For all $k \in \{ 1 \, , ... \, , K \}$, take $\omega_k^0 > 0$, $c_k \in (-1 \, , 1)$, $\sigma_k^0 \in \R$ and $\gamma_k^0 \in \R$. Assume that $c_k \neq c_\ell$ for any $1 \leqslant k \neq \ell \leqslant K$. Hence we can assume that $c_1 < c_2 < \cdots < c_K$. For all $k \in \{ 1 \, , ... \, , K \}$, define
\begin{align*}
& R_k^u = \sqrt{1-c_k^2} \, \phi_{\omega_k^0} (x-c_kt - \sigma_k^0) e^{i \left ( \frac{c_kx}{2} - \frac{c_k^2t}{4} + \omega_k^0 t + \gamma_k^0 \right )} , \\
& R_k^n = - \phi_{\omega_k^0}^2 (x-c_kt- \sigma_k^0 ) \\
\text{and} \, \, \, & R_k^v = - c_k \phi_{\omega_k^0}^2 (x-c_kt- \sigma_k^0 ). 
\end{align*}
Set $R^u = \sum\limits_{k=1}^K R_k^u$, $R^n = \sum\limits_{k=1}^K R_k^n$ and $R^v = \sum\limits_{k=1}^K R_k^v$. Define
\[ \omega_- = \frac{1}{2} \min\limits_{1 \leqslant k \leqslant K} \omega_k^0 >0 \midspace \text{and} \midspace \omega^+ = \frac{3}{2} \max\limits_{1 \leqslant k \leqslant K} \omega_k^0 >0 \]
and define $\theta_0 > 0$ such that
\[ \sqrt{\theta_0} = \frac{1}{16} \min \left ( c_2 - c_1 \, , ... \, , c_K - c_{K-1} \, , \sqrt{\omega_-} \right ). \]
Let $T_p \to + \infty$ and $(u_p \, , n_p \, , v_p)$ the solution of \eqref{ZL} with the initial condition $(u_p \, , n_p \, , v_p)(T_p) = (R^u \, , R^n \, , R^v)(T_p)$. Set $(U_p \, , N_p \, , V_p) = (u_p - R^u \, , n_p - R^n \, , v_p - R^v)$. \\
\\ Let $A_0>1$ and $T_0>0$ large enough (depending on $A_0$; in a sense to be defined later). Let $p_0$ large enough such that $T_{p_0}>T_0$ and $p \geqslant p_0$. \\
\\ \textbf{In all of sections 1, 2 and 3} (until Proposition 2), we assume that there exists $t^* \in [T_0 \, , T_p]$ such that
\begin{equation} 
||(U_p (t) \, , N_p (t) \, , V_p (t))||_{\HH} \leqslant A_0 e^{- \theta_0 t} 
\label{bootstrap}
\end{equation}
for all $t \in [t^* \, , T_p]$. \textbf{In sections 1, 2, 3 and 4.1}, we denote $(u_p \, , n_p \, , v_p)$ simply by $(u \, , n \, , v)$ and $(U_p \, , N_p \, , V_p)$ simply par $(U \, , N \, , V)$ to lighten the computations. \\
\\ We begin with the modulation of the parameters $\omega_k$, $\sigma_k$ and $\gamma_k$ in order to gain orthogonality relations. The interest of such relations is the following. Take $Q(y) = \frac{\sqrt{2}}{\text{cosh} (y)}$ the soliton of the nonlinear Schrödinger equation ($Q'' = Q - Q^3$); we have $\phi_\omega (x) = \sqrt{\omega} \, Q ( \sqrt{\omega} \, x)$. Set $L_+ = - \partial_y^2 + 1 - 3Q^2$ and $L_- = - \partial_y^2 + 1 - Q^2$ the operators that appear when linearizing the equation \eqref{ZL} or the nonlinear Schrödinger equation around the soliton $Q$ (splitting $u$ into real and imaginary parts). Then, for any $\eta = \eta_1 + i \eta_2 \in H^1 ( \R )$, we recall the following standard spectral property from \cite{We}:
\begin{equation}
\langle L_+ \eta_1 \, , \eta_1 \rangle + \langle L_- \eta_2 \, , \eta_2 \rangle \geqslant C || \eta ||_{H^1}^2 - C \left ( \langle \eta_1 \, , Q \rangle^2 + \langle \eta_1 \, , yQ \rangle^2 + \langle \eta_2 \, , \Lambda Q \rangle^2 \right ) 
\label{spectral}
\end{equation}
where $\Lambda Q = \frac{1}{2} (Q+yQ')$. The operators $L_+$ and $L_-$ will appear in section 3 when expanding the Weinstein functional. As a consequence, the three orthogonalities above ($Q$, $yQ$ and $\Lambda Q$) will enable to prove the conercivity of the Weinstein functional, a key step of the proof. \\
\\ We recall here some properties of $Q$ and $\phi_\omega$. We already introduced $\Lambda Q = \frac{1}{2} ( Q + yQ' ) = \frac{1 - y \, \text{tanh}}{\sqrt{2} \, \text{cosh}}$. Also set $\Lambda_\omega (x) = \frac{\partial \phi_\omega}{\partial \omega} (x) = \frac{1}{\sqrt{\omega}} \Lambda Q ( \sqrt{\omega} \, x)$. It is known that $L_+(Q')=0$, $L_- Q = 0$, $L_- (yQ) = -2Q'$ and $L_+ ( \Lambda Q ) = -Q$ thus $L_- ( \Lambda Q ) = - Q + 2 Q^2 \Lambda Q$. See for example \cite{Ch} for such properties and studies of $L_+$, $L_-$ and $Q$.

\begin{leftbar}
\noindent \textbf{Proposition 1.} If $T_0 > 0$ is large enough, there exist $\mathscr{C}^1$ functions $\omega_k \, : \, [t^* \, , T_p] \to (0 \, , + \infty )$, $\sigma_k \, : \, [t^* \, , T_p] \to \R$ and $\gamma_k \, : \, [t^* \, , T_p] \to \R$ such that, if we set $\varepsilon (t) = (\varepsilon_u \, , \varepsilon_n \, , \varepsilon_v) (t \, , x) = (u-S^u \, , n-S^n \, , v-S^v) (t \, , x)$ where
\begin{align*}
& (S^u , S^n , S^v ) = \sum\limits_{k=1}^K (S_k^u \, , S_k^n \, , S_k^v ) , \\
& S_k^u (t \, , x) = \sqrt{1-c_k^2} \, \phi_{\omega_k(t)} (x - c_kt- \sigma_k(t)) e^{i \Gamma_k (t  , x)} , \\
& S_k^n (t \, , x) = - \phi_{\omega_k(t)}^2 (x -c_kt- \sigma_k (t) ), \\
& S_k^v (t \, , x) = - c_k \phi_{\omega_k(t)}^2 (x -c_kt- \sigma_k (t) ) \\
\text{and} \, \, \, & \Gamma_k (t \, , x) = \frac{c_k   x}{2} - \frac{c_k^2 t}{4} + \omega_k^0 t + \gamma_k (t),
\end{align*}
then the following orthogonality relations hold:
\begin{align*}
& \text{Re} \int_{\R} \phi_{\omega_k (t)} (x-c_kt- \sigma_k (t)) e^{-i \Gamma_k (t  , x)} \varepsilon_u (t \, , x) \, \text{d}x =0 , \\
& \text{Re} \int_{\R} (x-c_kt - \sigma_k (t)) \phi_{\omega_k (t)} (x-c_kt - \sigma_k (t)) e^{-i \Gamma_k (t  , x)} \varepsilon_u (t \, , x) \, \text{d}x = 0 \\
\text{and} \, \, \, & \text{Im} \int_{\R} \Lambda_{\omega_k (t)} (x-c_kt - \sigma_k (t)) e^{-i \Gamma_k (t  , x)} \varepsilon_u (t \, , x) \, \text{d}x= 0
\end{align*}
for all $t \in [t^* \,  , T_p]$ and $k \in \{ 1 \, , ... \, , K \}$. Moreover, for all $t \in [t^* \, , T_p]$,
\[ || \varepsilon (t)||_{\HH} + \sum\limits_{k=1}^K \left ( | \omega_k(t) - \omega_k^0 | + | \sigma_k (t) - \sigma_k^0 | + | \gamma_k (t) - \gamma_k^0 | \right ) \leqslant 2 A_0 e^{- \theta_0 t} \]
and 
\[ \forall k \in \{ 1 \, , ... \, , K \} , \, \, \, | \dot{\omega}_k(t)| + | \dot{\sigma}_k (t) | + | \dot{\gamma}_k(t) - (\omega_k (t) - \omega_k^0) | \leqslant C || \varepsilon (t) ||_{\HH} + C e^{- \theta_0 t}. \]
\end{leftbar}

\noindent \textit{Proof.} Analogous results hold for other equations (NLS or gKdV for example) and are standard. The proof here is no different. The complete proof is detailed in Appendix. \hfill \qedsymbol

\noindent \\ In all of sections 1 to 3, $T_0>0$ is taken large enough (depending on $A_0$) such that
\begin{equation}
\forall k \in \{ 1 \, , ... \, , K \} , \, \, \, \, \, \omega_- \leqslant \omega_k (t) \leqslant \omega^+. 
\label{omega}
\end{equation}

\section{Control of local quantities}
\noindent Let $\psi \in \mathscr{C}^3 ( \R \, , \R )$ such that $0 \leqslant \psi \leqslant 1$ on $\R$, $\psi ' \geqslant 0$ on $\R$, $\psi = 0$ on $(- \infty \, , -1]$, $\psi = 1$ on $(1 \, , + \infty )$ and $(\psi ')^2 \leqslant \psi$ and $(\psi '' )^2 \leqslant C \psi '$ on $\R$. For all $k \in \{ 2 \, , ... \, , K \}$, define $\overline{c}_k = \frac{c_{k-1}+c_{k}}{2}$. For $L>0$ large enough (in a sense to be defined later), let
\begin{align*}
& \chi_1 (t \, , x) = 1 - \psi \left ( \frac{x - \overline{c}_2 t}{L} \right ), \\
& \chi_K (t \, , x) = \psi \left ( \frac{x - \overline{c}_K t}{L} \right ) \\
\text{and} \, \, \, & \chi_k (t \, , x) = \psi \left ( \frac{x - \overline{c}_k t}{L} \right ) - \psi \left ( \frac{x - \overline{c}_{k+1} t}{L} \right ) \, \, \, \, \text{for $k \in \{ 2 \, , ... \, , K-1 \}$.}
\end{align*}
The number $L$ will depend on $A_0$ but not on $T_0$; $T_0$ will depend on $A_0$ and $L$. At last, for $k \in \{ 1 \, , ... \, , K \}$ let
\begin{align*} & \mathcal{M}_k (t) = \int_{\R} |u(t \, , x)|^2 \chi_k (t \, , x) \, \text{d}x  \\
\text{and} \, \, \, &  \mathcal{P}_k (t) = \text{Im} \int_{\R} \partial_x u (t \, , x) \overline{u (t \, , x)} \chi_k(t \, , x) \, \text{d}x + \int_{\R} n(t \, , x) v(t \, , x) \chi_k (t \, , x) \, \text{d}x.
\end{align*}
These are versions of the mass $M$ and the momentum $P$ localized near the axis $x= \overline{c}_k t$. They are \textit{almost} preserved through time, in the following sense (the exponential error term represents the interaction between solitons, that go further and further from one another as time gets larger). See \cite{Me} or \cite{Mz} for similar virial computations on the Zakharov system.

\begin{leftbar}
\noindent \textbf{Lemma 1.} There exists $C>0$ such that, for $L>0$ and $T_0>0$ large enough,
\[ \forall t \in [t^* \, , T_p] , \, \forall k \in \{ 2 \, , ... \, , K \} , \, \, \, \, | \mathcal{M}_k (T_p) - \mathcal{M}_k (t)| + | \mathcal{P}_k ( T_p ) - \mathcal{P}_k (t)| \leqslant \frac{C A_0^2}{L} e^{-2 \theta_0 t}. \]
\end{leftbar}

\noindent \textit{Proof.} We integrate by parts in order to get a convenient expression for $\dot{\mathcal{M}}_k (t)$:
\begin{align*}
\dot{\mathcal{M}}_k (t) =& \, \, 2 \Ree \int_{\R} (i \partial_x^2 u - inu) \overline{u} \chi_k  + \int_{\R} |u|^2 \partial_t \chi_k  \\
=& \, -2 \Imm \int_{\R} (\partial_x^2 u) \overline{u} \chi_k  + \underbrace{2 \, \text{Im} \int_{\R} n|u|^2 \chi_k }_{= \, 0} + \int_{\R} |u|^2 \partial_t \chi_k  \\
=& \, \, \underbrace{2 \Imm \int_{\R} (\partial_x u)   ( \partial_x \overline{u}) \chi_k }_{= \, 0} + 2 \, \text{Im} \int_{\R} ( \partial_x u) \overline{u}   \partial_x \chi_k + \int_{\R} |u|^2 \partial_t \chi_k  \\
=& \, \, 2 \Imm \int_{\R} ( \partial_{x} u) \overline{u} \partial_{x} \chi_k  + \int_{\R} |u|^2 \partial_t \chi_k.
\end{align*}
We proceed similarly for $\dot{\mathcal{P}}_k (t)$:
\begin{align*}
\dot{\mathcal{P}}_k(t) =& \, \, \text{Im} \int_{\R} (\partial_t \overline{u} \partial_x u - \partial_t u \partial_x \overline{u}) \chi_k  - \text{Im} \int_{\R} (\partial_t u) \overline{u} \partial_x \chi_k   + \text{Im} \int_{\R} \overline{u} (\partial_x u) \partial_t \chi_k   \\ & \midspace + \int_{\R} (\partial_t (nv) \chi_k  + nv \partial_t \chi_k  ) \\
=& \, \, \text{Im} \int_{\R} \left ( -i \partial_x ( | \partial_x u |^2 ) + i n \partial_x (|u|^2) \right ) \chi_k   - \text{Im} \int_{\R} (i \partial_x^2 u-inu) \overline{u} \partial_x \chi_k   \\
& \midspace + \int_{\R} ( \text{Im} (\overline{u} \partial_x u ) + nv) \partial_t \chi_k  - \int_{\R} \left ( \partial_x \left ( \frac{v^2+n^2}{2} \right ) + n \partial_x (|u|^2) \right ) \chi_k   \\
=& \int_{\R} \partial_x ( | \partial_x u|^2 ) \chi_k   + \int_{\R} n \partial_x (|u|^2) \chi_k   - \text{Im} \int_{\R} (-i \partial_x u) ( \partial_x \overline{u}   \partial_x \chi_k   + \overline{u} \partial_x^2 \chi_k  ) \\ 
& \midspace + \text{Im} \int_{\R} in |u|^2 \partial_x  \chi_k   + \int_{\R} ( \text{Im} ( \overline{u} \partial_x u) + nv ) \partial_t \chi_k   + \int_{\R} \frac{v^2+n^2}{2} \partial_x \chi_k ) \\ 
& \midspace + \int_{\R} |u|^2 (n \partial_x (\chi_k  ) + \chi_k \partial_x n) \\
=& \int_{\R} | \partial_x u|^2 \partial_x \chi_k  - \int_{\R} |u|^2 ( n \partial_x \chi_k   + \chi_k   \partial_x n) + \int_{\R} | \partial_x u |^2 \partial_x  \chi_k   \\
& \midspace + \text{Re} \int_{\R} \overline{u} ( \partial_x u ) \partial_x^2 \chi_k  + \int_{\R} n |u|^2 \partial_x  \chi_k   + \int_{\R} ( \text{Im} ( \overline{u} \partial_x u) + nv ) \partial_t \chi_k   \\ 
& \midspace + \int_{\R} \left ( \frac{v^2+n^2}{2} \right ) \partial_x  \chi_k   + \int_{\R} |u|^2 ( n \partial_x \chi_k + \chi_k \partial_x n ) \\
=& \, \, 2 \int_{\R} | \partial_x u |^2 \partial_x  \chi_k   + \text{Re} \int_{\R} \overline{u} (\partial_x u) \partial_x^2  \chi_k   + \int_{\R} \left ( n |u|^2 + \frac{v^2+n^2}{2} \right ) \partial_x  \chi_k   \\
& \midspace + \int_{\R} ( \text{Im} ( \overline{u} \partial_x u ) + nv ) \partial_t \chi_k  .
\end{align*}
Note that $\partial_{x} \chi_k  = 0$ for $x \not\in \left [ \overline{c}_k t - L \, , \overline{c}_k t + L \right ]$. Thus,
\begin{align*}
| \dot{\mathcal{M}}_k (t)| + | \dot{\mathcal{P}}_k(t)| \leqslant & \, \,\frac{C}{L} \int_{\overline{c}_kt-L}^{\overline{c}_kt+L} ( |\partial_x u|^2 + |u|^2 + n^2 + v^2) \\
\leqslant & \, \, \frac{C}{L} \left [ \int_{\overline{c}_kt-L}^{\overline{c}_kt+L} \left ( | \partial_x S^u |^2 + |S^u|^2 + (S^n)^2 + (S^v)^2 \right ) + ||\varepsilon_u||_{H^1}^2 + ||\varepsilon_n||_{L^2}^2 + ||\varepsilon_v||_{L^2}^2 \right ]. 
\end{align*}
We know that $| \phi_\omega ' (x)| + | \phi_\omega (x)| \leqslant C e^{- \sqrt{\omega} |x|}$, hence (recalling \eqref{omega})
\[ \int_{\overline{c}_kt-L}^{\overline{c}_kt+L} \left ( | \partial_x S^u |^2 + |S^u|^2 + (S^n)^2 + (S^v)^2 \right ) \leqslant C e^{-8 \sqrt{\theta_0} ( \sqrt{\theta_0} \, t - L)} \leqslant C e^{-4 \theta_0 t} \]
by taking $T_0$ large enough (depending on $L$) so that $T_0 \sqrt{\theta_0} \geqslant 2L$. Taking $A_0 e^{- \theta_0 T_0}$ small enough, we obtain
\[ | \dot{\mathcal{M}}_k (t) | + | \dot{\mathcal{P}}_k (t)| \leqslant \frac{C A_0^2}{L} e^{-2 \theta_0 t}. \]
Integrating between $t$ and $T_p$ we get the desired result.
\hfill \qedsymbol

\section{Weinstein functional}
\noindent Consider the Weinstein functional
\[ G(t) = \int_{\R} \left ( |\partial_x u|^2 + n |u|^2 + \frac{n^2+v^2}{2} \right ) + \sum\limits_{k=1}^K \nu_k^0 \int_{\R} |u|^2 \chi_k  - \sum\limits_{k=1}^K c_k   \left ( \text{Im} \int_{\R} \overline{u} (\partial_x u) \chi_k + \int_{\R} nv \chi_k \right ) \]
with $\nu_k^0 = \omega_k^0 + \frac{c_k^2}{4}$. Writing $u = S^u + \varepsilon_u$, $v = S^v + \varepsilon_v$ and $n=S^n+ \varepsilon_n$, we have $G=G_0+G_1+G_2+G_3$ where 
\begin{align*}
G_0 =& \int_{\R} \left ( | \partial_x S^u|^2 + S^n |S^u|^2 + \frac{(S^n)^2}{2} + \frac{(S^v)^2}{2} \right ) \\
& \midspace + \sum\limits_{k=1}^K \int_{\R} \left ( \nu_k^0 |S^u|^2 \chi_k - c_k   \text{Im} ( \overline{S^u} (\partial_x S^u) \chi_k ) - c_k   S^n S^v \chi_k \right ), \\
G_1 =& \int_{\R} \left ( 2 \, \text{Re} ( \partial_x S^u   \partial_x \overline{\varepsilon_u} ) + 2 S^n \text{Re} ( S^u  \overline{\varepsilon_u} ) + \varepsilon_n |S^u|^2 + \varepsilon_nS^n +  \varepsilon_vS^v \right ) + \sum\limits_{k=1}^K \nu_k^0 \int_{\R} 2 \,  \text{Re} (S^u \overline{\varepsilon_u} ) \chi_k \\
& \midspace - \sum\limits_{k=1}^K c_k   \int_{\R} \left ( \text{Im} ( \overline{S^u} (\partial_x \varepsilon_u) \chi_k + \overline{\varepsilon_u} ( \partial_x S^u) \chi_k ) + S^n  \varepsilon_v \chi_k + \varepsilon_n S^v \chi_k \right ), \\
G_2 =& \int_{\R} \left ( | \partial_x \varepsilon_u|^2 + 2 \varepsilon_n \, \text{Re} (S^u \overline{\varepsilon_u}) + S^n |\varepsilon_u|^2 + \frac{\varepsilon_n^2}{2} + \frac{| \varepsilon_v|^2}{2} \right ) + \sum\limits_{k=1}^K \nu_k^0 \int_{\R} |\varepsilon_u|^2 \chi_k \\
& \midspace - \sum\limits_{k=1}^K c_k   \int_{\R} \left ( \text{Im} ( \overline{\varepsilon_u} ( \partial_x \varepsilon_u ) \chi_k ) + \varepsilon_n \varepsilon_v \chi_k \right ) \\
\text{and} \, \, \, G_3 =& \int_{\R} \varepsilon_n|\varepsilon_u|^2.
\end{align*}
We aim at showing that the quantity $|G(t) - G(T_p)|$ is quadratic in $||( \varepsilon_u \, , \varepsilon_n \, , \varepsilon_v )||_{\mathbf{H}}$, with exponential error terms. Let us control each of the terms above individually. See \cite{We}, \cite{Ma3} and \cite{Ma4} for the similar use of a suitable Weinstein functional for NLS multi-solitons.

\subsection{Control of zero-order term $G_0$}
\noindent Using $\sum\limits_{k=1}^K \chi_k = 1$, we see that
\begin{align*}
G_0 =& \, \sum\limits_{k=1}^K \int_{\R} \left ( \left ( | \partial_x S^u|^2 + S^n |S^u|^2 + \frac{(S^n)^2}{2} + \frac{(S^v)^2}{2} \right ) \chi_k + \nu_k^0 |S^u|^2 \chi_k - c_k   \text{Im} ( \overline{S^u} (\partial_x S^u) \chi_k ) - c_k   S^n S^v \chi_k \right ) \\
=& \, \sum\limits_{k=1}^K \int_{\R} \left ( |\partial_x S_k^u|^2 + S_k^n |S_k^u|^2 + \frac{(S_k^n)^2}{2} + \frac{(S_k^v)^2}{2} + \nu_k^0 |S_k^u|^2 - c_k   \text{Im} \left ( \overline{S_k^u} \partial_x S_k^u \right ) - c_k   S_k^n S_k^v \right ) + \mathcal{O} (e^{-2 \theta_0 t}).
\end{align*}
We compute 
\[ \partial_x S_k^u = \sqrt{1-c_k^2} \, e^{i \Gamma_k(t \, , x)} \left ( \phi_{\omega_k} ' + \frac{i c_k}{2} \phi_{\omega_k} \right )(x_k^t) \, \, \, \, \, \text{and} \, \, \, \, \, | \partial_x S_k^u |^2 = (1-c_k^2) \left ( ( \phi_{\omega_k} ' )^2 + \frac{c_k^2}{4} \phi_{\omega_k}^2 \right ) (x_k^t). \]
We obtain
\[ G_0 (t) = \sum\limits_{k=1}^K (1-c_k^2) \left ( ( \phi_{\omega_k} ' )^2 + \omega_k^0 \phi_{\omega_k}^2 - \frac{1}{2} \phi_{\omega_k}^4 \right ) + \mathcal{O} \left ( e^{-2 \theta_0 t} \right ). \]
Integrating the equation $\phi_\omega '' = \omega \phi_\omega - \phi_\omega^3$, we have $( \phi_\omega ' )^2 = \omega \phi_\omega^2 - \frac{1}{2} \phi_\omega^4$. This leads to
\[ G_0 (t) = \sum\limits_{k=1}^K (1-c_k^2) \int_{\R} \left ( (\omega_k^0 - \omega_k ) \phi_{\omega_k}^2 + \frac{1}{2} \phi_{\omega_k}^4 \right ) + \mathcal{O} (e^{-2 \theta_0 t}). \]
Since $\omega_k (T_p) = \omega_k^0$, it follows that
\[ G_0 (t) - G_0 (T_p) = \sum\limits_{k=1}^K (1-c_k^2) \int_{\R} \left ( ( \omega_k^0 - \omega_k ) \phi_{\omega_k}^2 + \frac{1}{2} \phi_{\omega_k}^4 - \frac{1}{2} \phi_{\omega_k^0}^4 \right ) + \mathcal{O} ( e^{-2 \theta_0 t} ). \]

\begin{leftbar}
\noindent \textbf{Lemma 2.} For all $t \in [t^* \, , T_p]$,
\[ | G_0 (t) - G_0 (T_p) | \leqslant C \sum\limits_{k=1}^K | \omega_k (t) - \omega_k^0 |^2 + C e^{-2 \theta_0 t}. \]
\end{leftbar}

\noindent \textit{Proof.} Changing variables, we find
\[ \int_{\R} \left ( (\omega_k^0 - \omega_k ) \phi_{\omega_k}^2 + \frac{1}{2} \phi_{\omega_k}^4 - \frac{1}{2} \phi_{\omega_k^0}^4 \right ) = \left ( \omega_k^0 \sqrt{\omega_k} - \omega_k \sqrt{\omega_k} \right ) \int_{\R} Q^2 + \frac{\omega_k \sqrt{\omega_k} - \omega_k^0 \sqrt{\omega_k^0}}{2} \int_{\R} Q^4. \]
We know $Q$ explicitly: $Q(y) = \frac{\sqrt{2}}{\text{cosh} (y)}$. It follows that $\int_{\R}Q^2 = 4$ and $\int_{\R} Q^4 = \frac{16}{3}$. This leads to
\begin{align*}
\int_{\R} \left ( (\omega_k^0 - \omega_k ) \phi_{\omega_k}^2 + \frac{1}{2} \phi_{\omega_k}^4 - \frac{1}{2} \phi_{\omega_k^0}^4 \right ) =& \, \, 4 \omega_k^0 \sqrt{\omega_k} - \frac{4}{3} \omega_k \sqrt{\omega_k} - \frac{8}{3} \omega_k^0 \sqrt{\omega_k^0} \\
=& \, - \frac{4}{3} \left ( \sqrt{\omega_k} - \sqrt{\omega_k^0} \right )^2 \left ( \sqrt{\omega_k} + 2 \sqrt{\omega_k^0} \right ).
\end{align*}
Since $\left | \sqrt{\omega_k} (t) - \sqrt{\omega_k^0} \right | \leqslant C | \omega_k (t) - \omega_k^0 |$, we deduce the desired result:
\[ | G_0 (t) - G_0 (T_p) | \leqslant C \sum\limits_{k=1}^K | \omega_k (t) - \omega_k^0 |^2 + C e^{-2 \theta_0 t}. \]
\hfill \qedsymbol

\noindent \\ Using Proposition 1, we see that $\sum\limits_{k=1}^K | \omega_k (t) - \omega_k^0 |^2 \leqslant C A_0^2 e^{-2 \theta_0 t}$ but this is not good enough. We need the following estimate, which is obtained by pseudo-conservation of the local mass introduced in Lemma 1. It is crucial that we are in the subcritical case: since $\int_{\R} \phi_\omega^2 = \sqrt{\omega} \int_{\R} Q^2$, the integral $\int_{\R} \phi_\omega^2$ gives access to $\omega$. 

\begin{leftbar}
\noindent \textbf{Lemma 3.} For all $t \in [t^* \, , T_p]$,
\[ | \omega_k (t) - \omega_k^0 | \leqslant C ||\varepsilon_u(t)||_{L^2}^2 + C \left ( \frac{A_0^2}{L} + 1 \right ) e^{-2 \theta_0 t}. \]
\end{leftbar}

\noindent \textit{Proof.} We recall the definition of $\mathcal{M}_k (t) = \int_{\R} |u(t)|^2 \chi_k (t) \, \text{d}x$. Expanding $u = \varepsilon_u + S^u$ we get
\[  \mathcal{M}_k (t) = \int_{\R} |\varepsilon_u(t)|^2 \chi_k (t) + \int_{\R} |S^u (t) |^2 \chi_k (t) + 2 \Ree \int_{\R} \varepsilon_u(t) \overline{S^u} (t) \chi_k (t). \]
Firstly, writing $|S^u|^2 = \sum\limits_{j=1}^K | S_j^u |^2 + \text{Re} \sum\limits_{1 \leqslant j \neq \ell \leqslant K} S_j^u \overline{S_\ell^u}$, we have
\[ \int_{\R} |S^u|^2 \chi_k = \int_{\R} |S_k^u|^2 + \underbrace{\int_{\R} |S_k^u|^2 (1- \chi_k ) + \sum\limits_{\substack{j = 1 \\ j \neq k}}^K \int_{\R} |S_j^u|^2 \chi_k + \text{Re}  \sum\limits_{1 \leqslant j \neq \ell \leqslant K} \int_{\R} S_j^u \overline{S_\ell^u} \chi_k}_{= \, \mathcal{O} (e^{-2 \theta_0 t})}. \]
Secondly,
\[ \text{Re} \int_{\R} \varepsilon_u \overline{S^u} \chi_k = \underbrace{\text{Re} \int_{\R} \varepsilon_u \overline{S_k^u}}_{= \, 0} + \underbrace{\text{Re} \int_{\R} \varepsilon_u \overline{S_k^u} (1 - \chi_k) + \sum\limits_{\substack{j=1 \\ j \neq k}}^K \text{Re} \int_{\R} \varepsilon_u \overline{S_j^u} \chi_k }_{= \, \mathcal{O} (e^{-2 \theta_0 t})}. \]
This first integral is indeed $0$ since this is the orthogonality choice made. Gathering these results, we have
\[ \mathcal{M}_k (t) = \int_{\R} |\varepsilon_u|^2 \chi_k + \int_{\R} |S_k^u|^2 + \mathcal{O} (e^{-2 \theta_0 t}) = \int_{\R} |\varepsilon_u|^2 \chi_k + 4(1- c_k^2 ) \sqrt{\omega_k (t)} + \mathcal{O} (e^{-2 \theta_0 t}) \]
since $\int_{\R} |S_k^u|^2 = (1-c_k^2) \int_{\R} \phi_{\omega_k (t)}^2 = (1-c_k^2) \sqrt{\omega_k (t)} \int_{\R} Q^2 = 4(1- c_k^2 ) \sqrt{\omega_k (t)}$. Applying for $t=T_p$ we find
\[ \mathcal{M}_k (T_p) = 4 (1- c_k^2 ) \sqrt{\omega_k^0} + \mathcal{O} (e^{-2 \theta_0 t}) \]
since $\varepsilon_u(T_p) = 0$ and $\omega_k (T_p) = \omega_k^0$. It follows that
\begin{align*}
4 (1- c_k^2) \left | \sqrt{\omega_k (t)} - \sqrt{\omega_k^0} \right  | =& \, \left | \mathcal{M}_k (t) - \mathcal{M}_k (T_p) - \int_{\R} |\varepsilon_u(t)|^2 \chi_k (t) \right | + \mathcal{O} \left ( e^{-2 \theta_0 t} \right ) \\
\leqslant & \, \left | \mathcal{M}_k (t) - \mathcal{M}_k (T_p) \right | + ||\varepsilon_u(t)||_{L^2}^2 + C e^{-2 \theta_0 t}  \\
\leqslant & \, \, \frac{C A_0^2}{L} e^{-2 \theta_0 t} + ||\varepsilon_u(t)||_{L^2}^2 + C e^{-2 \theta_0 t} 
\end{align*}
by Lemma 1. We deduce the desired result: $| \omega_k (t) - \omega_k^0 | \leqslant C ||\varepsilon_u(t)||_{L^2}^2 + C \left ( \frac{A_0^2}{L} + 1 \right ) e^{-2 \theta_0 t}$. \hfill \qedsymbol

\noindent \\ Gathering the two previous lemmas, we obtain the following control of $G_0$.

\begin{leftbar}
\noindent \textbf{Lemma 4.} For $T_0>0$ chosen large enough, for all $t \in [t^* \, , T_p ]$,
\[ |G_0 (t) - G_0(T_p)| \leqslant C e^{-2 \theta_0 t}. \]
\end{leftbar}

\noindent \textit{Proof.} Recall that $||(\varepsilon_u \, , \varepsilon_n \, , \varepsilon_v)||_{\mathbf{H}} \leqslant C A_0 e^{- \theta_0 t}$. Using Lemma 3 we get $| \omega_k (t) - \omega_k^0 | \leqslant C A_0^2 e^{-2 \theta_0 t} \leqslant C e^{- \theta_0 t}$ for $T_0>0$ large enough (depending on $A_0$). Then using Lemma 2 we obtain the desired result. \hfill \qedsymbol

\subsection{Control of first-order term $G_1$}
\noindent Let us now control $G_1$. The Weinstein functional is constructed so that this first-order term vanishes (up to exponential error terms).

\begin{leftbar}
\noindent \textbf{Lemma 5.} For all $t \in [t^* \, , T_p]$,
\[ |G_1(t)| \leqslant C e^{-2 \theta_0 t}. \]
\end{leftbar}

\noindent \textit{Proof.} First, as we did for $G_0$, 
\begin{align*} G_1 =& \, \sum\limits_{k=1}^K \int_{\R} \left ( 2 \, \text{Re} ( \partial_x S_k^u   \partial_x \overline{\varepsilon_u} ) + 2 S_k^n \text{Re} (S_k^u \overline{\varepsilon_u}) + \varepsilon_n |S_k^u|^2 + \varepsilon_n S_k^n +  \varepsilon_vS_k^v + 2 \nu_k^0 \text{Re} (S_k^u \overline{\varepsilon_u}) \chi_k \right. \\
& \midspace \midspace \left. -c_k   \text{Im} ( \overline{S_k^u} ( \partial_x \varepsilon_u ) \chi_k + \overline{\varepsilon_u} (\partial_x S_k^u ) \chi_k ) - c_k   S_k^n  \varepsilon_v \chi_k - c_k   \varepsilon_n S_k^v \chi_k \right ) + \mathcal{O} (e^{-2 \theta_0 t}).
\end{align*}
Recall that $\nu_k^0 = \omega_k^0 + \frac{c_k^2}{4}$. Using 
\begin{align*}
& \int_{\R} \partial_x S_k^u   \partial_x \overline{\varepsilon_u} = - \int_{\R} ( \partial_x^2 S_k^u ) \overline{\varepsilon_u} \\
\text{and} \, \, \, & \partial_x^2 S_k^u = \sqrt{1-c_k^2} \, e^{i \Gamma_k (t  , x)} \left [ \omega_k (t) \phi_{\omega_k(t)} - \phi_{\omega_k(t)}^3 + ic_k   \phi_{\omega_k (t)} ' - \frac{c_k^2}{4} \phi_{\omega_k(t)} \right ](x_k^t),
\end{align*}
we compute
\begin{align*} 
G_1 =& \, \sum\limits_{k=1}^K \int_{\R} \left [ -2 \, \text{Re} \left ( \left ( \sqrt{1-c_k^2} \, e^{-i \Gamma_k (t  , x)} \left ( \boxed{\omega_k} \, \phi_{\omega_k} - \phi_{\omega_k}^3 - i c_k   \phi_{\omega_k} ' - \frac{c_k^2}{4} \phi_{\omega_k} \right ) (x_k^t) \right. \right. \right. \\
& \midspace \midspace \midspace \midspace \midspace \midspace \left. \left. \left. + \phi_{\omega_k}^2 (x_k^t)   \sqrt{1-c_k^2} \, e^{-i \Gamma (t  , x)} \phi_{\omega_k} (x_k^t) \right ) \varepsilon_u \right ) \right.  \\
& \midspace \midspace \midspace + \left ( - \phi_{\omega_k}^2 (x_k^t) + (1-c_k^2) \phi_{\omega_k}^2 (x_k^t) \right ) \varepsilon_n - c_k \phi_{\omega_k}^2 (x_k^t)    \varepsilon_v \\
& \midspace \midspace \midspace + 2 \left ( \boxed{\omega_k^0} + \frac{c_k^2}{4} \right ) \text{Re} \left ( \sqrt{1-c_k^2} \, e^{-i \Gamma_k (t  , x)} \phi_{\omega_k} (x_k^t) \varepsilon_u \right ) \\
& \midspace \midspace \midspace + 2 c_k   \text{Im} \left ( \sqrt{1-c_k^2} \, e^{-i \Gamma_k(t \, , x)} \left ( \phi_{\omega_k} ' (x_k^t) - \frac{i c_k}{2} \phi_{\omega_k} (x_k^t) \right ) \varepsilon_u \right ) \\
& \midspace \midspace \midspace \left. + c_k \phi_{\omega_k}^2 (x_k^t)    \varepsilon_v + c_k^2 \phi_{\omega_k}^2 (x_k^t) \varepsilon_n \right ] \\
&  \midspace + \mathcal{O} (e^{-2 \theta_0 t}). 
\end{align*}
Almost every term of this expression vanishes, except for the framed terms. We eventually get
\[ G_1 = 2 \sum\limits_{k=1}^K ( \omega_k^0 - \omega_k (t)) \sqrt{1- c_k^2} \, \text{Re} \int_{\R} \phi_{\omega_k} (x_k^t) e^{-i \Gamma_k} \varepsilon_u \, \text{d}x + \mathcal{O} (e^{-2 \theta_0 t}). \]
Thanks to the orthogonality $\text{Re} \int_{\R} \phi_{\omega_k} (x_k^t) e^{-i \Gamma_k} \varepsilon_u \, \text{d}x = 0$, we obtain the desired result: $|G_1 (t)| \leqslant C e^{- 2 \theta_0 t}$. \hfill \qedsymbol

\subsection{Control of second-order term $G_2$}
\noindent The crucial point is the coercivity of $G_2$. Let us write $G_2 = G_{2,1} + G_{2,2}$ with
\begin{align*} 
G_{2,1} =& \int_{\R} \left ( | \partial_x \varepsilon_u|^2 + 2 \varepsilon_n \, \text{Re} (S^u \overline{\varepsilon_u}) + S^n |\varepsilon_u|^2 + \frac{\varepsilon_n^2}{2} + \frac{( \varepsilon_v )^2}{2} \right ) + \sum\limits_{k=1}^K \left ( \omega_k (t) + \frac{c_k^2}{4} \right ) \int_{\R} |\varepsilon_u|^2 \chi_k \\
& \midspace - \sum\limits_{k=1}^K c_k   \int_{\R} \left ( \text{Im} \left ( \overline{\varepsilon_u} ( \partial_x \varepsilon_u ) \chi_k \right ) + \varepsilon_n \varepsilon_v \chi_k \right )
\end{align*}
and $G_{2,2} = \sum\limits_{k=1}^K \left ( \omega_k^0 - \omega_k (t) \right ) \int_{\R}  |\varepsilon_u|^2 \chi_k$. Through a progression of three coercivity lemmas, let us control $G_{2,1}$. To do so, we adapt Lemmas 2.6 and 4.1 from \cite{Ma4}, where analogous estimates are established for NLS multi-solitons. First, here is the basic, one-soliton non-localized version.

\begin{leftbar}
\noindent \textbf{Lemma 6.} Let $c \in (-1 \, , 1)$, $\omega \in [\omega_- \, , \omega^+]$, $\sigma \in \R$ and $\gamma \in \R$. Define $\Gamma (t \, , x) = \frac{c   x}{2} - \frac{c^2t}{4} + \omega t + \gamma$, $\nu = \omega + \frac{c^2}{4}$ and
\begin{align*} 
H_2 (\eta_u \, , \eta_n \, ,  \eta_v) =& \int_{\R} \left [ | \partial_x \eta_u |^2 + 2 \sqrt{1-c^2} \, \eta_n \phi_\omega (x-ct- \sigma ) \text{Re} \left ( e^{i \Gamma (t  , x)} \overline{\eta_u} \right )  \right. \\
& \midspace \midspace \midspace \left. \, - \phi_\omega^2 (x-ct- \sigma ) |\eta_u|^2 + \frac{\eta_n^2 + \eta_v^2}{2} + \nu |\eta_u|^2 - c   \left ( \eta_n \eta_v + \text{Im} ( \overline{\eta_u} \partial_x \eta_u ) \right ) \right ]
\end{align*}
There exists $C>0$ such that, for all $(\eta_u \, , \eta_n \, ,  \eta_v) \in \mathbf{H}$ satisfying the following orthogonality relations:
\begin{align*}
& \text{Re} \int_{\R} \phi_\omega (x-ct- \sigma ) e^{-i \Gamma (t  , x)} \eta_u (x) \, \text{d}x = 0, \\
& \text{Re} \int_{\R} (x-ct- \sigma ) \phi_\omega (x-ct- \sigma ) e^{-i \Gamma (t  , x)} \eta_u (x) \, \text{d}x = 0 \\
\text{and} \, \, \, & \text{Im} \int_{\R} \Lambda_\omega (x-ct- \sigma ) e^{-i \Gamma (t  , x)} \eta_u (x) \, \text{d}x = 0;
\end{align*}
then
\[ ||(\eta_u \, , \eta_n \, ,  \eta_v)||_{\mathbf{H}}^2 \leqslant C H_2 (\eta_u \, , \eta_n \, ,  \eta_v). \]
The constant $C$ depends on $\omega_-$ and $\omega^+$ but does not depend on $\omega$ itself. It also does not depend on $c$, $\sigma$ or~$\gamma$.
\end{leftbar}

\noindent \textit{Proof.} Recalling that $\phi_\omega (x-ct- \sigma ) = \sqrt{\omega} \, Q ( \sqrt{\omega} (x-ct-\sigma ))$, define $\tilde{\eta}_u$, $\tilde{\eta}_n$ and $\tilde{\eta}_v$ such that $\eta_u (t \, , x) = \sqrt{\omega} \, \tilde{\eta}_u (t \, , \sqrt{\omega} (x-ct- \sigma )) e^{i \left ( \frac{c   x}{2} - \frac{c^2t}{4} + \omega t + \gamma \right )}$, $\eta_v(t \, , x) = \omega \tilde{\eta}_v (t \, , \sqrt{\omega} (x-ct - \sigma ))$ and $\eta_n = \omega \tilde{\eta}_n (t \, , \sqrt{\omega} (x-ct- \sigma ))$. With this scaling we have
\[ \omega^{- \frac{3}{2}} H_2 = \int_{\R} \left ( | \partial_x \tilde{\eta}_u |^2 + (1-Q^2) | \tilde{\eta}_u |^2 + \frac{\tilde{\eta}_n^2 + \tilde{\eta}_v^2}{2} + 2 \sqrt{1-c^2} \, Q \tilde{\eta}_n \, \text{Re} ( \tilde{\eta}_u ) - c   \tilde{\eta}_n \tilde{\eta}_v \right ). \]
Take $\delta > 0$ to be adjusted later on. Define $\alpha = 1 + \frac{\delta (1-c^2)}{2(4 + \delta ) c^2} > 1$ and $\beta = \frac{4}{4 + \delta} \in (0 \, , 1)$. According to Young's inequality, we have
\begin{align*}
& \left | c   \tilde{\eta}_n \tilde{\eta}_v \right | \leqslant \frac{\alpha c^2}{2} \tilde{\eta}_n^2 + \frac{1}{2 \alpha}  \tilde{\eta}_v^2 \\
\text{and} \, \, \, & \left | 2 \sqrt{1-c^2} \, Q \tilde{\eta}_n \, \text{Re} ( \tilde{\eta}_u ) \right | \leqslant \frac{\beta (1-c^2)}{2} \tilde{\eta}_n^2 + \frac{2Q^2}{\beta} \text{Re} ( \tilde{\eta}_u )^2.
\end{align*}
This leads to
\[ \omega^{- \frac{3}{2}} H_2 \geqslant \int_{\R} \left ( | \partial_x \tilde{\eta}_u |^2 + (1-Q^2) | \tilde{\eta}_u |^2 + \mu_n ( \delta ) \tilde{\eta}_n^2 + \mu_v ( \delta ) \tilde{\eta}_v^2 - \left ( 2 + \frac{\delta}{2} \right ) Q^2 \text{Re} ( \tilde{\eta}_u )^2 \right ) \]
where $\mu_n ( \delta ) := \frac{1- \beta + \beta c^2 - \alpha c^2}{2} = \frac{\delta (1-c^2)}{2 ( 4 + \delta )} >0$ and $\mu_v ( \delta ) := \frac{1 - \alpha^{-1}}{2} >0$. Since $\left | \frac{\delta}{2} Q^2 \text{Re} ( \tilde{\eta}_u )^2 \right | \leqslant \delta | \tilde{\eta}_u |^2$, we have
\begin{align*}
C H_2 \geqslant & \int_{\R} \left ( | \partial_x \tilde{\eta}_u |^2 + (1-Q^2) | \tilde{\eta}_u |^2 + \mu_n ( \delta ) \tilde{\eta}_n^2 + \mu_v ( \delta ) \tilde{\eta}_v^2 -2 Q^2 \text{Re} ( \tilde{\eta}_u )^2 - \delta | \tilde{\eta}_u |^2 \right ) \\
\geqslant & \int_{\R} \left ( | \partial_x \text{Re} ( \tilde{\eta}_u ) |^2 + (1-3Q^2) \text{Re} ( \tilde{\eta}_u )^2 \right ) + \int_{\R} \left ( | \partial_x \text{Im} ( \tilde{\eta}_u ) |^2 + (1-Q^2) \text{Im} ( \tilde{\eta}_u )^2 \right ) \\
& \midspace \midspace + \mu_n ( \delta ) || \tilde{\eta}_n ||_{L^2}^2 + \mu_v ( \delta ) || \tilde{\eta}_v ||_{L^2}^2 - \delta || \tilde{\eta}_u ||_{L^2}^2 \\
\geqslant & \, \, \langle L_+ \text{Re} ( \tilde{\eta}_u ) \, , \text{Re} ( \tilde{\eta}_u ) \rangle + \langle L_- \text{Im} ( \tilde{\eta}_u ) \, , \text{Im} ( \tilde{\eta}_u ) \rangle + \mu_n ( \delta ) || \tilde{\eta}_n ||_{L^2}^2 + \mu_v ( \delta ) || \tilde{\eta}_v ||_{L^2}^2 - \delta || \tilde{\eta}_u ||_{H^1}^2 . 
\end{align*}
The following orthogonality relations hold: $\langle \text{Re} ( \tilde{\eta}_u ) \, , Q \rangle = \langle \text{Re} ( \tilde{\eta}_u ) \, , yQ \rangle = \langle \text{Im} ( \tilde{\eta}_u ) \, , \Lambda Q \rangle = 0$. It follows from \eqref{spectral} that 
\[ \langle L_+ \text{Re} ( \tilde{\eta}_u ) \, , \text{Re} ( \tilde{\eta}_u ) \rangle + \langle L_- \text{Im} ( \tilde{\eta}_u ) \, , \text{Im} ( \tilde{\eta}_u ) \rangle \geqslant \lambda_0 || \tilde{\eta}_u ||_{H^1}^2 \]
where $\lambda_0 > 0$ is a universal constant. Thus,
\[ C H_2 \geqslant ( \lambda_0 - \delta ) || \tilde{\eta}_u ||_{H^1}^2 + \mu_n ( \delta ) || \tilde{\eta}_n ||_{L^2}^2 + \mu_v ( \delta ) || \tilde{\eta}_v ||_{L^2}^2 \geqslant C' || ( \tilde{\eta}_u \, , \tilde{\eta}_v \, , \tilde{\eta}_n ) ||_{\mathbf{H}}^2 \]
taking $\delta = \frac{\lambda_0}{2}$ and $C' = \min \left ( \frac{\lambda_0}{2} \, , \mu_n \left ( \frac{\lambda_0}{2} \right ) \, , \mu_v \left ( \frac{\lambda_0}{2} \right ) \right )$. \\
\\ Going back to $\eta_u$, we write that $\partial_x \eta_u (x \, , t) = \sqrt{\omega} \, e^{i \left ( \frac{c   x}{2} - \frac{c^2t}{4} + \omega t + \gamma \right )} \left ( \sqrt{\omega} \, \partial_x \tilde{\eta}_u + \frac{ic}{2} \right ) (t \, , \sqrt{\omega} (x-ct- \sigma ))$. Thus
\begin{align*}
\int_{\R} | \partial_x \eta_u |^2 =& \, \, \sqrt{\omega} \left ( \omega \int_{\R} | \partial_x \tilde{\eta}_u |^2 + \frac{c^2}{4} \int_{\R} | \tilde{\eta}_u|^2 + c   \sqrt{\omega} \, \text{Im} \int_{\R} \overline{\tilde{\eta}_u} \partial_x \tilde{\eta}_u \right ) \\
\leqslant & \, \, \omega^{\frac{3}{2}} \frac{3}{2} \int_{\R} | \partial_x \tilde{\eta}_u |^2 + \sqrt{\omega} \, \frac{3c^2}{4} \int_{\R} | \tilde{\eta}_u |^2. 
\end{align*}
We also have $\int_{\R} |\eta_u|^2 = \sqrt{\omega} \int_{\R} | \tilde{\eta}_u |^2$, which results in $||\eta_u||_{H^1}^2 \leqslant C || \tilde{\eta}_u ||_{H^1}^2$. In addition, $||\eta_n||_{L^2}^2 = \omega^{\frac{3}{2}} || \tilde{\eta}_n ||_{L^2}^2$ and $||\eta_v||_{L^2}^2 = \omega^{\frac{3}{2}} || \tilde{\eta}_v ||_{L^2}^2$. Hence, $|| (\eta_u  \, , \eta_n \, , \eta_v) ||_{\mathbf{H}}^2 \leqslant C || ( \tilde{\eta}_u  \, , \tilde{\eta}_n \, , \tilde{\eta}_v) ||_{\mathbf{H}}^2$. It follows that
\[ ||(\eta_u \, , \eta_n \, , \eta_v)||_{\mathbf{H}}^2 \leqslant C H_2 (\eta_u \, , \eta_n \, , \eta_v) \]
which is the desired result.
\hfill \qedsymbol

\noindent \\ Now take $\Phi \in \mathscr{C}^2 ( \R \, , \R )$ even, decreasing on $[0 \, , + \infty )$, such that $\Phi \equiv 1$ on $[0 \, , 1 ]$, $\Phi (x) = e^{-x}$ on $[2 \, , + \infty )$ and $e^{-|x|} \leqslant \Phi (x) \leqslant 3 e^{-|x|}$ on $\R$. Take $B>0$ and $\Phi_B (x) = \Phi \left ( \frac{x}{B} \right )$. We give the following adaptation of Lemma 6, where we still consider only one soliton, but we localize it using the function $\Phi_B$.

\begin{leftbar}
\noindent \textbf{Lemma 7.} Let $c \in (-1 \, , 1)$, $\omega \in [\omega_- \, , \omega^+]$, $\sigma \in \R$ and $\gamma \in \R$. Define $\Gamma (t \, , x) = \frac{c   x}{2} - \frac{c^2t}{4} + \omega t + \gamma$, $\nu = \omega + \frac{c^2}{4}$ and
\begin{align*}
H_{2B} (\eta_u \, , \eta_n \, , \eta_v) =& \int_{\R} \Phi_B (x-ct - \sigma ) \left ( | \partial_x \eta_u |^2 + \nu |\eta_u|^2 - c   \left ( \eta_n \eta_v + \text{Im} ( \overline{\eta_u} \partial_x \eta_u ) \right ) + \frac{\eta_n^2 + \eta_v^2}{2} \right ) \\
& \midspace + \int_{\R} \left ( 2 \sqrt{1-c^2} \, \phi_\omega (x-ct - \sigma ) \eta_n \, \text{Re} ( e^{i \Gamma (t  , x)} \overline{\eta_u} ) - \phi_\omega^2 (x-ct - \sigma ) |\eta_u|^2 \right ).
\end{align*}
There exists $C >0$ such that, for $B>0$ large enough (depending on $\omega_-$ and $\omega^+$ but not on $\omega$ itself), for all $(\eta_u \, , \eta_n \, , \eta_v) \in \mathbf{H}$ satisfying the following orthogonality relations:
\begin{align*}
& \text{Re} \int_{\R} \phi_\omega (x-ct- \sigma ) e^{-i \Gamma (t  , x)} \eta_u (x) \, \text{d}x = 0, \\
& \text{Re} \int_{\R} (x-ct- \sigma ) \phi_\omega (x-ct- \sigma ) e^{-i \Gamma (t  , x)} \eta_u (x) \, \text{d}x = 0 \\
\text{and} \, \, \, & \text{Im} \int_{\R} \Lambda_\omega (x-ct- \sigma ) e^{-i \Gamma (t  , x)} \eta_u (x) \, \text{d}x = 0;
\end{align*}
then
\[ H_{2B} (\eta_u \, , \eta_n \, , \eta_v) \geqslant C \int_{\R} \left ( |\eta_u|^2 + | \partial_x \eta_u |^2 + \eta_n^2 + \eta_v^2 \right ) \Phi_B (x-ct- \sigma ).  \]
The constant $C$ depends on $\omega_-$ and $\omega^+$ but does not depend on $\omega$ itself. It also does not depend on $c$, $\sigma$ or~$\gamma$.
\end{leftbar}

\noindent \textit{Proof.} Denote $x^t = x-ct- \sigma$ for simplicity. Let $z_u (t \, , x) = \eta_u (t \, , x) \sqrt{\Phi_B (x^t)}$, $z_v (t \, , x) = \eta_v (t \, , x) \sqrt{\Phi_B (x^t)}$ and $z_n(t \, , x) = \eta_n (t \, , x) \sqrt{\Phi_B (x^t)}$. We have
\[ \int_{\R} | \partial_x z_u |^2 - \frac{C}{B} \int_{\R} (| \partial_x z_u |^2 + |z_u|^2 ) \leqslant \int_{\R} | \partial_x \eta_u |^2 \Phi_B (x^t) \leqslant \int_{\R} | \partial_x z_u |^2 + \frac{C}{B} \int_{\R} ( | \partial_x z_u |^2 + |z_u|^2 ). \]
Moreover, 
\begin{align*} 
& \text{Im} \left ( \overline{\eta_u} \partial_x \eta_u \right ) \Phi_B (x^t) = \text{Im} \left ( \overline{z_u} \partial_x z_u - \frac{\Phi_B ' (x^t)}{2 \Phi_B (x^t)} |z_u|^2 \right ) = \text{Im} ( \overline{z_u} \partial_x z_u ), & & \int_{\R} |\eta_u|^2 \Phi_B (x^t) \, \text{d}x = \int_{\R} | z_u |^2, \\
& \int_{\R} \eta_n^2 \Phi_B (x^t) \, \text{d}x = \int_{\R} z_n^2, & & \int_{\R} \eta_v^2 \Phi_B (x^t) \, \text{d}x = \int_{\R} z_v^2 \\
\text{and} \, \, & \int_{\R} \eta_n \eta_v \Phi_B (x^t) \, \text{d}x = \int_{\R} z_n z_v. 
\end{align*}
For $|x| \leqslant B$, $\left | \frac{1}{\Phi_B} -1 \right | \phi_\omega (x)  =0$. For $|x| > B$, 
\begin{equation} 
\left | \frac{1}{\Phi_B} -1 \right | \phi_\omega (x) = \frac{1 - \Phi_B}{\Phi_B} \phi_\omega \leqslant \left ( 1 - e^{- |x|/B} \right ) e^{|x| / B} C e^{- \sqrt{\omega} |x|} \leqslant C e^{- \left ( \sqrt{\omega} - \frac{1}{B} \right ) |x|} \leqslant C e^{- \left ( \sqrt{\omega} B -1 \right )} \leqslant \frac{C}{B}
\label{phiB}
\end{equation}
for $B>0$ large enough (depending on $\omega_-$). It follows that
\begin{align*}
& \int_{\R} \left ( 2 \sqrt{1-c^2} \, \phi_\omega (x^t) \eta_n \, \text{Re} \left ( e^{-i \Gamma (t  , x)} \eta_u \right ) - \phi_\omega^2 (x^t) |\eta_u|^2 \right ) \\
=& \int_{\R} \left ( 2 \sqrt{1-c^2} \, z_n \, \text{Re} \left ( e^{-i \Gamma (t  , x)} z_u \right ) - \phi_\omega (x^t) |z_u|^2 \right ) \left ( \frac{\phi_\omega}{\Phi_B} \right ) (x^t) \\
\leqslant & \int_{\R} \left ( 2 \sqrt{1-c^2} \, \phi_\omega (x^t) z_n \, \text{Re} \left ( e^{-i \Gamma (t  , x)} z_u \right ) - \phi_\omega^2 (x^t) |z_u|^2 \right ) + \frac{C}{B} \int_{\R} (|z_nz_u| + |z_u|^2) \\
\leqslant & \int_{\R} \left ( 2 \sqrt{1-c^2} \, \phi_\omega (x^t) z_n \, \text{Re} \left ( e^{-i \Gamma (t  , x)} z_u \right ) - \phi_\omega^2 (x^t) |z_u|^2 \right ) + \frac{C}{B} \int_{\R} \left ( z_n^2 + |z_u|^2 \right ). 
\end{align*}
Hence
\[ H_{2B} (\eta_u \, , \eta_n \, , \eta_v) \geqslant H_2 (z_u \, , z_n \, , z_v ) - \frac{C}{B} \int_{\R} ( z_n^2 + |z_u|^2 + | \partial_x z_u |^2 ). \]
Now write 
\[ z_u(t \, , x) = z_u^\perp (t \, , x) + a \phi_\omega (x^t) e^{i \Gamma (t  , x)} + b   x^t \phi_\omega (x^t) e^{i \Gamma (t  , x)} + i \beta \Lambda_\omega (x^t) e^{i \Gamma} \]
with 
\begin{align*}
& a = \frac{1}{\int_{\R} \phi_\omega^2} \Ree \int_{\R}  \phi_\omega (x^t ) e^{-i \Gamma (t  , x)} z_u (x) \, \text{d}x, \\
& b = \frac{1}{\int_{\R} y^2 \phi_\omega^2} \Ree \int_{\R} x^t \phi_\omega (x^t) e^{-i \Gamma (t,x)} z_u (x) \, \text{d}x \\
\text{and} \, \, \, & \beta = \frac{1}{\int_{\R} \Lambda_\omega^2} \Imm \int_{\R} \Lambda_\omega (x^t) e^{-i \Gamma (t,x)} z_u (x) \, \text{d}x. 
\end{align*}
Let us show that $|a | \leqslant \frac{C}{B} || z_u ||_{H^1}$. Thanks to the orthogonality $\text{Re} \int_{\R} \phi_\omega (x^t) e^{-i \Gamma (t  , x)} \eta_u(x) = 0$, we have 
\begin{align*}
\left | \text{Re} \int_{\R} \phi_\omega (x^t ) e^{-i \Gamma (t  , x)} z_u (x) \right | \text{d}x =& \, \left | \text{Re} \int_{\R} \left ( \phi_\omega \sqrt{\Phi_B} \right ) (x^t) e^{-i \Gamma (t  , x)} \eta_u(x) \right | \text{d}x \\
= & \, \int_{\R}  \left ( \phi_\omega \left | \sqrt{\Phi_B} -1 \right | \right ) (x^t ) e^{-i \Gamma (t  , x)} |\eta_u|(x) \, \text{d}x \\
\leqslant & \int_{\R} \left ( \frac{\phi_\omega \left ( 1 - \sqrt{\Phi_B} \right )}{\sqrt{\Phi_B}} \right ) (x^t) |z_u|(x) \, \text{d}x  \\
\leqslant & \, \sqrt{\int_{\R} \left ( \frac{\phi_\omega \left ( 1 - \sqrt{\Phi_B} \right )}{\sqrt{\Phi_B}} \right )^2} \, ||z_u||_{L^2} .
\end{align*}
We have (see \eqref{phiB})
\[ \int_{\R} \frac{\phi_\omega^2 \left ( 1 - \sqrt{\Phi_B} \right )^2}{\Phi_B} \leqslant \int_{\R} \phi_\omega \frac{\phi_\omega ( 1 - \Phi_B )}{\Phi_B} \leqslant \frac{C}{B} \int_{\R}  \phi_\omega \leqslant \frac{C}{B} \]
for $B>0$ large enough (depending on $\omega_-$ and $\omega^+$). It follows that
\[ | a | = \frac{1}{C \sqrt{\omega}} \left | \text{Re} \int_{\R} \phi_\omega (x^t) e^{-i \Gamma (t  , x)} z_u (x) \, \text{d}x \right | \leqslant \frac{C}{B} ||z_u||_{L^2} \leqslant \frac{C}{B} || z_u ||_{H^1}. \]
Using the other orthogonalities $\text{Re} \int_{\R} x^t \phi_\omega (x^t) e^{-i \Gamma (t  , x)} \eta_u = 0$ and $\text{Im} \int_{\R} \Lambda_\omega (x^t) e^{-i \Gamma (t  , x)} \eta_u = 0$, we similarly prove that $| b | \leqslant \frac{C}{B} || z_u ||_{H^1}$ and $| \beta | \leqslant \frac{C}{B} || z_u ||_{H^1}$. Expanding $| \partial_x z_u (t \, , x) |^2$, it follows that
$\left | \int_{\R} | \partial_x z_u |^2 - \int_{\R} | \partial_x z_u^\perp |^2 \right | \leqslant \frac{C}{B} ||z_u||_{H^1}^2$. This leads to
\[ H_{2B} (\eta_u \, , \eta_n \, , \eta_v) \geqslant H_2 (z_u \, , z_n \, , z_v) - \frac{C}{B} ||(z_u \, , z_n \, , z_v )||_{\mathbf{H}}^2 \geqslant H_2 ( z_u^\perp \, , z_v \, , z_n ) - \frac{C}{B} ||(z_u \, , z_n \, , z_v)||_{\mathbf{H}}^2. \]
Now, $z_u^\perp$ satisfies the orthogonality relations 
\begin{align*} 
& \text{Re} \int_{\R} \phi_\omega (x^t) e^{-i \Gamma (t  , x)} z_u^\perp (x) \, \text{d}x = 0, 
\\ & \text{Re} \int_{\R} x^t \phi_\omega (x^t) e^{-i \Gamma (t  , x)} z_u^\perp (x) \, \text{d}x = 0 \\
\text{and} \, \, \, & \text{Im} \int_{\R} \Lambda_\omega (x^t) e^{-i \Gamma (t  , x)} z_u^\perp (x) \, \text{d}x = 0.
\end{align*}
Thus 
\[ H_2 ( z_u^\perp \, , z_n \, , z_v ) \geqslant C || (z_u^\perp \, , z_n \, , z_v) ||_{\mathbf{H}}^2  . \]
Recall that $|| z_u^\perp ||_{H^1}^2 \geqslant ||z_u||_{H^1}^2 - \frac{C}{B} ||z_u||_{H^1}^2$ thus $|| (z_u^\perp \, , z_n \, , z_v) ||_{\mathbf{H}}^2 \geqslant \left (1- \frac{C}{B} \right )  ||(z_u \, , z_n \, , z_v)||_{\mathbf{H}}^2$. Hence
\[ H_2 (z_u^\perp \, , z_n \, , z_v) \geqslant C \left ( 1 - \frac{C}{B} \right ) ||(z_u \, , z_n \, , z_v)||_{\mathbf{H}}^2  \geqslant C ||(z_u \, , z_n \, , z_v)||_{\mathbf{H}}^2  \]
for $B>0$ large enough (depending on $\omega_-$ and $\omega^+$). It follows that
\begin{align*}
H_{2B} (\eta_u \, , \eta_n \, , \eta_v) \geqslant & \, \, H_2 ( z_u^\perp \, , z_n \, , z_v) - \frac{C}{B} || (z_u \, , z_n \, , z_v) ||_{\mathbf{H}}^2  \\
\geqslant & \, \, C \left ( 1- \frac{C}{B} \right ) || (z_u \, , z_n \, , z_v ) ||_{\mathbf{H}}^2 \, \, \geqslant \, \, C ||(z_u \, , z_n \, , z_v)||_{\mathbf{H}}^2
\end{align*}
for $B>0$ large enough (depending on $\omega_-$ and $\omega^+$). \\
\\ Now, going back to $\eta_u$, recall that $z_u (x) = \eta_u(x) \sqrt{\Phi_B (x^t)}$ thus $\partial_x z_u (x) = \partial_x \eta_u (x) \sqrt{\Phi_B (x^t)} + \eta_u (x) \frac{\Phi_B ' (x^t )}{2 \sqrt{\Phi_B (x^t )}}$. We have
\[ ||z_u||_{H^1}^2 = \int_{\R} ( |\eta_u|^2 + | \partial_x \eta_u |^2 ) \Phi_B (x^t ) + \int_{\R} \left ( \frac{| \Phi_B ' |^2 (x^t)}{4 \Phi_B (x^t )} |\eta_u|^2 + \Phi_B ' (x^t )   \text{Re} ( \overline{\eta_u} \partial_x \eta_u ) \right ). \]
Since $| \Phi_B '| \leqslant \frac{C}{B} \Phi_B$, it follows that $\left | ||z_u||_{H^1}^2 - \int_{\R} \left ( |\eta_u|^2 + | \partial_x \eta_u |^2 \right ) \Phi_B (x^t) \right | \leqslant \frac{C}{B} \int_{\R} \left ( |\eta_u|^2 + | \partial_x \eta_u |^2 \right ) \Phi_B (x^t)$. Hence
\begin{align*}
H_{2B} (\eta_u \, , \eta_n \, , \eta_v) \geqslant & \, \, C \left ( 1 - \frac{C}{B} \right ) \int_{\R} \left ( |\eta_u|^2 +  | \partial_x \eta_u |^2 + \eta_n^2 + \eta_v^2 \right ) \Phi_B (x^t) \\
\geqslant & \, \, C \int_{\R} \left ( |\eta_u|^2 + | \partial_x \eta_u |^2 + \eta_n^2 + \eta_v^2 \right ) \Phi_B (x^t) \\ 
\end{align*}
for $B>0$ large enough (depending on $\omega_-$ and $\omega^+$). \hfill \qedsymbol

\noindent \\ Lastly, we give the version of Lemmas 6 and 7 with several solitons and the presence of the localization functions $\chi_k$. We establish the coercivity of $G_{2,1}$.

\begin{leftbar}
\noindent \textbf{Lemma 8.} There exists $C >0$ such that, for all $t \in [t^* \, , T_p ]$,
\[ G_{2,1} (t) \geqslant C || (\varepsilon_u \, , \varepsilon_n \, , \varepsilon_v)(t) ||_{\mathbf{H}}^2 . \]
\end{leftbar}

\noindent \textit{Proof.} Since $\sum\limits_{k=1}^K \chi_k = 1$, we have
\begin{align*}
G_{2,1} =& \, \sum\limits_{k=1}^K H_{2B,k} (\varepsilon_u \, , \varepsilon_n \, , \varepsilon_v)   \\ 
& \midspace + \sum\limits_{k=1}^K \int_{\R} \left ( \chi_k (t) - \Phi_B (x_k^t ) \right )  \left ( | \partial_x \varepsilon_u |^2 + \nu_k (t) |\varepsilon_u|^2 + \frac{\varepsilon_n^2 + \varepsilon_v^2}{2} - c_k   ( \varepsilon_n\varepsilon_v + \text{Im} ( \overline{\varepsilon_u} \partial_x \varepsilon_u)) \right ) 
\end{align*}
where, for any $k \in \{ 1 \, , ... \, , K \}$, $\nu_k (t) = \omega_k (t) + \frac{c_k^2}{4}$ and
\begin{align*}
H_{2B,k} (\varepsilon_u \, , \varepsilon_n \, , \varepsilon_v) =& \int_{\R} \Phi_B (x_k^t ) \left ( | \partial_x \varepsilon_u |^2 + \nu_k (t) |\varepsilon_u|^2 - c_k   \left ( \varepsilon_n\varepsilon_v + \text{Im} ( \overline{\varepsilon_u} \partial_x \varepsilon_u ) \right ) + \frac{\varepsilon_n^2 + \varepsilon_v^2}{2} \right ) \\
& \midspace + \int_{\R} \left ( 2 \sqrt{1-c_k^2} \, \phi_{\omega_k} (x_k^t ) \varepsilon_n \, \text{Re} ( e^{i \Gamma_k (t  , x)} \overline{\varepsilon_u} ) - \phi_\omega^2 (x_k^t) |\varepsilon_u|^2 \right ).
\end{align*}
The following convenient orthogonality relations hold (see Proposition 1):
\begin{align*}
& \text{Re} \int_{\R} \phi_{\omega_k (t)} (x_k^t) e^{-i \Gamma_k (t , x)} \varepsilon_u(x) \, \text{d}x = 0 , \\
& \text{Re} \int_{\R} x_k^t \phi_{\omega_k (t)} (x_k^t) e^{-i \Gamma_k (t  , x)} \varepsilon_u(t \, , x) \, \text{d}x = 0 \\
\text{and} \, \, \, & \text{Im} \int_{\R} \Lambda_{\omega_k (t)} (x_k^t) e^{-i \Gamma_k (t , x)} \varepsilon_u(t \, , x) \, \text{d}x = 0
\end{align*}
This enables to apply Lemma 7 with $B>0$ large enough (depending on $\omega_-$ and $\omega^+$ but not on $t$). Recall \eqref{omega}. The constant $C>0$ depends as well on $\omega_-$ and $\omega^+$ but not on $t$. The estimate writes as follows:
\[ H_{2B,k} (\varepsilon_u \, , \varepsilon_n \, , \varepsilon_v) \geqslant C \int_{\R} \Phi_B (x_k^t ) \left ( | \varepsilon_u |^2 + | \partial_x \varepsilon_u  |^2 + \varepsilon_n^2 + \varepsilon_v^2 \right )  \]
for all $k \in \{ 1 \, , ... \, , K \}$. \\
\\ Take $L>0$ to be fixed later. Take $k \in \{ 1 \, , ... \, , K \}$. For $t>t^*>T_0$ large enough, we have 
\[ | \chi_k (t) - \Phi_B (x- c_k t - \sigma_k (t) ) | \leqslant e^{- \frac{L}{4B}} \midspace \text{thus} \midspace \chi_k (t) - \Phi_B (x- c_k t - \sigma_k (t)) \geqslant - e^{- \frac{L}{4B}}. \]
Moreover, set $\alpha_k (t) = 2 - \frac{4 \omega_k (t)}{2 \omega_k (t) + c_k^2} = \frac{2c_k^2}{2 \omega_k + c_k^2} \in (0 \, , 2)$. If $c_k=0$, we have
\[ | \partial_x \varepsilon_u |^2 + \nu_k (t) |\varepsilon_u|^2 + \frac{\varepsilon_n^2 + \varepsilon_v^2}{2} - c_k   \left ( \varepsilon_n\varepsilon_v + \text{Im} \left ( \overline{\varepsilon_u} \partial_x \varepsilon_u \right ) \right ) \geqslant \mu_1 \left ( | \partial_x \varepsilon_u |^2 + |\varepsilon_u|^2 + \varepsilon_n^2 + \varepsilon_v^2 \right ) \]
with $\mu_1 = \min\limits_{1 \leqslant k \leqslant K} \left ( \frac{\omega_k^0}{2} + \frac{c_k^2}{4} \, , \frac{1}{2} \right )$. Now, if $c_k \neq 0$, we have $\alpha_k (t) \neq 0$ and $| \text{Im} ( \overline{\varepsilon_u} \partial_x \varepsilon_u) | \leqslant \frac{\alpha_k | \partial_x \varepsilon_u |^2}{2} + \frac{c_k^2 |\varepsilon_u|^2}{2 \alpha_k}$ by Young's inequality. Thus
\begin{align*}
& \, \, | \partial_x \varepsilon_u |^2 + \nu_k (t) |\varepsilon_u|^2 + \frac{\varepsilon_n^2 + \varepsilon_v^2}{2} - c_k   ( \varepsilon_n \varepsilon_v + \text{Im} ( \overline{\varepsilon_u} \partial_x \varepsilon_u ) ) \\
\geqslant & \, \, | \partial_x \varepsilon_u |^2 + \left ( \omega_k (t) + \frac{c_k^2}{4} \right ) |\varepsilon_u|^2 + \frac{\varepsilon_n^2 + \varepsilon_v^2}{2} - |c_k| \frac{\varepsilon_n^2 + \varepsilon_v^2}{2} - \frac{\alpha_k | \partial_x \varepsilon_u |^2}{2} - \frac{c_k^2 |\varepsilon_u|^2}{2 \alpha_k}  \\
\geqslant & \, \, \left ( 1 - \frac{\alpha_k}{2} \right ) | \partial_x \varepsilon_u |^2 + (1-|c_k|) \frac{\varepsilon_n^2 + \varepsilon_v^2}{2} + \left ( \omega_k(t) + \frac{c_k^2}{4} - \frac{c_k^2}{2 \alpha_k} \right ) |\varepsilon_u|^2 \\
\geqslant & \, \, \frac{2 \omega_k(t)}{2 \omega_k(t) + c_k^2} | \partial_x \varepsilon_u |^2 + (1-|c_k|) \frac{\varepsilon_n^2 + \varepsilon_v^2}{2} + \frac{\omega_k(t)}{2} |\varepsilon_u|^2  \\
\geqslant & \, \, \mu \left ( | \partial_x \varepsilon_u |^2 + |\varepsilon_u|^2 + \varepsilon_n^2 + \varepsilon_v^2 \right )
\end{align*}
where $\mu_2 := \min\limits_{1 \leqslant k \leqslant K} \left ( \frac{2 \frac{\omega_k^0}{2}}{2 \frac{\omega_k^0}{2} + c_k^2} \, , \frac{\omega_k^0 /2}{2} \, , \frac{1-|c_k|}{2} \right ) > 0$. Recall indeed \eqref{omega}. Take $\mu = \min ( \mu_1 \, , \mu_2 )$. Hence,
\[ \forall k \in \{ 1 \, , ... \, , K \} , \, \, \, | \partial_x \varepsilon_u |^2 + \nu_k (t) |\varepsilon_u|^2 + \frac{\varepsilon_n^2 + \varepsilon_v^2}{2} - c_k   \left ( \varepsilon_n \varepsilon_v + \text{Im} \left ( \overline{\varepsilon_u} \partial_x \varepsilon_u \right ) \right ) \geqslant \mu \left (  | \partial_x \varepsilon_u |^2 +  |\varepsilon_u|^2 + \varepsilon_n^2 + \varepsilon_v^2 \right ) \]
and this estimate still holds if $\mu$ is replaced by a smaller positive value. The bounds
\begin{align*}
& \chi_k(t) - \Phi_B (x-c_kt - \sigma_k (t)) \geqslant - e^{- \frac{L}{4B}} \\
\text{and} \, \, \, & 0 \leqslant \mu \left ( | \partial_x \varepsilon_u |^2 + |\varepsilon_u|^2 + \varepsilon_n^2 + \varepsilon_v^2 \right ) \leqslant | \partial_x \varepsilon_u |^2 + \nu_k (t) |\varepsilon_u|^2 + \frac{\varepsilon_n^2 + \varepsilon_v^2}{2} - c_k   ( \varepsilon_n \varepsilon_v + \text{Im} ( \overline{\varepsilon_u} \partial_x \varepsilon_u ) ) \\ 
& \midspace \midspace \midspace \midspace \midspace \midspace \midspace \midspace \midspace \midspace \midspace \midspace \midspace \leqslant C \left ( | \partial_x \varepsilon_u |^2 + |\varepsilon_u|^2 + \varepsilon_n^2 + \varepsilon_v^2 \right )
\end{align*}
show that
\begin{align*}
& \int_{\R} ( \chi_k (t) - \Phi_B (x_k^t) ) \left ( | \partial_x \varepsilon_u |^2 + \nu_k (t) |\varepsilon_u|^2 + \frac{\varepsilon_n^2 +  \varepsilon_v^2}{2} - c_k   (\varepsilon_n \varepsilon_v + \text{Im} (  \overline{\varepsilon_u} \partial_x \varepsilon_u )) \right ) \\
\geqslant & \, \, \mu \int_{\R} ( \chi_k (t) - \Phi_B (x_k^t) ) \left ( | \partial_x \varepsilon_u |^2 + |\varepsilon_u|^2 + \varepsilon_n^2 +\varepsilon_v^2 \right ) - C e^{- \frac{L}{4B}} \int_{\R} \left ( | \partial_x \varepsilon_u |^2 + |\varepsilon_u|^2 + \varepsilon_n^2 + \varepsilon_v^2 \right ).
\end{align*}
It follows that
\begin{align*} 
G_{2,1} (t) \geqslant & \, \sum\limits_{k=1}^K C \int_{\R} \Phi_B (x_k^t) \left ( | \partial_x \varepsilon_u |^2 + |\varepsilon_u|^2 + \varepsilon_n^2 + \varepsilon_v^2 \right )  \\
& \midspace + \sum\limits_{k=1}^K \left [ \mu \int_{\R} \left ( \chi_k (t) - \Phi_B (x_k^t) \right ) \left ( | \partial_x \varepsilon_u |^2 + |\varepsilon_u|^2 + \varepsilon_n^2 + \varepsilon_v^2 \right ) \right. \\
& \midspace \midspace \midspace \midspace \midspace \left. - C e^{- \frac{L}{4B}} \int_{\R} \left ( | \partial_x \varepsilon_u |^2 + |\varepsilon_u|^2 + \varepsilon_n^2 + \varepsilon_v^2 \right ) \right ] \\
\geqslant & \, \, (C - \mu ) \sum\limits_{k=1}^K \int_{\R} \Phi_B (x_k^t) \left ( | \partial_x \varepsilon_u |^2 + |\varepsilon_u|^2 + \varepsilon_n^2 + \varepsilon_v^2 \right ) \\
& \midspace + \mu \sum\limits_{k=1}^K \int_{\R} \chi_k (t) \left ( | \partial_x \varepsilon_u |^2 + |\varepsilon_u|^2 + \varepsilon_n^2 + \varepsilon_v^2 \right ) - C e^{- \frac{L}{4B}} \int_{\R} \left ( | \partial_x \varepsilon_u |^2 + |\varepsilon_u|^2 + \varepsilon_n^2 + \varepsilon_v^2 \right ).
\end{align*}
We can replace $\mu$ by a smaller $\mu >0$ so that $C - \mu > 0$. Using $\sum\limits_{k=1}^K \chi_k = 1$, we find
\[ G_{2,1} (t) \geqslant \left ( \mu - C e^{- \frac{L}{4B}} \right ) \int_{\R} \left ( | \partial_x \varepsilon_u |^2 + |\varepsilon_u|^2 + \varepsilon_n^2 + \varepsilon_v^2 \right ) \geqslant \frac{\mu}{2} || (\varepsilon_u \, , \varepsilon_n \, , \varepsilon_v) ||_{\mathbf{H}}^2 \]
taking $L>0$ large enough. \hfill \qedsymbol

\subsection{Conclusion}
\noindent Let us gather the controls of $G_0$, $G_1$, $G_2$ and $G_3$ in order to get a control of the whole Weinstein functional $G(t)$. This leads to a good estimate on $||( \varepsilon_u \, , \varepsilon_n \, , \varepsilon_v ) ||_{\mathbf{H}}$. 

\begin{leftbar}
\noindent \textbf{Lemma 9.} For all $t \in [t^* \, , T_p]$,
\[ ||( \varepsilon_u \, , \varepsilon_n \, , \varepsilon_v ) ||_{\mathbf{H}}^2 \leqslant C \left ( \frac{A_0^2}{L} + 1 \right ) e^{-2 \theta_0 t} \]
provided that $T_0>0$ is chosen large enough (depending on $A_0$ and $L$).
\end{leftbar}

\noindent \textit{Proof.} From Lemma 1 we know that $| \mathcal{M}_k (t) - \mathcal{M}_k (T_p)| + | \mathcal{P}_k (t) - \mathcal{P}_k (T_p)| \leqslant \frac{CA_0^2}{L} e^{-2 \theta_0 t}$ for all $k \in \{ 1 \, , ... \, , K \}$. On the other hand, we know that the energy $E(t) = \int_{\R} \left ( | \partial_x u |^2 + n |u|^2 + \frac{n^2 + v^2}{2} \right )$ is preserved: $E (t) = E (T_p)$. Thus, the Weinstein functional $G(t) = E (t) + \sum\limits_{k=1}^K \left ( \nu_k^0 \mathcal{M}_k (t) - c_k   \mathcal{P}_k (t) \right )$ satisfies
\[ G(t) - G(T_p) \leqslant | G(t) - G(T_p) | \leqslant \frac{CA_0^2}{L} e^{-2 \theta_0 t}. \]
Since $\varepsilon_u (T_p) = 0$, $\varepsilon_n (T_p) = 0$ and $\varepsilon_v (T_p)=0$, we have $G(T_p) = G_0 (T_p)$. We write
\[ \frac{CA_0^2}{L} e^{-2 \theta_0 t} \geqslant G(t)-G(T_p) \geqslant G_{2,1} (t) - | G_0 (t)-G_0 (T_p)| - |G_1(t)| - |G_{2,2} (t)| - |G_3(t)|. \]
The quantities $G_{2,1} (t)$, $|G_0 (t) - G_0(T_p)|$ and $|G_1(t)|$ are estimates respectively in Lemmas 8, 4 and 5. The quantity $|G_{2,2} (t)|$ is estimated as follows:
\[ |G_{2,2}| \leqslant C \sum\limits_{k=1}^K | \omega_k^0 - \omega_k (t)| \, || \varepsilon_u ||_{L^2}^2 \leqslant C A_0^3 e^{-3 \theta_0 t} \leqslant C e^{-2 \theta_0 t} \]
for $T_0>0$ large enough. The quantity $| G_3 (t)|$ is estimated as follows:
\[ |G_3| \leqslant C || ( \varepsilon_u \, , \varepsilon_n \, , \varepsilon_v ) ||_{\mathbf{H}}^3 \leqslant C A_0^3 e^{-3 \theta_0 t} \leqslant C e^{-2 \theta_0 t} \]
for $T_0>0$ large enough. Gathering all these estimates we find
\begin{align*} ||(\varepsilon_u \, , \varepsilon_n \, , \varepsilon_v)||_{\mathbf{H}}^2 \leqslant & \, \, C G_{2,1} (t)  \\
\leqslant & \, \, | G_0 (t)-G_0 (T_p)| + |G_1(t)| + |G_{2,2} (t)| + |G_3(t)| + \frac{CA_0^2}{L} e^{-2 \theta_0 t} \\
\leqslant & \, \, C e^{-2 \theta_0 t} + C e^{-2 \theta_0 t} + C e^{-2 \theta_0 t} + C e^{-2 \theta_0 t} + \frac{CA_0^2}{L} e^{-2 \theta_0 t} \\
\leqslant & \, \, C \left ( \frac{A_0^2}{L} +1 \right ) e^{-2 \theta_0 t}
\end{align*}
which is the desired result.
\hfill \qedsymbol

\noindent \\ From this control on $|| \varepsilon ||_{\mathbf{H}}$ we deduce the following control on the parameters.

\begin{leftbar}
\noindent \textbf{Lemma 10.} For all $t \in [t^* \, , T_p ]$, 
\[ \sum\limits_{k=1}^K \left ( | \omega_k (t) - \omega_k^0 |^2 + | \sigma_k (t) - \sigma_k^0 |^2 + | \gamma_k (t) - \gamma_k^0 |^2 \right )  \leqslant C \left ( \frac{A_0^2}{L} + 1 \right ) e^{-2 \theta_0 t} \]
provided that $T_0$ is chosen large enough (depending on $A_0$). 
\end{leftbar}

\noindent \textit{Proof.} Take $k \in \{ 1 \, , ... \, , K \}$. From Proposition 1 we know that
\[ | \dot{\omega}_k (t) |^2 + | \dot{\sigma}_k (t) |^2 \leqslant C || \varepsilon (t) ||_{\mathbf{H}}^2 + C e^{-2 \theta_0 t} \leqslant C \left ( \frac{A_0^2}{L} + 1 \right ) e^{-2 \theta_0 t} \]
thanks to Lemma 9. Integrating this inequality between $t$ and $T_p$ (recall that $\omega_k (T_p) = \omega_k^0$ and $\sigma_k (T_p) = \sigma_k^0$) we find $| \omega_k (t) - \omega_k^0 |^2 + | \sigma_k (t) - \sigma_k^0 |^2 \leqslant C \left ( \frac{A_0^2}{L} +1 \right ) e^{-2 \theta_0 t}$. Alternatively, to control $\omega_k - \omega_k^0$ we could have used Lemma 3 and Lemma 9. Now, from Proposition 1 we know that
\begin{align*} 
| \dot{\gamma}_k (t) |^2 \leqslant \, \, & 2 | \dot{\gamma}_k (t) - ( \omega_k (t) - \omega_k^0 ) |^2 + 2 | \omega_k (t) - \omega_k^0 |^2 \\
\leqslant \, \, & C || \varepsilon (t) ||_{\mathbf{H}}^2 + C e^{-2 \theta_0 t} + C | \omega_k (t) - \omega_k^0 |^2 \\
\leqslant \, \, & C \left ( \frac{A_0^2}{L} +1 \right ) e^{-2 \theta_0 t}.
\end{align*}
Integrating this inequality between $t$ and $T_p$ (recall that $\gamma_k (T_p) = \gamma_k^0$), we find $| \gamma_k (t) - \gamma_k^0 |^2 \leqslant C \left ( \frac{A_0^2}{L} +1 \right ) e^{-2 \theta_0 t}$ which concludes the proof. \hfill \qedsymbol

\noindent \\ \textbf{From the beginning of section 1 until now}, we have assumed that $||(U_p \, , N_p \, , V_p)(t)||_{\mathbf{H}} \leqslant A_0 e^{- \theta_0 t}$ for all $t \in [t^* \, , T_p]$. Proposition 1 and Lemmas 1 to 10 hold under this hypothesis. We cease to assume it from now on. The following proposition gather the results of the previous and current sections in order to prove a bootstrap argument for the control of $||(U_p \, , N_p \, , V_p)(t)||_{\mathbf{H}}$.

\begin{leftbar}
\noindent \textbf{Proposition 2.} There exist $A_0>0$, $T_0>0$, $p_0>0$ such that, for all $p \geqslant p_0$ and for all $t^* \in [T_0 \, , T_p]$,
\[ \left [ \, \, \forall t \in [t^* \, , T_p] , \, \, || (U_p \, , N_p \, , V_p)(t) ||_{\mathbf{H}} \leqslant A_0 e^{- \theta_0 t} \, \, \right ] \, \, \, \, \, \Longrightarrow \, \, \, \, \, \left [ \, \, \forall t \in [t^* \, , T_p ] , \, \, || (U_p \, , N_p \, , V_p )(t) ||_{\mathbf{H}} \leqslant \frac{A_0}{2} e^{- \theta_0 t} \, \, \right ] . \]
\end{leftbar}

\noindent \textit{Proof.} Assume that $|| (U_p \, , N_p \, , V_p )(t) ||_{\mathbf{H}} \leqslant A_0 e^{- \theta_0 t}$ for all $t \in [t^* \, , T_p]$. Then the results of sections 1 to 3 hold - in particular Lemmas 9 and 10 hold. For all $t \in [t^* \, , T_p ]$, we have
\begin{align*}
|| (R^u - S^u \, , R^n - S^n \, , R^v - S^v ) (t) ||_{\mathbf{H}}^2 \leqslant & \, \, C \sum\limits_{k=1}^K \left ( | \omega_k (t) - \omega_k^0 |^2 + | \sigma_k (t) - \sigma_k^0 |^2 + | \gamma_k (t) - \gamma_k^0 |^2 \right ) \\
\leqslant & \, \, C \left ( \frac{A_0^2}{L} +1 \right ) e^{-2 \theta_0 t}.
\end{align*}
Thus
\[ || (U_p \, , N_p \, , V_p )(t) ||_{\mathbf{H}}^2 \leqslant 2 || \varepsilon_p (t) ||_{\mathbf{H}}^2 + 2 || (R^u - S^u \, , R^n - S^n \, , R^v - S^v ) (t) ||_{\mathbf{H}}^2 \leqslant C \left ( \frac{A_0^2}{L} +1 \right ) e^{-2 \theta_0 t} \]
where the constant $C>0$ does not depend on $A_0$. Choosing $A_0^2 > 8C$ then $L=A_0^2$ and then $T_0>0$ large enough accordingly (see previous Lemmas), it follows that
\[ || (U_p \, , N_p \, , V_p) (t) ||_{\mathbf{H}}^2 \leqslant 2C e^{-2 \theta_0 t} \leqslant \frac{A_0^2}{4} e^{-2 \theta_0 t} \]
for all $t \in [t^* \, , T_p ]$. This is the desired result. \hfill \qedsymbol

\noindent \\ Thanks to this bootstrap lemma, we conclude the following uniform exponential bound on $||(U_p \, , N_p \, , V_p )||_{\mathbf{H}}$ via a classical connexity argument.

\begin{leftbar}
\noindent \textbf{Proposition 3.} There exist constants $C_0 >0$ and $T_0>0$ such that, for all $p \geqslant 1$ and all $t \in [T_0 \, , T_p ]$,
\begin{equation}
|| (U_p \, , N_p \, , V_p) (t) ||_{\mathbf{H}} \leqslant C_0 e^{- \theta_0 t}. 
\label{unifbound}
\end{equation}
\end{leftbar}

\noindent \textit{Proof.} Take $A_0$, $T_0$ and $p_0$ as in Proposition 2. Take $p \geqslant p_0$. Since $t \mapsto || (U_p \, , N_p \, , V_p ) (t) ||_{\mathbf{H}}$ is continuous and $(U_p \, , N_p \, , V_p)(T_p) = 0$, we have $|| (U_p \, , N_p \, , V_p) (t) ||_{\mathbf{H}} \leqslant A_0 e^{- \theta_0 t}$ on $[T_p - \tau \, , T_p ]$ for a certain $\tau >0$. Define
\[ t^* = \inf \left \{ t_* \in [T_0 \, , T_p] \, \, | \, \, \forall t \in [t_* \, , T_p] , \, || (U_p \, , N_p \, , V_p) (t) ||_{\mathbf{H}} \leqslant A_0 e^{- \theta_0 t} \right \} < T_p . \]
Suppose that $t^* > T_0$. By continuity, $|| ( U_p \, , N_p \, , V_p) (t) ||_{\mathbf{H}} \leqslant A_0 e^{- \theta_0 t}$ for all $t \in [t^* \, , T_0 ]$. By Proposition 2, $|| (U_p \, , N_p \, , V_p) (t) ||_{\mathbf{H}} \leqslant \frac{A_0}{2} e^{- \theta_0 t}$ for all $t \in [t^* \, , T_0 ]$. By continuity, we have $|| (U_p \, , N_p \, , V_p)(t) ||_{\mathbf{H}} \leqslant A_0 e^{- \theta_0 t}$ on $[t^* - \tau ' \, , T_p ]$, for a certain $\tau ' >0$. This contradicts the definition of $t^*$. Hence $t^* = T_0$. \\
\\ Take $C_0 = A_0$. For all $p \geqslant p_0$, \eqref{unifbound} holds on $[T_0 \, , T_p]$. We can replace $C_0$ and $T_0$ by larger values to ensure that \eqref{unifbound} holds on $[T_0 \, , T_p]$ for all $p \geqslant 1$. \hfill \qedsymbol

\noindent \\ Since the $\mathbf{H}$ norm of $(R^u \, , R^n \, , R^v)(t)$ is uniformly bounded, Proposition 3 shows that there exists $C>0$ such that
\begin{equation}
\forall p \geqslant 1 , \, \forall t \in [T_0 \, , T_p] , \, \, \, || (u_p \, , n_p \, , v_p) (t) ||_{\mathbf{H}} \leqslant C. 
\label{bounded}
\end{equation}

\section{Compactness argument}
\noindent We now aim at using compactness arguments to extract a converging subsequence $(u_{\psi (p)} \, , n_{\psi (p)} \, , v_{\psi (p)} )$. The issue is the following: since $(u_p \, , n_p \, , v_p)$ are uniformly bounded in $H^1 \times L^2 \times L^2$, Rellich's compactness theorem enables to extract converging subsequences for example in $L_{\text{loc}}^2 \times H_{\text{loc}}^{-1} \times H_{\text{loc}}^{-1}$. To conclude we would need a Cauchy theory in a space of lower regularity than $H^1 \times L^2 \times L^2$. Such a theory exists: see \cite{Gin}, the Cauchy problem for the Zakharov system is well-posed in $H^{1/2} \times H^{-1/2} \times H^{-1/2}$ (in dimension $1$) for example. However working in such a space would require additional estimates. See \cite{Gin}, \cite{Bou}, \cite{Be} or \cite{Sa} for Cauchy theory in lower regularity spaces. Instead of using the Cauchy theory in a lower regularity space, we choose to prove higher regularity estimates in order to use the Cauchy theory in the energy space $H^1 \times L^2 \times L^2$. As a side consequence, we gain additional regularity on the solution. The point of section 4.1 is to establish uniform exponential bounds in $H^2 \times H^1 \times H^1$ thanks to the use of modified energies. See \cite{Ka} and \cite{Ma1} for a similar use of modified energies for the gKdV equation and \cite{Ber} for modified energies for the NLS equation; in both cases the use of suitable modified energies gives access to better regularity on the solution. Such modified energies may have applications in a different framework, for example to prove global wellposedness in higher Sobolev spaces. 

\subsection{Modified energies}
\noindent We keep on denoting $(U_p \, , N_p \, , V_p )$ simply by $(U \, , N \, , V)$ in this section 4.1. The functions $(U \, , N \, , V)$ satisfy the system
\[ \left \{ \begin{array}{ccl} \partial_t U &=& i \partial_x^2 U - i ( R^u N + R^n U + NU ) - i \Psi_R^u \\ \partial_t N &=& - \partial_x V \\ \partial_t V &=& - \partial_x N - \partial_x (|U|^2) - 2 \partial_x \text{Re} (R^u \overline{U}) - \Psi_R^v \end{array} \right. \]
where $\Psi_R^u := \sum\limits_{1 \leqslant j \neq k \leqslant K} R_j^u R_k^n$ and $\Psi_R^v := \partial_x \left ( \left | \sum\limits_{k=1}^K R_k^u \right |^2 - \sum\limits_{k=1}^K |R_k^u|^2 \right ) = \text{Re} \sum\limits_{1 \leqslant j \neq k \leqslant K} \partial_x ( R_j^u \overline{R_k^u} )$. 

\begin{leftbar}
\noindent \textbf{Lemma 11.} Define the functionals
\begin{align*} 
\mathbb{H} (t) =& \int_{\R} \left ( | \partial_x^2 U |^2 + \frac{( \partial_x N )^2}{2} + \frac{( \partial_x V )^2}{2} \right ) \\
\text{and} \, \, \, \mathbb{G} (t) =& \int_{\R} \left ( | \partial_x^2 U |^2 + \frac{( \partial_x N )^2}{2} + \frac{( \partial_x V )^2}{2} \right ) + 2 \int_{\R} N | \partial_x U |^2 + 2 \Ree \int_{\R} U \left ( \partial_x N   \partial_x \overline{U} \right )  \\
& \midspace + 2 \Ree \int_{\R} R^u \left ( \partial_x N   \partial_x \overline{U} \right ) - 2 \Ree \int_{\R} \overline{U} \left (  \partial_x R^u   \partial_x N \right ) . 
\end{align*}
For all $t \in [T_0 \, , T_p]$, 
\begin{equation}
| \mathbb{G} '(t) | \leqslant C \left ( |\mathbb{G} (t)|^{3/4} + 1 \right ) e^{- \theta_0 t}
\label{edoG}
\end{equation}
and
\begin{equation}
\mathbb{H} \leqslant C \left ( |\mathbb{G}| + e^{- \theta_0 t} \right ).
\label{GH}
\end{equation}
\end{leftbar}

\noindent \textit{Proof.} The quantity we seek to control is $\mathbb{H}$, but no such conservation law holds. Differentiating $\mathbb{H}$ with regards to $t$, nonlinear terms of higher derivative order appear in $\mathbb{H} '(t)$, for example $\partial_x^3 U$ or $\partial_x^2 N$. To counterweight these unwanted terms, we add progressively the other terms in $\mathbb{G}$. In return, other nonlinear terms naturally appear when differentiating $\mathbb{G}$ with regards to $t$, but with reasonable derivative orders, and they can all be controled via Sobolev and Gagliardo-Nirenberg inequalities. Take $\kappa >0$. We compute
\begin{align*} 
& \frac{\text{d}}{\text{d}t} \int_{\R} \left ( | \partial_x^2 U |^2 + \frac{( \partial_x N )^2}{\kappa} + \frac{( \partial_x V )^2}{\kappa} \right ) \\
=& \, \, 2 \Ree \int_{\R} \partial_x^2 \left ( i \partial_x^2 U - i (R^u N + R^n U + NU) - i \Psi_R^u \right ) \overline{\partial_x^2 U} \underbrace{- \frac{2}{\kappa} \int_{\R} \partial_x (\partial_x V)   \partial_x N - \frac{2}{\kappa} \int_{\R} (\partial_x^2 N) \partial_x V}_{= \, 0} \\
& \midspace - \frac{2}{\kappa} \int_{\R} ( \partial_x V ) \partial_x   \left ( \partial_x (|U|^2) + 2 \partial_x \Ree (R^u \overline{U}) + \Psi_R^v \right ) \\
=& -2 \Imm \int_{\R} (\partial_x^4 U) \overline{\partial_x^2 U} +2 \Imm \int_{\R} ((\partial_x^2 R^u) N + 2 (\partial_x R^u )   ( \partial_x N) + ( \partial_x^2 R^n) U + 2 \partial_x R^n   \partial_x U ) \overline{\partial_x^2 U} \\
& \midspace + 2 \Imm \int_{\R} R^u ( \partial_x^2 N) \overline{\partial_x^2 U} \underbrace{+2 \Imm \int_{\R} (R^n+N) ( \partial_x^2 U ) \overline{\partial_x^2 U}}_{= \, 0} +2 \Imm \int_{\R} ( (\partial_x^2 N) U + 2 \partial_x N   \partial_x U) \overline{\partial_x^2 U} \\
& \midspace + 2 \Imm \int_{\R} ( \partial_x^2 \Psi_R^u ) \overline{\partial_x^2 U}  - \frac{2}{\kappa} \int_{\R} (\partial_x V) ( 2 | \partial_x U |^2 + 2 \Ree ( \overline{U} \partial_x^2 U )) - \frac{4}{\kappa} \Ree \int_{\R} \partial_x^2 (R^u \overline{U} ) \partial_x V \\
& \midspace - \frac{2}{\kappa} \int_{\R} ( \partial_x V ) ( \partial_x   \Psi_R^v ) \\
=& \, \, \underbrace{2 \Imm \int_{\R} | \partial_x^3 U |^2}_{= \, 0} + 2 \Imm \int_{\R} ( \partial_x^2 R^u) N \overline{\partial_x^2 U} + 4 \Imm \int_{\R} ( \partial_x R^u   \partial_x N ) \overline{\partial_x^2 U} + 2 \Imm \int_{\R} ( \partial_x^2 R^n ) U \overline{\partial_x^2 U} \\
& \midspace +4 \Imm \int_{\R} ( \partial_x R^n   \partial_x U ) \overline{\partial_x^2 U} - 2 \Imm \int_{\R} (\partial_x R^u   \partial_x N) \overline{\partial_x^2 U} - 2 \Imm \int_{\R} R^u \partial_x N   \left ( \overline{\partial_x^3 U} \right ) \\
& \midspace +2 \Imm \int_{\R} (\partial_x N   \partial_x U ) \overline{ \partial_x^2 U}  - 2 \Imm \int_{\R} U \partial_x N   \left ( \overline{\partial_x^3 U} \right ) + 2 \Imm \int_{\R} ( \partial_x^2 \Psi_R^u) \overline{\partial_x^2 U} - \frac{4}{\kappa} \int_{\R} ( \partial_x V ) | \partial_x U |^2 \\
& \midspace - \frac{4}{\kappa} \Ree \int_{\R} ( \partial_x V ) \overline{U} \partial_x^2 U - \frac{4}{\kappa} \Ree \int_{\R} (\partial_x^2 R^u) \overline{U} (\partial_x V) - \frac{8}{\kappa} \Ree \int_{\R} (\partial_x R^u   \partial_x \overline{U}) ( \partial_x V ) \\
& \midspace - \frac{4}{\kappa} \Ree \int_{\R} R^u \overline{\partial_x^2 U} ( \partial_x V ) - \frac{2}{\kappa} \int_{\R} ( \partial_x V ) ( \partial_x   \Psi_R^v ) \\
=& \, \, 2 \Imm \int_{\R} ( \partial_x^2 R^u) N \overline{\partial_x^2 U} \, \boxed{+ 2 \Imm \int_{\R} ( \partial_x R^u   \partial_x N ) \overline{\partial_x^2 U}} + 2 \Imm \int_{\R} ( \partial_x^2 R^n ) U \overline{\partial_x^2 U} +4 \Imm \int_{\R} ( \partial_x R^n   \partial_x U ) \overline{\partial_x^2 U}  \\
& \midspace \, \boxed{- 2 \Imm \int_{\R} R^u \partial_x N   \left (  \overline{\partial_x^3 U} \right )} \,  \boxed{+2 \Imm \int_{\R} (\partial_x N   \partial_x U ) \overline{ \partial_x^2 U}} \, \boxed{- 2 \Imm \int_{\R} U \partial_x N   \left ( \overline{\partial_x^3 U} \right )} \\
& \midspace + 2 \Imm \int_{\R} ( \partial_x^2 \Psi_R^u) \overline{\partial_x^2 U} - \frac{4}{\kappa} \int_{\R} ( \partial_x V ) | \partial_x U |^2 \, \boxed{- \frac{4}{\kappa} \Ree \int_{\R} ( \partial_x V ) \overline{U} \partial_x^2 U} - \frac{4}{\kappa} \Ree \int_{\R} (\partial_x^2 R^u) \overline{U} (\partial_x V) \\
& \midspace - \frac{8}{\kappa} \Ree \int_{\R} (\partial_x R^u   \partial_x \overline{U}) ( \partial_x V )  \, \boxed{- \frac{4}{\kappa} \Ree \int_{\R} R^u \overline{\partial_x^2 U} ( \partial_x V )} - \frac{2}{\kappa} \int_{\R} ( \partial_x V ) ( \partial_x   \Psi_R^v ).
\end{align*}
We compute
\begin{align*} 
& \frac{\text{d}}{\text{d}t} \int_{\R} N | \partial_x U |^2 \\
=& \int_{\R} (- \partial_x V) | \partial_x U |^2 +2 \Ree \int_{\R} N \partial_x( i \partial_x^2 U - i (R^u N + R^n U + NU ) - i \Psi_R^u )   \partial_x \overline{U} \\
=& - \int_{\R} (\partial_x V) | \partial_x U |^2 - 2 \Imm \int_{\R} N \partial_x ( \partial_x^2 U )   \partial_x \overline{U} + 2 \Imm \int_{\R} N^2 \partial_x R^u   \partial_x \overline{U} + 2 \Imm \int_{\R} R^u N \partial_x N   \partial_x \overline{U} \\
& \midspace + 2 \Imm \int_{\R} NU \partial_x R^n   \partial_x \overline{U} \underbrace{+2 \Imm \int_{\R} NR^n \partial_x U   \partial_x \overline{U}}_{= \, 0} + 2 \Imm \int_{\R} NU \partial_x N   \partial_x \overline{U} \underbrace{+ 2 \Imm \int_{\R} N^2 \partial_x U   \partial_x \overline{U}}_{= \, 0}  \\
& \midspace + 2 \Imm \int_{\R} N \partial_x \Psi_R^u   \partial_x \overline{U} \\
=& - \int_{\R} ( \partial_x V ) | \partial_x U |^2 \, \boxed{+ 2 \Imm \int_{\R} ( \partial_x^2 U ) \partial_x N   \partial_x \overline{U}} \, \underbrace{+2 \Imm \int_{\R} N \partial_x^2 U \left ( \overline{\partial_x^2 U} \right )}_{= \, 0} + 2 \Imm \int_{\R} N^2 \partial_x R^u   \partial_x \overline{U} \\
& \midspace +2 \Imm \int_{\R} NR^u \partial_x N   \partial_x \overline{U} + 2 \Imm \int_{\R} NU \partial_x R^n   \partial_x \overline{U} + 2 \Imm \int_{\R} NU \partial_x N   \partial_x \overline{U} + 2 \Imm \int_{\R} N \partial_x \Psi_R^u   \partial_x \overline{U}.
\end{align*}
We compute
\begin{align*}
& \frac{\text{d}}{\text{d}t} \Ree \int_{\R} U \left ( \partial_x N   \partial_x \overline{U} \right ) \\
=& \, \Ree \int_{\R} U (-\partial_x ( \partial_x V))   \partial_x \overline{U} + \Ree \int_{\R} i ( \partial_x^2 U - R^u N - R^n U - NU - \Psi_R^u) \partial_x N   \partial_x \overline{U} \\
& \midspace + \Ree \int_{\R} U \partial_x N   \partial_x \left ( -i \overline{\partial_x^2 U} + i \overline{R^u} N + i R^n \overline{U} + i N \overline{U} + i \overline{\Psi_R^u} \right ) \\
=& \, \Ree \int_{\R} ( \partial_x V ) \partial_x U   \partial_x \overline{U} + \Ree \int_{\R} U \left ( \overline{\partial_x^2 U} \right ) ( \partial_x V ) - \Imm \int_{\R} (\partial_x^2 U ) \partial_x N   \partial_x \overline{U} + \Imm \int_{\R} R^u N \partial_x N   \partial_x \overline{U} \\
& \midspace + \Imm \int_{\R} R^n U \partial_x N   \partial_x \overline{U} + \Imm \int_{\R} N U \partial_x N   \partial_x \overline{U} + \Imm \int_{\R} \Psi_R^u \partial_x N   \partial_x \overline{U} + \Imm \int_{\R} U \partial_x N   \partial_x \left ( \overline{\partial_x^2 U} \right ) \\
& \midspace - \Imm \int_{\R} NU \partial_x N   \partial_x \overline{R^u} - \Imm \int_{\R} U \overline{R^u} \partial_x N   \partial_x N \underbrace{- \Imm \int_{\R} |U|^2 \partial_x N   \partial_x R^n}_{= \, 0} - \Imm \int_{\R} U R^n \partial_x N   \partial_x \overline{U} \\
& \midspace \underbrace{- \Imm \int_{\R} |U|^2 ( \partial_x N )^2}_{= \, 0} - \Imm \int_{\R} NU \partial_x N   \partial_x \overline{U} - \Imm \int_{\R} U \partial_x N   \partial_x \overline{\Psi_R^u} \\
=& \int_{R^d} (\partial_x V ) | \partial_x U |^2 \, \boxed{+ \Ree \int_{\R} U \left ( \overline{\partial_x^2 U} \right ) ( \partial_x V )} \, \boxed{- \Imm \int_{\R} (\partial_x^2 U ) \partial_x N   \partial_x \overline{U}} + \Imm \int_{\R} R^u N \partial_x N   \partial_x \overline{U} \\
& \midspace + \Imm \int_{\R} \Psi_R^u \partial_x N   \partial_x \overline{U} \, \boxed{+ \Imm \int_{\R} U \partial_x N   \left (  \overline{\partial_x^3 U} \right )} - \Imm \int_{\R} NU \partial_x N   \partial_x \overline{R^u} - \Imm \int_{\R} ( \partial_x N )^2 U \overline{R^u} \\
& \midspace - \Imm \int_{\R} U \partial_x N   \partial_x \overline{\Psi_R^u}.
\end{align*}
We compute
\begin{align*}
& \frac{\text{d}}{\text{d}t} \Ree \int_{\R} R^u \left ( \partial_x N   \partial_x \overline{U} \right ) \\
=& \, \Ree \int_{\R} R^u (- \partial_x ( \partial_x V ))   \partial_x \overline{U} + \Ree \int_{\R} (i \partial_x^2 R^u - i R^n R^u + i \Psi_R^u ) \partial_x N   \partial_x \overline{U} \\
& \midspace + \Ree \int_{\R} R^u \partial_x N   \partial_x \left ( -i \overline{\partial_x^2 U} + i \overline{R^u} N + i R^n \overline{U} + i N \overline{U} + i \overline{\Psi_R^u} \right ) \\
=& \, \Ree \int_{\R} (\partial_x V ) \partial_x R^u   \partial_x \overline{U} + \Ree \int_{\R} (\partial_x V) R^u \overline{\partial_x^2 U} - \Imm \int_{\R} (\partial_x^2 R^u) \partial_x N   \partial_x \overline{U} + \Imm \int_{\R} R^n R^u \partial_x N   \partial_x \overline{U} \\
& \midspace - \Imm \int_{\R} \Psi_R^u \partial_x N   \partial_x \overline{U} + \Imm \int_{\R} R^u \partial_x N   \left ( \overline{\partial_x^3 U} \right ) - \Imm \int_{\R} R^u N \partial_x N   \partial_x \overline{R^u} \underbrace{- \Imm \int_{\R} R^u \overline{R^u} ( \partial_x N )^2}_{= \, 0} \\
& \midspace - \Imm \int_{\R} R^u \overline{U} \partial_x N   \partial_x R^n - \Imm \int_{\R} R^u R^n \partial_x N   \partial_x \overline{U} - \Imm \int_{\R} R^u \overline{U} ( \partial_x N )^2 - \Imm \int_{\R} R^u N \partial_x N   \partial_x \overline{U} \\
&  \midspace - \Imm \int_{\R} R^u \partial_x N   \partial_x \overline{\Psi_R^u} \\
=& \, \Ree \int_{\R} (\partial_x V ) \partial_x R^u   \partial_x \overline{U} \, \boxed{+ \Ree \int_{\R} (\partial_x V) R^u \overline{\partial_x^2 U}} - \Imm \int_{\R} (\partial_x^2 R^u) \partial_x N   \partial_x \overline{U} - \Imm \int_{\R} \Psi_R^u \partial_x N   \partial_x \overline{U} \\
& \midspace  \boxed{+ \Imm \int_{\R} R^u \partial_x N   \left ( \overline{\partial_x^3 U} \right )} - \Imm \int_{\R} R^u N \partial_x N   \partial_x \overline{R^u} - \Imm \int_{\R} R^u \overline{U} \partial_x N   \partial_x R^n - \Imm \int_{\R} R^u \overline{U} ( \partial_x N )^2 \\
& \midspace  - \Imm \int_{\R} R^u N \partial_x N   \partial_x \overline{U}  - \Imm \int_{\R} R^u \partial_x N   \partial_x \overline{\Psi_R^u}.
\end{align*}
At last, we compute
\begin{align*}
& \frac{\text{d}}{\text{d}t} \Ree \int_{\R} \overline{U} \left ( \partial_x R^u   \partial_x N \right ) \\
=& \, \Ree \int_{\R} \overline{U} \partial_x \left ( i \partial_x^2 R^u - i R^u R^n + i \Psi_R^u \right )   \partial_x N + \Ree \int_{\R} \overline{U} \partial_x ( - \partial_x V )   \partial_x R^u \\
& \midspace + \Ree \int_{\R} \left ( - i \overline{\partial_x^2 U} + i \overline{R^u} N + i R^n \overline{U} + iN \overline{U} + i \overline{\Psi_R^u} \right ) \partial_x R^u   \partial_x N \\
=& - \Imm \int_{\R} \overline{U} \left ( \partial_x^3 R^u - R^u \partial_x R^n \right )   \partial_x N + \Imm \int_{\R} \overline{U} R^n \partial_x R^u   \partial_x N - \Imm \int_{\R} \overline{U} \partial_x \Psi_R^u   \partial_x N \\
& \midspace + \Ree \int_{\R} \overline{U} (\partial_x^2 R^u) (\partial_x V) + \Ree \int_{\R} (\partial_x V) \partial_x R^u   \partial_x \overline{U} \, \boxed{+ \Imm \int_{\R} \left ( \overline{\partial_x^2 U} \right ) \partial_x R^u   \partial_x N} \\
& \midspace - \Imm \int_{\R} \overline{R^u} N \partial_x R^u   \partial_x N - \Imm \int_{\R} \overline{U} R^n \partial_x R^u   \partial_x N - \Imm \int_{\R} N \overline{U} \partial_x R^u   \partial_x N \\
& \midspace - \Imm \int_{\R} \overline{\Psi_R^u} \partial_x R^u   \partial_x N.
\end{align*}
Taking $\kappa = 2$, the function $\mathbb{G}$ is defined so that all framed terms (in the previous computations of the derivatives) cancel one another. Other terms also cancel one another but the framed ones are those which would not be well controlled below. This leads to
\[ \mathbb{G}'(t) = \mathbb{G}_R (t) + \mathbb{G}_\Psi (t) + \mathbb{G}_\star (t) \]
where
\begin{align*}
\mathbb{G}_R (t) =& \, \, 2 \Imm \int_{\R} N ( \partial_x^2 R^u ) \left ( \overline{\partial_x^2 U} \right ) + 2 \Imm \int_{\R} U ( \partial_x^2 R^n ) \left ( \overline{\partial_x^2 U} \right ) + 4 \Imm \int_{\R} ( \partial_x R^n   \partial_x U ) \left ( \overline{\partial_x^2 U} \right ) \\
& \midspace - 4 \Ree \int_{\R} ( \partial_x^2 R^u ) \overline{U} ( \partial_x V ) -4 \Ree \int_{\R} \left ( \partial_x R^u   \partial_x \overline{U} \right ) ( \partial_x V ) + 4 \Imm \int_{\R} \left ( \partial_x R^u   \partial_x \overline{U} \right ) N^2 \\
& \midspace + 4 \Imm \int_{\R} R^u N \left ( \partial_x N   \partial_x \overline{U} \right ) + 4 \Imm \int_{\R} \left ( \partial_x R^n   \overline{\partial_x U} \right ) NU - 4 \Imm \int_{\R} \left ( \partial_x \overline{R^u}   \partial_x N \right ) NU   \\
& \midspace -2 \Imm \int_{\R} ( \partial_x^2 R^u ) \left ( \partial_x N   \partial_x \overline{U} \right ) - 4 \Imm \int_{\R} R^u N \left ( \partial_x \overline{R^u}   \partial_x N \right ) - 4 \Imm \int_{\R} R^u \overline{U} \left ( \partial_x R^n   \partial_x N \right ) \\
& \midspace  + 2 \Imm \int_{\R} \overline{U} \left ( \partial_x^3 R^u \right ) \partial_x N , \\
\mathbb{G}_\Psi (t) =& \, \, 2 \Imm \int_{\R} ( \partial_x^2 \Psi_R^u ) \left ( \overline{\partial_x^2 U} \right ) - \int_{\R} ( \partial_x   \Psi_R^v ) ( \partial_x V ) + 4 \Imm \int_{\R} \left ( \partial_x \Psi_R^u   \partial_x \overline{U} \right ) N \\
& \midspace - 4 \Imm \int_{\R} \left ( \partial_x \overline{\Psi_R^u}   \partial_x N \right ) U + 2 \Imm \int_{\R} \left ( \partial_x \Psi_R^u   \partial_x N \right ) \overline{R^u} + 2 \Imm \int_{\R} \overline{\Psi_R^u} \left ( \partial_x R^u \partial_x N \right ) , \\
\mathbb{G}_\star (t) =&  -2 \int_{\R} ( \partial_x V ) | \partial_x U |^2 + 4 \Imm \int_{\R} NU \left ( \partial_x N   \partial_x \overline{U} \right ).
\end{align*}
Let us deal with $\mathbb{G}_\star$ first. Recall the following Sobolev and Gagliardo-Nirenberg inequalities:
\[ || N ||_{L^\infty} \leqslant C || \partial_x N ||_{L^2}^{1/2} || N ||_{L^2}^{1/2} , \midspace  || \partial_x U ||_{L^\infty} \leqslant C || \partial_x^2 U ||_{L^2}^{1/2} || \partial_x U ||_{L^2}^{1/2} \midspace \text{and} \midspace || U ||_{L^\infty} \leqslant C || U ||_{H^1}. \]
Also recall the estimate \eqref{bounded}. It follows that
\[ \left | \int_{\R} ( \partial_x V ) | \partial_x U |^2 \right | \leqslant  || \partial_x V ||_{L^2} || \partial_x U ||_{L^\infty} || \partial_x U ||_{L^2} \leqslant C \mathbb{H}^{1/2} || \partial_x^2 U ||_{L^2}^{1/2} || \partial_x U ||_{L^2}^{3/2} \leqslant C \mathbb{H}^{3/4} e^{- \frac{3 \theta_0}{2} t}. \]
Similarly,
\[ \left | \int_{\R} NU \left ( \partial_x N   \partial_x \overline{U} \right ) \right | \leqslant || N ||_{L^\infty} || U ||_{L^\infty} || \partial_x N ||_{L^2} || \partial_x \overline{U} ||_{L^2} \leqslant C \mathbb{H}^{3/4} e^{- \frac{5 \theta_0}{2} t}. \]
Thus
\[ | \mathbb{G}_\star (t) | \leqslant C ( \mathbb{H}^{3/4} +1 ) e^{- \theta_0 t}. \]
Now deal with $\mathbb{G}_R$. We know that functions $R^u$ and $R^n$, as well as all of their derivatives, are bounded (regardless of $t$). Many terms in $\mathbb{G}_R$ are estimated using simple Cauchy-Schwarz inequalities:
\begin{align*} 
& \left | \int_{\R} N ( \partial_x^2 R^u ) \left ( \overline{\partial_x^2 U} \right ) \right | + \left | \int_{\R} U ( \partial_x^2 R^n ) \left ( \overline{\partial_x^2 U} \right ) \right | + \left | \int_{\R} \left ( \partial_x R^n   \partial_x U \right ) \left ( \overline{\partial_x^2 U} \right ) \right | \\
& \midspace + \left | \int_{\R} ( \partial_x^2 R^u ) \overline{U} ( \partial_x V ) \right | + \left | \int_{\R} \left ( \partial_x R^u   \partial_x \overline{U} \right ) ( \partial_x V ) \right | + \left | \int_{\R} ( \partial_x^2 R^u ) \left ( \partial_x N   \partial_x \overline{U} \right ) \right | \\
& \midspace + \left | \int_{\R} R^u N \left ( \partial_x \overline{R^u}   \partial_x N \right ) \right | + \left | \int_{\R} R^u \overline{U} \left ( \partial_x R^n   \partial_x N \right ) \right | + \left | \int_{\R} \overline{U} \left ( \partial_x^3 R^u \right ) \partial_x N \right | & & \leqslant \, \, C \mathbb{H}^{1/2} e^{- \theta_0 t}. 
\end{align*}
To estimate the remaining terms in $\mathbb{G}_R$ we use Sobolev/Gagliardo-Nirenberg inequalities:
\begin{align*} 
& \left | \int_{\R} R^u N \left ( \partial_x N   \partial_x \overline{U} \right ) \right | \leqslant C ||N||_{L^\infty} || \partial_x N ||_{L^2} || \partial_x U ||_{L^2} \leqslant C \mathbb{H}^{3/4} e^{- \frac{3 \theta_0}{2} t} , \\
&  \left | \int_{\R} \left ( \partial_x \overline{R^u}   \partial_x N \right ) NU \right | \leqslant C || \partial_x N ||_{L^2} ||N||_{L^2} ||U||_{L^\infty} \leqslant C \mathbb{H}^{1/2} e^{-2 \theta_0 t} , \\
&  \left | \int_{\R} \left ( \partial_x R^u   \partial_x \overline{U} \right ) N^2 \right | \leqslant C || \partial_x U ||_{L^\infty} ||N||_{L^2}^2 \leqslant C \mathbb{H}^{1/4} e^{- \frac{5 \theta_0}{2} t} \\
\text{and} \, \, \, &  \left | \int_{\R} \left ( \partial_x R^n   \partial_x \overline{U} \right ) NU \right | \leqslant C || \partial_x U ||_{L^2} || N ||_{L^2} || U ||_{L^\infty} \leqslant C e^{-3 \theta_0 t} .
\end{align*}
Gathering these estimates we find
\[ | \mathbb{G}_R (t) | \leqslant C ( \mathbb{H}^{3/4} +1 ) e^{- \theta_0 t}. \]
Now deal with $\mathbb{G}_\Psi$. Simply using the fact that $\partial_x \Psi_R^u$ is bounded, we have
\begin{align*}
& \left | \int_{\R} \left ( \partial_x \overline{\Psi_R^u}   \partial_x N \right ) U \right | \leqslant C || \partial_x N ||_{L^2} ||U||_{L^2} \leqslant C \mathbb{H}^{1/2} e^{- \theta_0 t} \\
\text{and} \, \, \, & \left | \int_{\R} \left ( \partial_x \Psi_R^u   \partial_x \overline{U} \right ) N \right | \leqslant || \partial_x U ||_{L^2} ||N||_{L^2} \leqslant C e^{-2 \theta_0 t}.
\end{align*}
To control the other terms in $\mathbb{G}_\Psi$ we recall that $|| \partial_x \Psi_R^u ||_{L^2} + || \partial_x^2 \Psi_R^u ||_{L^2} + || \partial_x   \Psi_R^v ||_{L^2} \leqslant C e^{- \theta_0 t}$. It follows that
\[ \left | \int_{\R} ( \partial_x^2 \Psi_R^u ) \left ( \overline{\partial_x^2 U} \right ) \right | + \left | \int_{\R} ( \partial_x   \Psi_R^v ) ( \partial_x V ) \right | + \left | \int_{\R} \left ( \partial_x \Psi_R^u   \partial_x N \right ) \overline{R^u} \right | + \left | \int_{\R} \overline{\Psi_R^u} \left ( \partial_x N \partial_x R^u \right ) \right | \leqslant C \mathbb{H}^{1/2} e^{- \theta_0 t}. \]
Hence, 
\[ | \mathbb{G}_\Psi (t) | \leqslant C ( \mathbb{H}^{1/2} + 1 ) e^{- \theta_0 t}. \]
Eventually, gathering all these estimates, we find
\[ | \mathbb{G} '(t) | \leqslant C ( \mathbb{H}^{3/4} + 1 ) e^{- \theta_0 t} . \]
Now let us show that $\mathbb{H} \leqslant C ( \mathbb{G} + 1 )$. Let us control the terms in $\mathbb{G}$ (other than $\mathbb{H}$). We have
\begin{align*}
& \left | \int_{\R} N | \partial_x U |^2 \right | \leqslant ||N||_{L^\infty} || \partial_x U ||_{L^2}^2 \leqslant C  \mathbb{H}^{1/4} e^{- \frac{5 \theta_0}{2} t} , \\
&  \left | \int_{\R} U \left ( \partial_x N   \partial_x \overline{U} \right ) \right | \leqslant C ||U||_{L^\infty} || \partial_x N ||_{L^2} || \partial_x U ||_{L^2} \leqslant C \mathbb{H}^{1/2} e^{-2 \theta_0 t} , \\
&  \left | \int_{\R} R^u \left ( \partial_x N   \partial_x \overline{U} \right ) \right | \leqslant C || \partial_x N ||_{L^2} || \partial_x U ||_{L^2} \leqslant C \mathbb{H}^{1/2} e^{- \theta_0 t}  \\
\text{and} \, \, \, & \left | \int_{\R} \overline{U} \left ( \partial_x R^u   \partial_x N \right ) \right | \leqslant C ||U||_{L^2} || \partial_x N ||_{L^2} \leqslant C \mathbb{H}^{1/2} e^{- \theta_0 t}.
\end{align*}
Hence,
\[ |\mathbb{G}| \geqslant \mathbb{H} - C e^{- \theta_0 t} ( \mathbb{H} +1 ) \geqslant (1 - C e^{- \theta_0 T_0}) \mathbb{H} - C e^{- \theta_0 t} \geqslant \frac{1}{2} \mathbb{H} - C e^{- \theta_0 t} \]
for $T_0$ chosen large enough. Thus \eqref{GH} is established:
\[ \mathbb{H} \leqslant C \left ( |\mathbb{G}| + e^{- \theta_0 t} \right ). \]
This leads to $\mathbb{H}^{3/4} \leqslant C ( |\mathbb{G}| +1)^{3/4} \leqslant C ( |\mathbb{G}|^{3/4} +1 )$. Ultimately, $\mathbb{G}$ satisfies the following differential inequality:
\[ | \mathbb{G} '(t) | \leqslant C \left ( |\mathbb{G} (t)|^{3/4} + 1 \right ) e^{- \theta_0 t}. \]
which is precisely \eqref{edoG}. \hfill \qedsymbol

\noindent \\ Using the differential inequality verified by $\mathbb{G}$, we conduct a bootstrap argument to obtain the following control on $\mathbb{G}$ and thus on $\mathbb{H}$, using \eqref{GH}.

\begin{leftbar}
\noindent \textbf{Lemma 12.} For all $t \in [T_0 \, , T_p]$,
\begin{equation}
|| \partial_x^2 U ||_{L^2}^2 + || \partial_x N ||_{L^2}^2 + || \partial_x V ||_{L^2}^2 \leqslant C e^{- \theta_0 t}.
\label{controlH}
\end{equation}
\end{leftbar}

\noindent \textit{Proof.} Let $T_\star = \inf \left \{ t \in [T_0 \, , T_p] \, \, | \, \, |\mathbb{G} (t)| \leqslant e^{- \theta_0 t /2}. \right \}$. Since $\mathbb{G} (T_p) = 0$, we have $T_\star < T_p$ by continuity. Assume $T_\star \neq T_0$. Then $| \mathbb{G}(t)| > e^{- \theta_0 t /2}$ for some $t$ as close as $T_\star^-$ as one wants; thus $| \mathbb{G} (T_\star)| = e^{- \theta_0 T_\star /2}$ by continuity. Meanwhile, on $[T_\star \, , T_p ]$ we have (thanks to \eqref{edoG})
\[ | \mathbb{G} '(t)| \leqslant C \left ( | \mathbb{G} (t) |^{3/4} +1 \right ) e^{- \theta_0 t} \leqslant C \left ( e^{- 3 \theta_0 t /8} +1 \right ) e^{- \theta_0 t} \leqslant C e^{- \theta_0 t}. \]
Integrating on $[T_\star \, , T_p]$ and recalling that $\mathbb{G} (T_p) = 0$, we find that
\[ | \mathbb{G} (T_\star) | \leqslant \int_{T_\star}^{T_p} | \mathbb{G} ' | \leqslant C e^{- \theta_0 T_\star}. \]
The constant $C$ here does not depend on $T_0$. Take $T_0$ large enough such that $C e^{- \theta_0 t} \leqslant \frac{1}{2} e^{- \theta_0 t /2}$ for all $t \geqslant T_0$ (namely, $T_0 \geqslant 2 \ln (2C) / \theta_0$). We have
\[ e^{- \theta_0 T_\star /2} = | \mathbb{G} (T_\star)| \leqslant C e^{- \theta_0 T_\star} \leqslant \frac{1}{2} e^{- \theta_0 T_\star / 2} \]
which is absurd. We conclude that $T_\star = T_0$. Hence, $| \mathbb{G} (t) | \leqslant e^{- \theta_0 t /2}$ for all $t \in [T_0 \, , T_p ]$. Using \eqref{edoG} and integrating back, we even find that $| \mathbb{G} (t) | \leqslant C e^{- \theta_0 t}$ for all $t \in [T_0 \, , T_p]$. Using \eqref{GH} we ultimately find that
\[ || \partial_x^2 U ||_{L^2}^2 + || \partial_x N ||_{L^2}^2 + || \partial_x V ||_{L^2}^2 \leqslant 2 \mathbb{H} \leqslant C e^{- \theta_0 t} \]
which is precisely \eqref{controlH}. \hfill \qedsymbol

\noindent \\ Recall that $(U \, , N \, , V)$ is actually $(U_p \, , N_p \, , V_p)$. Gathering Proposition 3, Lemma 12 and the fact that the $H^2 \times H^1 \times H^1$ norm of $(R^u \, ,R^n \, , R^v)(t)$ is uniformly bounded, we deduce the following result.

\begin{leftbar}
\noindent \textbf{Corollary 1.} For all $t \in [T_0 \, , T_p]$,
\[ ||U_p (t)||_{H^2}^2 + || N_p (t) ||_{H^1}^2 + ||V_p(t)||_{H^1}^2  \leqslant C e^{- \theta_0 t} \]
and
\[ || u_p(t) ||_{H^2} + ||n_p(t)||_{H^1} + || v_p(t)||_{H^1}  \leqslant C. \]
\end{leftbar}

\noindent \textbf{Remark.} The use of higher-order modified energies would certainly enable to prove similar bounds in higher regularity spaces; we do not dig deeper into this matter here.

\subsection{Conclusion}
\noindent Rellich's compactness theorem enables to extract an $(H_{\text{loc}}^1 \times L_{\text{loc}}^2 \times L_{\text{loc}}^2)$-converging subsequence $(u_{\psi (p)} \, , n_{\psi (p)} \, , v_{\psi (p)})$. In order to go from $H_{\text{loc}}^1 \times L_{\text{loc}}^2 \times L_{\text{loc}}^2$ to $H^1 \times L^2 \times L^2$, we need to control (uniformly in $p$) $|u_p|^2$, $|\partial_x u_p|^2$, $n_p^2$ and $v_p^2$ for $|x|$ large. This is the point of the following two lemmas, based on the conservation of the mass and energy.

\begin{leftbar}
\noindent \textbf{Lemma 13.} Let $\eta >0$. There exists $K_0 = K_0 ( \eta ) > 0$ such that, for all $p \geqslant 1$,
\[ \int_{|x| > K_0} |u_p(T_0 \, , x)|^2 \, \text{d}x \leqslant C \eta \midspace \text{and} \midspace \int_{|x| > K_0} \left ( | \partial_x u_p |^2 + n_p |u_p|^2 + \frac{n_p^2 + v_p^2}{2} \right ) (T_0 \, , x) \, \text{d}x \leqslant C \eta . \]
\end{leftbar}

\noindent \textit{Proof.} Start with the first integral. Using Proposition 3, we have
\[ ||(U_p \, , N_p \, , V_p) (t_\infty)||_{\mathbf{H}}^2 \leqslant C e^{-2 \theta_0 t_\infty} \leqslant \eta \]
for $t_\infty$ large enough (depending on $\eta$). Moreover, due to the localization property of $(R^u \, , R^n \, , R^v)$, there exists $A > 0$ (depending on $\eta$) such that $\int_{|x| > A} \left ( |R^u(t_\infty)|^2 + | \partial_x R^u (t_\infty)|^2 + R^n (t_\infty)^2 + R^v ( t_\infty )^2 \right ) \, \text{d}x \leqslant \eta$. It follows that
\begin{equation}
\int_{|x| > A} \left ( |u_p (t_\infty)|^2 + |\partial_x u_p (t_\infty)|^2 + n_p (t_\infty)^2 + v_p (t_\infty)^2 \right ) \, \text{d}x \leqslant 4 \eta.
\label{K1}
\end{equation}
Now take a cut-off function $g \in \mathscr{C}^1 (\R \, , [0 \, , 1] )$ such that $g \equiv 0$ on $(- \infty \, , 1]$, $g \equiv 1$ on $[2 \, , + \infty )$ and $0 < g' < 2$ on $\R$. Take $B > 0$ to be adjusted later. Set $\text{sgn} (x) = \frac{x}{|x|} \in \R$ if $x \neq 0$ and $0$ if $x=0$. Using the identity
\[ \partial_t \left ( |u_p|^2 \right ) = 2 \Ree \left ( \overline{u_p} \partial_t u_p \right ) = 2 \Ree \left ( i \overline{u_p} \partial_x^2 u_p - i n_p u_p \overline{u_p} \right ) = -2 \Imm \left ( \overline{u_p} \partial_x^2 u_p \right ), \]
we compute
\begin{align*}
\frac{\text{d}}{\text{d}t} \int_{\R} |u_p(t)|^2 g \left ( \frac{|x|-A}{B} \right ) \, \text{d}x =& -2 \Imm \int_{\R} \overline{u_p} \partial_x^2 u_p \, g \left ( \frac{|x|-A}{B} \right ) \, \text{d}x \\
=& \, 2 \Imm \int_{\R} \partial_x u_p   \left ( \overline{u_p} \frac{\sgn (x)}{B} g' \left ( \frac{|x|-A}{B} \right ) + \overline{\partial_x u_p} \, g \left ( \frac{|x|-A}{B} \right ) \right ) \, \text{d}x \\
=& \, \frac{2}{B} \Imm \int_{\R} \overline{u_p} \partial_x u_p \,  \sgn (x) g' \left ( \frac{|x|-A}{B} \right ) \, \text{d}x.
\end{align*}
Thus, using \eqref{bounded},
\[ \left | \frac{\text{d}}{\text{d}t} \int_{\R} |u_p(t)|^2 g \left ( \frac{|x|-A}{B} \right ) \, \text{d}x \right | \leqslant \frac{4}{B} ||u_p(t)||_{H^1}^2 \leqslant \frac{4C}{B} \leqslant \frac{\eta}{t_\infty} \]
for $B>0$ chosen large enough. Now we integrate between $T_0$ and $t_\infty$:
\begin{align*} \int_{\R} |u_p (T_0)|^2 g \left ( \frac{|x|-A}{B} \right ) \, \text{d}x - \int_{\R} |u_p (t_\infty)|^2 g \left ( \frac{|x|-A}{B} \right ) \, \text{d}x \leqslant & \int_{T_0}^{t_\infty} \left | \frac{\text{d}}{\text{d}t} \int_{\R} |u_p(t)|^2 g \left ( \frac{|x|-A}{B} \right ) \, \text{d}x \right | \, \text{d}t \\
\leqslant & \, (t_\infty - T_0) \frac{\eta}{t_\infty} \leqslant \eta . 
\end{align*}
Using the fact that $g \left ( \frac{|x|-A}{B} \right ) = 1$ if $|x|>2B+A$ and $g \left ( \frac{|x|-A}{B} \right ) = 0$ if $|x| \leqslant A$, it follows that
\begin{align*} 
\int_{|x| > 2B+A} |u_p(T_0)|^2 \leqslant & \int_{\R} |u_p(T_0)|^2 g \left ( \frac{|x|-A}{B} \right ) \leqslant \eta + \int_{\R} |u_p(t_\infty)|^2 g \left ( \frac{|x|-A}{B} \right ) \\
\leqslant & \, \eta + \int_{|x| > A} |u_p (t_\infty)|^2 \leqslant C \eta
\end{align*}
thanks to \eqref{K1}. Now deal with the second integral of the lemma. We compute
\begin{align*}
& \frac{\text{d}}{\text{d}t} \int_{\R} \left ( |\partial_x u_p |^2 + n_p |u_p|^2 + \frac{n_p^2 + v_p^2}{2} \right  )(t) \, g \left ( \frac{|x|-A}{B} \right ) \, \text{d}x \\
=& \, 2 \Ree \int_{\R} i ( \partial_x \partial_x^2 u_p - i u_p \partial_x n_p - i n_p \partial_x u_p )   \overline{\partial_x u_p} \, g \left ( \frac{|x|-A}{B} \right ) + \int_{\R} (- \partial_x v_p ) |u_p|^2 g \left ( \frac{|x|-A}{B} \right ) \\
& \midspace +2 \Ree \int_{\R} n_p i (\partial_x^2 u_p - n_pu_p) \overline{u_p} \, g \left ( \frac{|x|-A}{B} \right ) + \int_{\R} (- \partial_x v_p ) n_p \, g \left ( \frac{|x|-A}{B} \right ) \\
& \midspace + \int_{\R} \left ( - \partial_x n_p - \partial_x (|u_p|^2) \right )   v_p \, g \left ( \frac{|x|-A}{B} \right ) \\
=& \, -2 \Imm \int_{\R} \partial_x \partial_x^2 u_p   \overline{\partial_x u_p} \, g \left ( \frac{|x|-A}{B} \right ) + 2 \Imm \int_{\R} u_p \partial_x n_p   \overline{\partial_x u_p} \, g \left ( \frac{|x|-A}{B} \right ) \\
& \midspace - \int_{\R} \left ( ( \partial_x v_p ) |u_p|^2 + v_p   \partial_x (|u_p|^2) \right ) g \left ( \frac{|x|-A}{B} \right ) -2 \Imm \int_{\R} n_p (\partial_x^2 u_p) \overline{u_p} \, g \left ( \frac{|x|-A}{B} \right ) \\
& \midspace - \int_{\R} ( (\partial_x v_p) n_p + (\partial_x n_p )  v_p ) g \left ( \frac{|x|-A}{B} \right ) \\
=& \, 2 \Imm \int_{\R} \partial_x^2 u_p \overline{\partial_x^2 u_p} \, g \left ( \frac{|x|-A}{B} \right ) + 2 \Imm \int_{\R} \partial_x^2 u_p \overline{\partial_x u_p}   \frac{\sgn (x)}{B} g' \left ( \frac{|x|-A}{B} \right ) \\
& \midspace + 2 \Imm \int_{\R} u_p \partial_x n_p   \overline{\partial_x u_p} \, g \left ( \frac{|x|-A}{B} \right ) + \int_{\R} v_p |u_p|^2   \frac{\sgn (x)}{B} g' \left ( \frac{|x|-A}{B} \right ) \\
& \midspace +2 \Imm \int_{\R} ( \partial_x n_p   \partial_x u_p ) \overline{u_p} \, g \left ( \frac{|x|-A}{B} \right ) + 2 \Imm \int_{\R} n_p \partial_x u_p    \overline{\partial_x u_p} \, g \left ( \frac{|x|-A}{B} \right ) \\
& \midspace +2 \Imm \int_{\R} n \overline{u_p} \partial_x u_p   \frac{\sgn (x)}{B} g' \left ( \frac{|x|-A}{B} \right ) + \int_{\R} n_pv_p   \frac{\sgn (x)}{B} g' \left ( \frac{|x|-A}{B} \right ) \\
=& \, \frac{1}{B} \int_{\R} \left ( 2 \Imm \left ( \partial_x^2 u_p \overline{\partial_x u_p} \right ) + v_p |u_p|^2 + 2 \Imm \left ( n_p \overline{u_p} \partial_x u_p \right ) + n_pv_p \right )   \sgn (x) g' \left ( \frac{|x|-A}{B} \right ) \, \text{d}x.
\end{align*}
It follows that
\begin{align*}
& \left | \frac{\text{d}}{\text{d}t} \int_{\R} \left ( |\partial_x u_p |^2 + n_p |u_p|^2 + \frac{n_p^2 + v_p^2}{2} \right ) (t) \, g \left ( \frac{|x|-A}{B} \right ) \, \text{d}x \right | \\
\leqslant & \, \frac{C}{B} \left ( ||\partial_x^2 u_p||_{L^2} || \partial_x u_p ||_{L^2} + ||v_p||_{L^\infty} ||u_p||_{L^2}^2 + ||n_p||_{L^2} ||\partial_x u_p ||_{L^2} ||u_p||_{L^\infty} + ||n_p||_{L^2} ||v_p||_{L^2} \right ).
\end{align*}
Using Corollary 1, it leads to
\[ \left | \frac{\text{d}}{\text{d}t} \int_{\R} \left ( |\partial_x u_p |^2 + n_p |u_p|^2 + \frac{n_p^2 + v_p^2}{2} \right ) (t) \, g \left ( \frac{|x|-A}{B} \right ) \, \text{d}x \right | \leqslant \frac{C}{B} \leqslant \frac{\eta}{t_\infty} \]
for $B>0$ chosen large enough. We conclude as previously:
\begin{align*}
\int_{|x|>2B+A} \left ( | \partial_x u_p |^2 + n_p |u_p|^2 + \frac{n_p^2 + v_p^2}{2} \right ) (T_0) \leqslant & \int_{\R} \left ( | \partial_x u_p |^2 + n_p |u_p|^2 + \frac{n_p^2 + v_p^2}{2} \right ) (T_0) g \left ( \frac{|x|-A}{B} \right ) \\
\leqslant & \, \eta + \int_{\R} \left ( | \partial_x u_p |^2 + n_p |u_p|^2 + \frac{n_p^2 + v_p^2}{2} \right ) (t_\infty) g \left ( \frac{|x|-A}{B} \right ) \\
\leqslant & \, \eta + \int_{|x|>A} \left ( | \partial_x u_p |^2 + n_p |u_p|^2 + \frac{n_p^2 + v_p^2}{2} \right ) (t_\infty) \leqslant C \eta
\end{align*}
thanks to \eqref{K1}. We conclude the proof by taking $K_0 = 2B+A$. \hfill \qedsymbol

\begin{leftbar}
\noindent \textbf{Lemma 14.} Take $T_0>0$ large enough. For all $\eta > 0$, there exists $K_0 = K_0 ( \eta ) > 0$ such that, for all $p \geqslant 1$,
\[ \int_{|x| > K_0} \left ( |u_p|^2 + | \partial_x u_p |^2 + n_p^2 + v_p^2 \right ) (T_0 \, , x) \, \text{d}x \leqslant \eta . \]
\end{leftbar}

\noindent \textit{Proof.} By Lemma 13, there exists $K_0 ( \eta ) > 0$ such that
\[ \int_{|x| < K_0} |u_p (T_0 \, , x) |^2 \, \text{d}x \leqslant \frac{\eta}{4} \midspace \text{and} \midspace \int_{|x| > K_0} \left ( | \partial_x u_p |^2 + n_p |u_p|^2 + \frac{n_p^2 + v_p^2}{2} \right ) (T_0 \, , x) \, \text{d}x \leqslant \frac{\eta}{8}. \]
Now write that $\left | \int_{|x| > K_0} (n_p |u_p|^2) (T_0) \right | \leqslant \int_{|x| > K_0} |n_pu_pU_p| (T_0) + \int_{|x| > K_0} |n_pu_pR^u| (T_0)$. 
\begin{itemize}
	\item First, we can replace $K_0$ by a larger value (depending on $\eta$ and $T_0$) such that $|R^u (T_0 \, , x)| \leqslant \frac{1}{4}$ for all $x \in \R$ such that $|x|>K_0$. Then,
	\[ \int_{|x|>K_0} |n_pu_pR^u| (T_0) \leqslant \frac{1}{4} \int_{|x|>K_0} \frac{n_p^2 + |u_p|^2}{2} (T_0) \leqslant \frac{1}{8} \int_{|x|>K_0} n_p(T_0)^2 + \frac{\eta}{32}. \]
	\item Second, thanks to Sobolev embedding, $||U_p (T_0)||_{L^\infty} \leqslant C ||U_p (T_0)||_{H^1} \leqslant C e^{- \theta_0 T_0} \leqslant \frac{1}{4}$ provided that $T_0$ is chosen large enough. Then
	\[ \int_{|x|>K_0} |n_pu_pU_p| (T_0) \leqslant \frac{1}{8} \int_{|x|>K_0} n_p(T_0)^2 + \frac{\eta}{32} \]
	as well.
\end{itemize}
\noindent Hence, $\left | \int_{|x|>K_0} (n_p|u_p|^2) (T_0) \right | \leqslant \frac{1}{4} \int_{|x|>K_0} n_p (T_0)^2 + \frac{\eta}{16}$. It follows that
\[ \frac{\eta}{8} \geqslant \int_{|x| > K_0} \left ( | \partial_x u_p |^2 + n_p |u_p|^2 + \frac{n_p^2 + v_p^2}{2} \right ) (T_0) \geqslant \int_{|x| > K_0} \left ( | \partial_x u_p |^2 - \frac{n_p^2}{4} + \frac{n_p^2 + v_p^2}{2} \right ) (T_0) - \frac{\eta}{16} \]
and thus
\[ \int_{|x|>K_0} \left ( | \partial_x u_p |^2 + \frac{n_p^2}{4} + \frac{v_p^2}{2} \right ) (T_0) \leqslant \frac{\eta}{16}. \]
This concludes the proof. \hfill \qedsymbol

\noindent \\ Now we can finally properly extract an $(H^1 \times L^2 \times L^2)$-converging subsequence for the initial data. We will then conclude the proof of Theorem 1 by using the wellposedness of the Cauchy problem for \eqref{ZL} in $H^1 \times L^2 \times L^2$: the convergence of the initial data leads to the convergence of the solution.

\begin{leftbar}
\noindent \textbf{Lemma 15.} There exist $u_\infty^0 \in H^2$, $n_\infty^0 \in H^1$, $v_\infty^0 \in H^1$ and a subsequence $((u_{\psi (p)} \, , n_{\psi (p)} \, , v_{\psi (p)}))$ of $((u_p \, , n_p \, , v_p))$ such that
\[ u_{\psi (p)} (T_0) \, \underset{p \to + \infty}{\overset{H^1}{\longrightarrow}} \, u_\infty^0, \midspace n_{\psi (p)} (T_0) \, \underset{p \to + \infty}{\overset{L^2}{\longrightarrow}} \, n_\infty^0 \midspace \text{and} \midspace v_{\psi (p)} (T_0) \, \underset{p \to + \infty}{\overset{L^2}{\longrightarrow}} \, v_\infty^0. \]
\end{leftbar}

\noindent \textit{Proof.} Since $||(u_p \, , n_p \, , v_p) (T_0)||_{H^2 \times H^1 \times H^1} \leqslant C$ (see Corollary 1), Rellich's compactness theorem ensures that there exist $(u_\infty^0 \, , n_\infty^0 \, , v_\infty^0 ) \in H^2 \times H^1 \times H^1$ such that
\[ u_{\psi (p)} (T_0) \, \underset{p \to + \infty}{\overset{H_{\text{loc}}^1}{\longrightarrow}} \, u_\infty^0, \midspace n_{\psi (p)} (T_0) \, \underset{p \to + \infty}{\overset{L_{\text{loc}}^2}{\longrightarrow}} \, n_\infty^0 \midspace \text{and} \midspace v_{\psi (p)} (T_0) \, \underset{p \to + \infty}{\overset{L_{\text{loc}}^{2}}{\longrightarrow}} \, v_\infty^0. \]
Take $\eta > 0$ and the associated $K_0$ in Lemma 14. Possibly replacing $K_0$ by a larger value, we can also assume that $\int_{|x|>K_0} \left ( |u_\infty^0|^2 + |\partial_x u_\infty^0|^2 \right ) \leqslant \eta$. Moreover, by the convergence above, $\int_{|x| \leqslant K_0} \left ( |u_{\psi (p)}|^2 + | \partial_x u_{\psi (p)} |^2 \right ) (T_0) \leqslant \eta$ for $p$ large enough. Combining these two bounds with Lemma 14, it follows that
\[ || u_{\psi (p)} (T_0) - u_\infty^0 ||_{H^1}^2 \leqslant 3 \eta \]
for $p$ large enough. This proves that $u_{\psi (p)} (T_0) \, \underset{p \to + \infty}{\longrightarrow} \, u_\infty^0$ in $H^1 ( \R )$. \\
\\ We proceed similarly to obtain the convergences $n_{\psi (p)} (T_0) \, \underset{p \to + \infty}{\longrightarrow} \, n_\infty^0$ in $L^2 ( \R )$ and $v_{\psi (p)} (T_0) \, \underset{p \to + \infty}{\longrightarrow} \, v_\infty^0$ in $L^2 ( \R )$.
\hfill \qedsymbol

\noindent \\ \textbf{It is time to prove Theorem 1.} Consider the global $H^1 \times L^2 \times L^2$ solution $(u_\infty \, , n_\infty \, , v_\infty )$ of
\[ \left \{ \begin{array}{l} \partial_t u_\infty = i \partial_x^2 u_\infty - in_\infty u_\infty \\ \partial_t n_\infty = - \partial_x v_\infty \\ \partial_t v_\infty = - \partial_x^2 n_\infty - \partial_x^2 (|u_\infty|^2 ) \\ u_\infty (T_0) = u_\infty^0 , \, \, n_\infty (T_0) = n_\infty^0 , \, \, v_\infty (T_0) = v_\infty^0. \end{array} \right. \]
Take $t \geqslant T_0$ and $p \geqslant 1$ large enough such that $T_p > t$. The Cauchy problem for the Zakharov system is well-posed for in $(u \, , n \, , v ) \in H^1 \times L^2 \times L^2$ (see \cite{Gin}). By the continuous dependence of the solution with regards to the initial data, we have
\[ u_{\psi (p)} (t) \, \underset{p \to + \infty}{\overset{H^1}{\longrightarrow}} \, u_\infty (t) , \midspace n_{\psi (p)} (t) \, \underset{p \to + \infty}{\overset{L^2}{\longrightarrow}} \, n_\infty (t) \midspace \text{and} \midspace v_{\psi (p)} (t) \, \underset{p \to + \infty}{\overset{L^2}{\longrightarrow}} \, v_\infty (t). \]
Moreover, recall from Corollary 1 that the sequence $((u_p (t) \, , n_p (t) \, , v_p (t)))_{p \geqslant 1}$ is uniformly bounded in $H^2 \times H^1 \times H^1$. As a consequence, possibly replacing $\psi (p)$ by another subsequence, we have
\[ u_{\psi (p)} (t) \, \underset{p \to + \infty}{\overset{H^2}{\rightharpoonup}} \, u_\infty (t) , \midspace n_{\psi (p)} (t) \, \underset{p \to + \infty}{\overset{H^1}{\rightharpoonup}} \, n_\infty (t) \midspace \text{and} \midspace v_{\psi (p)} (t) \, \underset{p \to + \infty}{\overset{H^1}{\rightharpoonup}} \, v_\infty (t). \] 
Using Proposition 3 and the inferior semi-continuity of the norm, we have
\begin{align*}
& || u_\infty (t) - R^u (t) ||_{H^2} \leqslant \liminf\limits_{p \to + \infty} || u_{\psi (p)} (t) - R^u (t) ||_{H^2} \leqslant C e^{- \theta_0 t} , \\
& || n_\infty (t) - R^n (t) ||_{H^1} \leqslant \liminf\limits_{p \to + \infty} || n_{\psi (p)} (t) - R^n (t) ||_{H^1} \leqslant C e^{- \theta_0 t}  \\
\text{and} \, \, \, & || v_\infty (t) -  R^v (t) ||_{H^1} \leqslant \liminf\limits_{p \to + \infty} || v_{\psi (p)} -  R^v (t) ||_{H^1} \leqslant C e^{- \theta_0 t}. 
\end{align*}
We can replace $C$ by a greater constant so that the inequalities above hold for all $t \geqslant 0$. This concludes the proof of Theorem 1 (in the statement of Theorem 1, $(u_\infty \, , n_\infty \, , v_\infty )$ are simply denoted $(u \, , n \, , v)$). \hfill \qedsymbol

\section*{Appendix}

\noindent \textit{Proof of Proposition 1.} Instead of a proof based on the implicit function theorem, we construct here directly the functions $\omega_k$, $\sigma_k$ and $\gamma_k$ as solutions of a well chosen differential system, then verify that the orthogonality relations are deduced from this system. See \cite{Com,MaMe2} for a similar approach for the gKdV equation. \\
\\ Take a set of parameters $\Pi = (\omega_1 \, , ... \, , \omega_K \, , \sigma_{1} \, , ... \, , \sigma_{K} \, , \gamma_1 \, , ... \, , \gamma_K ) \in (0 \, , + \infty )^K \times \R^{2K}$. For all $k \in [\![ 1 \, , K ]\!]$, define 
\begin{align*}
& x_k^{t,\Pi} := x-c_kt - \sigma_k , \quad \Gamma_k^\Pi (t \, , x) = \frac{c_k   x}{2} - \frac{c_k^2 t}{4} + \omega_k^0 t + \gamma_k, \\
& S_k^u (t \, , \Pi \, , x ) = \sqrt{1-c_k^2} \, \phi_{\omega_k} (x_k^{t,\Pi}) e^{i \Gamma_k^\Pi (t,x)} , \quad S_k^n (t \, , \Pi \, , x) = - \phi_{\omega_k}^2 (x_k^{t,\Pi}) \quad \text{and} \quad S_k^v (t \, , \Pi \, , x) = -c_k \phi_{\omega_k}^2 (x_k^{t,\Pi}).
\end{align*}
Define 
\begin{align*} 
& S^u = \sum\limits_{k=1}^K S_k^u (t \, , \Pi \, , x) , \quad S^n = \sum\limits_{k=1}^K S_k^n (t \, , \Pi \, , x) \quad \text{and} \quad S^v = \sum\limits_{k=1}^K S_k^v (t \, , \Pi \, , x); \\
& \varepsilon_u (t \, , \Pi ) = u(t) - S^u (t \, , \Pi) , \quad \varepsilon_n (t \, , \Pi) = n(t) - S^n (t \, , \Pi ) \quad \text{and} \quad \varepsilon_v (t) = v(t) - S^v (t \, , \Pi ).
\end{align*} 
Set $\Pi^0 = (\omega_1^0 \, , ... \, , \omega_K^0 \, , \sigma_1^0 \, , ... \, , \sigma_K^0 \, , \gamma_1^0 \, , ... \, , \gamma_K^0 )$. Note that $(S^u \, , S^n \, , S^v)(t \, , \Pi^0 ) = (R^u \, , R^n \, , R^v)(t)$. Also note that there exists $C_1 > 0$ such that, for all $t \in [t^* \, , T_p ]$ and all $\Pi \in ( 0 \, , + \infty )^K \times \R^{2K}$,
\begin{equation}
||(S^u \, , S^n \, , S^v )(t \, , \Pi ) - (R^u \, , R^n \, , R^v) (t) ||_{\mathbf{H}} \leqslant C_1 | \Pi - \Pi^0 |. 
\label{RS}
\end{equation}
We prove Proposition 1 in three steps: first, we introduce, out of a hat, a suitable differential system which we prove to have a solution on a certain interval $[T_* \, , T_p]$. Then we show that this differential system is equivalent to the orthogonality relations we desire. Last, we check that $T_* = t^*$. \\
\\ \textit{(Differential system.)} In the following lines we define a lot of quantities. The reason for their introduction is contained in system \eqref{modulEDO} below. See the remark just after \eqref{modulEDO}. The choice made here is to write properly, prove the wellposedness of and solve the differentiel system \eqref{modulEDO}, and then to explain where this system comes from and deduce from it the orthogonality relations we seek to establish. \\
\\ Define
\[ \Sigma^0 (t \, , \Pi \, , x) = -i S^n S^u + \sum\limits_{j=1}^K \sqrt{1-c_j^2} \left ( -i \phi_{\omega_j}^3 (x_j^{t,\Pi}) \right ) e^{i \Gamma_j^\Pi} = i \sum\limits_{1 \leqslant j \neq \ell \leqslant K} \sqrt{1-c_j^2} \, \phi_{\omega_j} (x_j^{t,\Pi}) \phi_{\omega_\ell}^2 (x_\ell^{t,\Pi}) e^{i \Gamma_j^\Pi}. \]
Also define 
\begin{align*} 
& \Sigma_{j}^\omega (t \, , \Pi) = - \sqrt{1-c_j^2} \, \Lambda_{\omega_j} (x_j^{t,\Pi}) e^{i \Gamma_j^\Pi}, \\
& \Sigma_{j}^\sigma (t \, , \Pi ) = \sqrt{1-c_j^2} \, \phi_{\omega_j} ' (x_j^{t,\Pi}) e^{i \Gamma_j^\Pi} \\
\text{and} \, \, \, & \Sigma_{j}^\gamma (t \, , \Pi) = - \sqrt{1-c_j^2} \, i \phi_{\omega_j} (x_j^{t,\Pi}) e^{i \Gamma_j^\Pi}.
\end{align*}
Also define
\[ \Sigma_k^\varepsilon (t \, , \Pi ) = i \varepsilon_u \sum\limits_{\substack{j=1 \\ j \neq k}}^K \phi_{\omega_j}^2 (x_j^{t,\Pi}) - i \varepsilon_n \sum\limits_{\substack{j=1 \\ j \neq k}}^K \sqrt{1-c_j^2} \, \phi_{\omega_j} (x_j^{t,\Pi}) e^{-i \Gamma_j^\Pi}. \]
Set
\begin{align*} 
& \Upsilon_j^{\omega,0} (t \, , \Pi ) = \text{Re} \int_{\R} \phi_{\omega_k} (x_k^{t,\Pi}) e^{-i \Gamma_k^\Pi} \Sigma^0, & & \Upsilon_{k,j}^{\omega,\omega} (t \, , \Pi) = \text{Re} \int_{\R}  \phi_{\omega_k} (x_k^{t,\Pi}) e^{-i \Gamma_k^\Pi} \Sigma_{j}^\omega \\ 
& \Upsilon_{k,j}^{\omega,\sigma} (t \, , \Pi ) = \text{Re} \int_{\R}  \phi_{\omega_k} (x_k^{t,\Pi}) e^{-i \Gamma_k^\Pi} \Sigma_{j}^\sigma , & & \Upsilon_{k,j}^{\omega,\gamma} (t \, , \Pi) = \text{Re} \int_{\R}  \phi_{\omega_k} (x_k^{t,\Pi}) e^{-i \Gamma_k^\Pi} \Sigma_{j}^\gamma , \\
& \Upsilon_k^{\omega , \varepsilon} (t \, , \Pi) = \text{Re} \int_{\R}  \phi_{\omega_k} (x_k^{t,\Pi}) e^{-i \Gamma_k^\Pi} \Sigma_k^\varepsilon , \\
& \Upsilon_j^{\sigma,0} (t \, , \Pi ) = \text{Re} \int_{\R} x_k^{t,\Pi} \phi_{\omega_k} (x_k^{t,\Pi}) e^{-i \Gamma_k^\Pi} \Sigma^0, & & \Upsilon_{k,j}^{\sigma,\omega} (t \, , \Pi) = \text{Re} \int_{\R} x_k^{t,\Pi} \phi_{\omega_k} (x_k^{t,\Pi}) e^{-i \Gamma_k^\Pi} \Sigma_{j}^\omega \\ 
& \Upsilon_{k,j}^{\sigma,\sigma} (t \, , \Pi ) = \text{Re} \int_{\R} x_k^{t,\Pi} \phi_{\omega_k} (x_k^{t,\Pi}) e^{-i \Gamma_k^\Pi} \Sigma_{j}^\sigma , & & \Upsilon_{k,j}^{\sigma,\gamma} (t \, , \Pi) = \text{Re} \int_{\R} x_k^{t,\Pi} \phi_{\omega_k} (x_k^{t,\Pi}) e^{-i \Gamma_k^\Pi} \Sigma_{j}^\gamma , \\
& \Upsilon_k^{\sigma , \varepsilon} (t \, , \Pi) = \text{Re} \int_{\R} x_k^{t,\Pi} \phi_{\omega_k} (x_k^{t,\Pi}) e^{-i \Gamma_k^\Pi} \Sigma_k^\varepsilon , \\
& \Upsilon_j^{\gamma,0} (t \, , \Pi ) = \text{Im} \int_{\R} \Lambda_{\omega_k} (x_k^{t,\Pi}) e^{-i \Gamma_k^\Pi} \Sigma^0, & & \Upsilon_{k,j}^{\gamma,\omega} (t \, , \Pi) = \text{Im} \int_{\R} \Lambda_{\omega_k} (x_k^{t,\Pi}) e^{-i \Gamma_k^\Pi} \Sigma_{j}^\omega \\ 
& \Upsilon_{k,j}^{\gamma,\sigma} (t \, , \Pi ) = \text{Im} \int_{\R} \Lambda_{\omega_k} (x_k^{t,\Pi}) e^{-i \Gamma_k^\Pi} \Sigma_{j}^\sigma , & & \Upsilon_{k,j}^{\gamma,\gamma} (t \, , \Pi) = \text{Im} \int_{\R} \Lambda_{\omega_k} (x_k^{t,\Pi}) e^{-i \Gamma_k^\Pi} \Sigma_{j}^\gamma , \\
& \Upsilon_k^{\gamma , \varepsilon} (t \, , \Pi) = \text{Im} \int_{\R} \Lambda_{\omega_k} (x_k^{t,\Pi}) e^{-i \Gamma_k^\Pi} \Sigma_k^\varepsilon .
\end{align*}
Define the three constants $\mathbf{d}^\omega = \int_{\R} Q \Lambda Q = 1 > 0$, $\mathbf{d}^\sigma = -\int_{\R} yQQ' = \frac{1}{2} \int_{\R} Q^2 = 2 > 0$ and $\mathbf{d}^\gamma = \int_{\R} Q \Lambda Q = 1 > 0$. Now define the following functions of $t$ and $\Pi$: 
\begin{align*}
& \Theta_{k,k}^{\omega , \omega}  = - \frac{\sqrt{\omega_k}}{\sqrt{1-c_k^2}} \Ree \int_{\R} \Lambda_{\omega_k} (x_k^{t,\Pi}) e^{-i \Gamma_k^\Pi} \varepsilon_u , & & \Theta_{k,k}^{\omega , \sigma}   = \frac{\sqrt{\omega_k}}{\sqrt{1-c_k^2}} \Ree \int_{\R}  \phi_{\omega_k} ' (x_k^{t,\Pi}) e^{-i \Gamma_k^\Pi} \varepsilon_u , \\
& \Theta_{k,k}^{\omega , \gamma}   = - \frac{\sqrt{\omega_k}}{\sqrt{1-c_k^2}} \Imm \int_{\R} \phi_{\omega_k} (x_k^{t,\Pi}) e^{-i \Gamma_k^\Pi} \varepsilon_u , & & \Theta_{k,j}^{\omega , \omega}   = - \frac{\sqrt{\omega_k}}{\sqrt{1-c_k^2}} \Upsilon_{k,j}^{\omega , \omega} \, \, \text{for} \, \, j \neq k , \\
& \Theta_{k,j}^{\omega , \sigma}   = - \frac{\sqrt{\omega_k}}{\sqrt{1-c_k^2}} \Upsilon_{k,j}^{\omega , \sigma} \, \, \text{for} \, \, j \neq k , & & \Theta_{k,j}^{\omega , \gamma}   = - \frac{\sqrt{\omega_k}}{\sqrt{1-c_k^2}} \Upsilon_{k,j}^{\omega , \gamma} \, \, \text{for} \, \, j \neq k , \\
& \Theta_k^{\omega,q}   = \frac{\sqrt{\omega_k}}{\sqrt{1-c_k^2}} \Imm \int_{\R} \phi_{\omega_k} (x_k^{t,\Pi}) e^{-i \Gamma_k^\Pi} \varepsilon_u \varepsilon_n , & & \Theta_k^{\omega , \Upsilon}   = \frac{\sqrt{\omega_k}}{\sqrt{1-c_k^2}} \left (  \Upsilon_k^{\omega , 0} + \Upsilon_k^{\omega , \varepsilon} \right ). 
\end{align*}
Also, set
\begin{align*}
& \Theta_{k,k}^{\sigma , \omega}  = - \frac{1}{\sqrt{\omega_k} \sqrt{1-c_k^2}} \Ree \int_{\R} x_k^{t,\Pi} \Lambda_{\omega_k} (x_k^{t,\Pi}) e^{-i \Gamma_k^\Pi} \varepsilon_u , \\
& \Theta_{k,k}^{\sigma , \sigma} = \frac{1}{\sqrt{\omega_k} \sqrt{1-c_k^2}} \Ree \int_{\R} \left ( x_k^{t,\Pi} \phi_{\omega_k} ' (x_k^{t,\Pi})  + \phi_{\omega_k} (x_k^{t,\Pi}) \right ) e^{-i \Gamma_k^\Pi} \varepsilon_u , \\
& \Theta_{k,k}^{\sigma , \gamma} = - \frac{1}{\sqrt{\omega_k} \sqrt{1-c_k^2}} \Imm \int_{\R} x_k^{t , \Pi} \phi_{\omega_k} (x_k^{t,\Pi}) e^{-i \Gamma_k^\Pi} \varepsilon_u , & & \Theta_{k,j}^{\sigma , \omega} = - \frac{1}{\sqrt{\omega_k} \sqrt{1-c_k^2}} \Upsilon_{k,j}^{\sigma , \omega} \, \, \text{for} \, \, j \neq k , \\
& \Theta_{k,j}^{\sigma , \sigma} = - \frac{1}{\sqrt{\omega_k} \sqrt{1-c_k^2}} \Upsilon_{k,j}^{\sigma , \sigma} \, \, \text{for} \, \, j \neq k , & & \Theta_{k,j}^{\sigma , \gamma} = - \frac{1}{\sqrt{\omega_k} \sqrt{1-c_k^2}} \Upsilon_{k,j}^{\sigma , \gamma} \, \, \text{for} \, \, j \neq k, \\
& \Theta_k^{\sigma , q} = - \frac{2}{\sqrt{\omega_k} \sqrt{1-c_k^2}} \Imm \int_{\R} \phi_{\omega_k} ' (x_k^{t,\Pi}) e^{-i \Gamma_k^\Pi} \varepsilon_u  & & \Theta_k^{\sigma , \Upsilon} = \frac{1}{\sqrt{\omega_k} \sqrt{1-c_k^2}} \left ( \Upsilon_k^{\sigma , 0} + \Gamma_k^{\sigma , \varepsilon} \right ). \\
& \midspace \midspace \midspace + \frac{1}{\sqrt{\omega_k} \sqrt{1-c_k^2}} \Imm \int_{\R} x_k^{t,\Pi} \phi_{\omega_k} (x_k^{t,\Pi}) e^{-i \Gamma_k^\Pi} \varepsilon_u \varepsilon_n , 
\end{align*}
Also, set
\begin{align*}
& \Theta_{k,k}^{\gamma , \omega} = \frac{\sqrt{\omega_k}}{\sqrt{1-c_k^2}} \Imm \int_{\R} \left. \frac{\partial \Lambda_\omega}{\partial \omega} \right |_{\omega = \omega_k} (x_k^{t,\Pi}) e^{-i \Gamma_k} \varepsilon_u , & & \Theta_{k,k}^{\gamma , \sigma} = - \frac{\sqrt{\omega_k}}{\sqrt{1-c_k^2}} \Imm \int_{\R} \Lambda_{\omega_k} ' (x_k^{t,\Pi}) e^{-i \Gamma_k} \varepsilon_u , \\
& \Theta_{k,k}^{\gamma , \gamma} = - \frac{\sqrt{\omega_k}}{\sqrt{1-c_k^2}} \Ree \int_{\R} \Lambda_{\omega_k} (x_k^{t,\Pi}) e^{-i \Gamma_k} \varepsilon_u , & & \Theta_{k,j}^{\gamma , \omega} = \frac{\sqrt{\omega_k}}{\sqrt{1 - c_k^2}} \Upsilon_{k,j}^{\gamma , \omega} \, \, \text{for} \, \, j \neq k , \\
& \Theta_{k,j}^{\gamma , \sigma} = \frac{\sqrt{\omega_k}}{\sqrt{1-c_k^2}} \Upsilon_{k,j}^{\gamma , \sigma} \, \, \text{for} \, \, j \neq k , & & \Theta_{k,j}^{\gamma , \gamma} = \frac{\sqrt{\omega_k}}{\sqrt{1-c_k^2}} \Upsilon_{k,j}^{\gamma , \gamma} \, \, \text{for} \, \, j \neq k , \\
& \Theta_k^{\gamma , q} = \frac{\sqrt{\omega_k}}{\sqrt{1-c_k^2}} \Ree \int_{\R} \left ( - \phi_{\omega_k}  + 2 \phi_{\omega_k}^2 \Lambda_{\omega_k} \right ) (x_k^{t , \Pi}) e^{-i \Gamma_k^\Pi} \varepsilon_u & & \Theta_k^{\gamma , \Upsilon} = \frac{\sqrt{\omega_k}}{\sqrt{1-c_k^2}} \left ( \Upsilon_k^{\gamma , 0} + \Upsilon_k^{\gamma , \varepsilon} \right ). \\
& \midspace \midspace \midspace + \sqrt{\omega_k} \, \text{Re} \int_{\R} (\phi_{\omega_k} \Lambda_{\omega_k} )(x_k^{t,\Pi}) e^{-i \Gamma_k} \varepsilon_n \\
& \midspace \midspace \midspace + \frac{\sqrt{\omega_k}}{\sqrt{1-c_k^2}} \Ree \int_{\R} \Lambda_{\omega_k} (x_k^{t,\Pi}) e^{-i \Gamma_k^\Pi} \varepsilon_u \varepsilon_n ,
\end{align*}
For letters $\mathfrak{a}$,$\mathfrak{b} \in \{ \omega \, , \sigma \, , \gamma \}$, set $\Theta^{\mathfrak{a} , \mathfrak{b}} (t \, , \Pi) = \left ( \Theta_{k,j}^{\mathfrak{a} , \mathfrak{b}} (t \, , \Pi) \right )_{1 \leqslant k,j \leqslant K}$. For letters $\mathfrak{a} \in \{ \omega \, , \sigma \, , \gamma \}$ and $\mathfrak{b} \in \{ q \, , \Upsilon \}$, set $\Theta^{\mathfrak{a} , \mathfrak{b}} (t \, , \Pi) = \left ( \Theta_k^{\mathfrak{a} , \mathfrak{b}} (t \, , \Pi) \right )_{1 \leqslant k \leqslant K}$. We consider the following differential system:
\begin{equation}
\left \{ \begin{array}{l} \left ( \mathbf{I}_{3K} + \mathbf{A} (t \, , \Pi) \right ) \left ( \dot{\Pi} (t) - \mathbf{N} (\Pi - \Pi^0) \right ) = \mathbf{B} (t \, , \Pi ) \\ \Pi(T_p) = \Pi^0. \end{array} \right.
\label{modulEDO}
\end{equation}
where 
\[ \mathbf{A} (t \, , \Pi ) = \left [ \begin{array}{ccc} \frac{1}{\mathbf{d}^\omega} \Theta^{\omega , \omega}(t \, , \Pi) & \frac{1}{\mathbf{d}^\omega} \Theta^{\omega , \sigma}(t \, , \Pi) & \frac{1}{\mathbf{d}^\omega} \Theta^{\omega , \gamma}(t \, , \Pi) \\ 
\frac{1}{\mathbf{d}^\sigma} \Theta^{\sigma , \omega}(t \, , \Pi) & \frac{1}{\mathbf{d}^\sigma} \Theta^{\sigma , \sigma}(t \, , \Pi) &  \frac{1}{\mathbf{d}^\sigma} \Theta^{\sigma , \gamma}(t \, , \Pi) \\ 
\frac{1}{\mathbf{d}^\gamma} \Theta^{\gamma , \omega}(t \, , \Pi) & \frac{1}{\mathbf{d}^\gamma} \Theta^{\gamma , \sigma}(t \, , \Pi) &  \frac{1}{\mathbf{d}^\gamma} \Theta^{\gamma , \gamma}(t \, , \Pi) \end{array} \right ], \]
\[ \mathbf{B} (t \, , \Pi) = \left [ \begin{array}{c} \frac{1}{\mathbf{d}^\omega} \left ( \Theta^{\omega,q } (t \, , \Pi) + \Theta^{\omega , \Upsilon} (t \, , \Pi) \right ) \\ 
\frac{1}{\mathbf{d}^\sigma} \left ( \Theta^{\sigma , q } (t \, , \Pi) + \Theta^{\sigma , \Upsilon} (t \, , \Pi) \right ) \\ 
\frac{1}{\mathbf{d}^\gamma} \left ( \Theta^{\gamma , q } (t \, , \Pi) + \Theta^{\gamma , \Upsilon} (t \, , \Pi) \right ) \end{array} \right ] \]
and
\[ \mathbf{N} = \left [ \begin{array}{ccc}  0_{K,K} & 0_{K,K} & 0_{K,K} \\  0_{K,K} & 0_{K,K} & 0_{K,K} \\  I_K & 0_{K,K} & 0_{K,K} \end{array} \right ]. \]
The reason for the introduction of such quantities and such a system is that \eqref{modulEDO} above is equivalent to
\[ \forall k \in \{ 1 \, , ... \, , K \} , \, \, \, \frac{\text{d}}{\text{d}t} \left ( \begin{array}{c} \text{Re} \int_{\R} \phi_{\omega_k (t)} (x-c_kt - \sigma_k (t)) e^{- i \Gamma_k (t,x)} \varepsilon_u (t \, , x) \, \text{d}x \\ \text{Re} \int_{\R} (x-c_kt - \sigma_k (t)) \phi_{\omega_k (t)} (x-c_kt - \sigma_k (t)) e^{-i \Gamma_k (t,x)} \varepsilon_u (t \, , x) \, \text{d}x \\ \text{Im} \int_{\R} \Lambda_{\omega_k (t)} (x-c_kt - \sigma_k (t)) e^{-i \Gamma_k (t,x)} \varepsilon_u (t \, , x) \, \text{d}x \end{array} \right ) = \left ( \begin{array}{c} 0 \\ 0 \\ 0 \end{array} \right ) \]
(assuming that the functions $\omega_k$, $\sigma_k$ and $\gamma_k$ are $\mathscr{C}^1$), which is the set of orthogonalities we seek to establish. In this proof, we begin to verify that the system \eqref{modulEDO} is well-posed and then, considering its solution, we will establish the equivalence between \eqref{modulEDO} and the orthogonality relations we desire; see equations \eqref{explomega}, \eqref{explsigma} and \eqref{explgamma}. \\
\\ It is clear that $\mathbf{A}$ and $\mathbf{B}$ are continuous with regards to $t \in [t^* \, , T_p]$ and $\mathscr{C}^1$ (even $\mathscr{C}^\infty$) with regards to $\Pi \in (0 \, , + \infty )^K \times \R^{2K}$. Let us introduce the following hypothesis \eqref{modulHYP}:
\begin{equation}
|| ( \varepsilon_u \, , \varepsilon_n \, , \varepsilon_v ) (t \, , \Pi ) ||_{\mathbf{H}} + | \Pi - \Pi^0 | \leqslant C_0 A_0 e^{- \theta_0 t}
\label{modulHYP}
\end{equation}
where $C_0 >0$ is constant whose expression shall be introduced later on and which does not depend on $t$, $A_0$ or $T_0$. See the very end of the proof for the suitable choice of $C_0$. \\
\\ We can take $T_0>0$ large enough (depending on $A_0$) such that, under hypothesis \eqref{modulHYP}, $\omega_- \leqslant \omega_k \leqslant \omega^+$ for all $k \in \{ 1 \, , ... \, , K \}$. Under hypothesis \eqref{modulHYP}, the matrix $\mathbf{A} (t \, , \Pi )$ is invertible. Indeed, first, the terms $\Upsilon$ can all be estimated exponentially. Let $f$,$g$ be two functions that satisfies $|f(x)| + |g(x)| \leqslant Ce^{- \frac{\sqrt{\omega_-} |x|}{4}}$ and let us estimate $\int_{\R} f(x_k^{t,\Pi}) g(x_j^{t,\Pi}) \, \text{d}x$ for $j \neq k$. We have
\begin{align*}
\left | \int_{\R} f(x_k^{t,\Pi}) g(x_j^{t,\Pi}) \, \text{d}x \right | \leqslant & \, \,  C \int_{\R} e^{- \frac{\sqrt{\omega_-} |x-c_kt - \sigma_k|}{4}} e^{- \frac{\sqrt{\omega_-} |x-c_jt - \sigma_j|}{4}} \, \text{d}x \\
\leqslant & \, \, C \int_{|x-c_jt - \sigma_j| \leqslant \frac{|(c_k-c_j)t + \sigma_k - \sigma_j|}{2}} e^{- \frac{\sqrt{\omega_-} |x-c_kt - \sigma_k |}{4}} e^{- \frac{\sqrt{\omega_-} |x-c_jt - \sigma_j|}{4}} \, \text{d}x \\
& \midspace + C \int_{|x-c_jt - \sigma_j| > \frac{|(c_k-c_j)t + \sigma_k - \sigma_j|}{2}} e^{- \frac{\sqrt{\omega_-} |x-c_kt - \sigma_k|}{4}} e^{- \frac{\sqrt{\omega_-} |x-c_jt - \sigma_j|}{4}} \, \text{d}x  \\
\leqslant & \, \, C e^{- \frac{\sqrt{\omega_-}}{8} |(c_k-c_j)t + \sigma_k  - \sigma_j |} \int_{\R} e^{- \frac{\sqrt{\omega_-} |x-c_jt - \sigma_j |}{4}} \, \text{d}x \\
& \midspace + C e^{- \frac{\sqrt{\omega_-}}{8} |(c_k-c_j)t + \sigma_k  - \sigma_j |} \int_{\R} e^{- \frac{\sqrt{\omega_-} |x-c_kt - \sigma_k |}{4}} \, \text{d}x \\
\leqslant & \, \, C e^{- \frac{\sqrt{\omega_-}}{8} |(c_k-c_j)t + \sigma_k - \sigma_j |}. 
\end{align*}
Since $| \sigma_k - \sigma_k^0 | \leqslant 2 A_0 e^{- \theta_0 t} \leqslant C$, $| \sigma_k - \sigma_j | \leqslant |\sigma_k^0 - \sigma_j^0| + C \leqslant C$. By definition of $\theta_0$, $|(c_k-c_j)t + \sigma_k  - \sigma_j| \geqslant |c_k-c_j|t - C \geqslant 8 \sqrt{\theta_0} \, t$ for $T_0>0$ large enough. It follows that
\[ \left | \int_{\R} f(x_k^{t,\Pi}) g(x_j^{t,\Pi}) \, \text{d}x \right | \leqslant C e^{- \sqrt{\omega_-} \sqrt{\theta_0} t} \leqslant C e^{- 16 \theta_0 t}. \]
Moreover, under hypothesis \eqref{modulHYP}, $|| ( \varepsilon_u \, , \varepsilon_n \, , \varepsilon_v) (t \, , \Pi ) ||_{\mathbf{H}} \leqslant C$. It follows that $| \Upsilon_{k,j}^{\mathfrak{a} , \mathfrak{b}} (t \, , \Pi ) | \leqslant C e^{-16 \theta_0 t}$ for all $j \neq k$ and $\mathfrak{a}$,$\mathfrak{b} \in \{ \omega \, , \sigma \, , \gamma \}$ and $|\Upsilon_k^{\mathfrak{a} , \mathfrak{b}} (t \, , \Pi ) | \leqslant C e^{-16 \theta_0 t}$ for all $\mathfrak{a} \in \{ \omega \, , \sigma \, , \gamma \}$ and $\mathfrak{b} \in \{ 0 \, , \varepsilon \}$. Hence $| \Theta_{k,j}^{\mathfrak{a} , \mathfrak{b}} (t \, , \Pi ) | \leqslant C e^{-16 \theta_0 t}$ for all $j \neq k$ and $\mathfrak{a}$,$\mathfrak{b} \in \{ \omega \, , \sigma \, , \gamma \}$ and $|\Theta_k^{\mathfrak{a} , \Upsilon} (t \, , \Pi ) | \leqslant C e^{-16 \theta_0 t}$ for all $\mathfrak{a} \in \{ \omega \, , \sigma \, , \gamma \}$. On the other hand, $| \Theta_{k,k}^{\mathfrak{a} , \mathfrak{b}} (t \, , \Pi ) | + | \Theta_k^{\mathfrak{a} , q} (t \, , \Pi ) | \leqslant C || ( \varepsilon_u \, , \varepsilon_n \, , \varepsilon_v) (t \, , \Pi ) ||_{\mathbf{H}}$ by simple Cauchy-Schwarz inequalities, for all $k \in \{ 1 \, , ... \, , K \}$ and $\mathfrak{a}$,$\mathfrak{b} \in \{ \omega \, , \sigma \, , \gamma \}$. It follows that, under hypothesis \eqref{modulHYP},
\begin{equation}
| \mathbf{A} (t \, , \Pi ) | + | \mathbf{B} (t \, , \Pi )|  \leqslant C || ( \varepsilon_u \, , \varepsilon_n \, , \varepsilon_v ) (t \, , \Pi ) ||_{\mathbf{H}} + C e^{-16 \theta_0 t} \leqslant C A_0 e^{- \theta_0 t}. 
\label{modulAB}
\end{equation}
\noindent Taking $T_0>0$ large enough (depending on $A_0$), the norm of $\mathbf{A}$ can be taken as small as one wants, provided that hypothesis \eqref{modulHYP} holds. \\
\\ Define the compact $\mathcal{K} = [t^* \, , T_p ] \times \overline{B} \left ( \Pi^0 \, , \frac{1}{4} \right )$. The function $\mathbf{A}$ is uniformly continuous on $\mathcal{K}$: there exists a constant $\delta_A >0$ such that $|\mathbf{A}(t_1 \, , \Pi_1) - \mathbf{A} (t_2 \, , \Pi_2 )| \leqslant \frac{1}{4}$ if $|t_1 - t_2| \leqslant \delta_A$, $|\Pi_1 - \Pi^0| \leqslant \frac{1}{4}$, $|\Pi_2 - \Pi^0| \leqslant \frac{1}{4}$ and $|\Pi_1 - \Pi_2 | \leqslant \delta_A$. We can assume that $\delta_A < \frac{1}{4}$ and $\delta_A < \omega_-$. Set
\[ T_* = \inf \left \{ t_* \in [t^* \, , T_p ] \, \, | \, \, \exists \Pi \in \mathscr{C}^1 ([t_* \, , T_p] \, , \R^{3K})  \, \, \text{solution of} \, \, \eqref{modulEDO} \, \, \text{that satisfies \eqref{modulHYP} for all} \, \, t \in [t_* \, , T_p] \right \} . \]
First, let us show that $T_* < T_p$. Note that $(\varepsilon_u \, , \varepsilon_n \, , \varepsilon_v) (T_p \, , \Pi^0) = 0$ thus \eqref{modulHYP} holds for $t=T_p$ and $\Pi = \Pi^0$. Consequently $|\mathbf{A} (T_p \, , \Pi^0)| \leqslant C || ( \varepsilon_u \, , \varepsilon_n \, , \varepsilon_v ) (T_p \, , \Pi^0) ||_{\mathbf{H}} + C e^{-16 \theta_0 T_p} \leqslant C e^{-16 \theta_0 t} \leqslant \frac{1}{4}$ for $T_0>0$ chosen large enough. For $t \in [T_p - \delta_A \, , T_p]$ and $\Pi \in \overline{B} \left ( \Pi^0 \, , \delta_A \right )$, it follows that $| \mathbf{A} (t \, , \Pi ) | \leqslant \frac{1}{2} < 1$ and the matrix $\mathbf{I}_{3K} + \mathbf{A} (t \, , \Pi )$ is invertible. Hence, there exists $\Pi \in \mathscr{C}^1 \left ( [T_p - \delta_A \, , T_p] \, , \overset{o}{B} \left ( \Pi^0 \, , \delta_A \right ) \right )$ solution of \eqref{modulEDO}. Since the function $t \mapsto || \varepsilon (t \, , \Pi (t))||_{\mathbf{H}} + | \Pi (t) - \Pi^0 |$ is continuous on $[T_p - \delta_A \, , T_p]$ and vanishes at $t=T_p$, we can take $\delta_A>0$ small enough (depending on $A_0$ and $T_p$) such that, for all $t \in [T_p - \delta_A \, , T_p]$, \eqref{modulHYP} holds. This demonstrates that $T_* \leqslant T_p - \delta_A < T_p$. We shall prove later on that $T_* = t^*$ (see the end of the proof). Note that, since \eqref{modulHYP} holds on $]T_* \, , T_p]$, the solution $\Pi (t)$ cannot blow up at $T_*^+$ and, by continuity, \eqref{modulHYP} still holds for $t=T_*^+$. \\
\\ \textit{(Orthogonality relations.)} We have proved the existence of a unique solution $\Pi \in \mathscr{C}^1 ([T_* \, , T_p] \, , \R^{3K})$ of the system \eqref{modulEDO}. This solution satisfies \eqref{modulHYP} and \eqref{modulAB} for all $t \in [T_* \, , T_p]$. We simply denote $\Gamma_k (t) = \Gamma_k^{\Pi(t)} (t)$, $S_k^* (t) = S_k^* (t \, , \Pi (t))$, $\varepsilon (t) = \varepsilon (t \, , \Pi (t))$, $\Sigma_*^* (t) = \Sigma_*^* (t \, , \Pi (t))$, $\Upsilon_*^* (t) = \Upsilon_*^* (t \, , \Pi (t))$, $\Theta_*^* (t) = \Theta_*^* (t \, , \Pi (t))$, $\mathbf{A} (t) = \mathbf{A} (t \, , \Pi (t))$, $\mathbf{B} (t) = \mathbf{B} (t \, , \Pi (t))$ and $x_k^t = x_k^{t,\Pi(t)}$. \\
\\ The equation of $\varepsilon_u$ can be written as
\begin{equation}
\partial_t \varepsilon_u = i \partial_x^2 \varepsilon_u - iS^u \varepsilon_n - i S^n \varepsilon_u - i \varepsilon_u \varepsilon_n + \Sigma (t)
\label{eqeps}
\end{equation}
where
\begin{align*}
\Sigma (t) =& \, -i S^n S^u + \sum\limits_{j=1}^K \sqrt{1- c_j^2} \left [ - \dot{\omega}_j (t) \Lambda_{\omega_j (t)} (x_j^t)  - \dot{\sigma}_j (t)   \phi_{\omega_j (t)} ' (x_j^t) \right. \\
& \midspace \midspace \midspace \midspace \midspace \midspace \midspace \midspace \midspace \left. + i \left ( \omega_j (t) - \omega_j^0 - \dot{\gamma}_j (t) \right ) \phi_{\omega_j (t)} (x_j^t) - i \phi_{\omega_j (t)}^3 (x_j^t) \right ] e^{i \Gamma_j (t \, , x)} .  
\end{align*}
Let us rewrite $\Sigma$ as
\begin{equation}
\begin{split}
\Sigma =& \, \, \sqrt{1-c_k^2} \left [ - \dot{\omega}_k (t) \Lambda_{\omega_k (t)} (x_k^t) - \dot{\sigma}_k (t)  \phi_{\omega_k (t)} '  (x_k^t) + i \left ( \omega_k (t) - \omega_k^0 - \dot{\gamma}_k (t) \right ) \phi_{\omega_k (t)} (x_k^t) \right ] e^{i \Gamma_k (t  , x)} \\
& \midspace + \Sigma^0 (t) +\sum\limits_{\substack{j=1 \\ j \neq k}}^K \left ( \dot{\omega}_j (t) \Sigma_{j}^{\omega} (t) + \dot{\sigma}_j (t)   \Sigma_{j}^{\sigma} (t) +  ( \dot{\gamma}_j(t) - \omega_j(t) + \omega_j^0 ) \Sigma_{j}^{\gamma} (t) \right )
\end{split} 
\label{defSigma}
\end{equation}
where $\Sigma^0$, $\Sigma_j^\omega$, $\Sigma_j^\sigma$ and $\Sigma_j^\gamma$ have been introduced at the beginning of the proof. We also write
\begin{equation}
-iS^u \varepsilon_n - i S^n \varepsilon_u = -i \sqrt{1 - c_k^2} \, \phi_{\omega_k(t)} (x_k^t) e^{i \Gamma_k (t  , x)} \varepsilon_n + i \phi_{\omega_k(t)}^2 (x_k^t) \varepsilon_u + \Sigma_k^\varepsilon
\label{eqSeps}
\end{equation}
where $\Sigma_k^\varepsilon$ has been introduced at the beginning of the proof. \\
\\ First, let us use the equation of $\dot{\omega}_k (t)$ to deduce the first orthogonality. See that
\begin{equation}
\begin{split}
& \, \frac{\text{d}}{\text{d}t} \text{Re} \int_{\R} \phi_{\omega_k (t)} (x-c_kt - \sigma_k (t)) e^{-i \Gamma_k (t  , x)} \varepsilon_u (t \, , x) \, \text{d}x \\
=& \, \, \text{Re} \int_{\R} \phi_{\omega_k} (x_k^t) e^{-i \Gamma_k} ( \partial_t \varepsilon_u ) + \dot{\omega}_k (t) \text{Re} \int_{\R}  \Lambda_{\omega_k (t)} (x_k^t) e^{-i \Gamma_k} \varepsilon_u \\
& \midspace - (c_k + \dot{\sigma}_k (t))   \text{Re} \int_{\R} \phi_{\omega_k} ' (x_k^t) e^{-i \Gamma_k} \varepsilon_u + \left ( \dot{\gamma}_k (t) + \omega_k^0 - \frac{c_k^2}{4} \right ) \text{Im} \int_{\R} \phi_{\omega_k} (x_k^t) e^{-i \Gamma_k} \varepsilon_u .
\end{split} 
\label{modul1}
\end{equation}
\noindent Let us compute $\text{Re} \int_{\R} \phi_{\omega_k} (x_k^t) e^{-i \Gamma_k} ( \partial_t \varepsilon_u )$ using equation \eqref{eqeps}. Let us start with $\Sigma$. First,
\[ \text{Re} \int_{\R} \phi_{\omega_k} (x_k^t) e^{-i \Gamma_k} \Lambda_{\omega_k} (x_k^t) e^{i \Gamma_k} = \int_{\R} \phi_{\omega_k}  \Lambda_{\omega_k} = \frac{1}{\sqrt{\omega_k}} \int_{\R} Q \Lambda Q \]
since $\Lambda_{\omega} (x) = \frac{1}{\sqrt{\omega}} \Lambda Q ( \sqrt{\omega} \, x)$. Second,
\[ \int_{\R} \phi_{\omega_k} (x_k^t) e^{-i \Gamma_k} \phi_{\omega_k} ' (x_k^t) e^{i \Gamma_k} = \int_{\R} \phi_{\omega_k} \phi_{\omega_k} ' = 0. \]
Third,
\[ \text{Re} \int_{\R} \phi_{\omega_k} (x_k^t) e^{-i \Gamma_k} i \phi_{\omega_k} (x_k^t) e^{i \Gamma_k} = - \text{Im} \int_{\R} \phi_{\omega_k}^2 = 0. \]
Thus, we get from \eqref{defSigma} that
\begin{equation}
\begin{split}
\text{Re} \int_{\R} \phi_{\omega_k} (x_k^t) e^{-i \Gamma_k} \Sigma =&  - \dot{\omega}_k (t) \frac{\sqrt{1-c_k^2}}{\sqrt{\omega_k}} \int_{\R} Q \Lambda Q + \Upsilon_k^{\omega,0} (t)  \\
& \midspace + \sum\limits_{\substack{j=1 \\ j \neq k}}^K \left ( \dot{\omega}_j (t) \Upsilon_{k,j}^{\omega,\omega} (t) + \dot{\sigma}_j (t)   \Upsilon_{k,j}^{\omega,\sigma} (t) + \left ( \dot{\gamma}_j (t) - \omega_j (t) + \omega_j^0 \right ) \Upsilon_{k,j}^{\omega,\gamma} (t) \right ).
\end{split}
\label{modul2}
\end{equation}
Now, from \eqref{eqSeps},
\begin{align*} 
& \, \, \text{Re} \int_{\R} \phi_{\omega_k} (x_k^t) e^{-i \Gamma_k} (i \partial_x^2 \varepsilon_u - i S^u \varepsilon_n - i S^n \varepsilon_u ) \\
=&  \, \, \text{Im} \int_{\R} \left ( - \partial_x^2 \varepsilon_u - \phi_{\omega_k}^2 (x_k^t) \varepsilon_u \right ) \phi_{\omega_k} (x_k^t) e^{-i \Gamma_k} + \sqrt{1-c_k^2} \underbrace{\text{Im} \int_{\R}  \phi_{\omega_k}^2 (x_k^t) \varepsilon_n}_{= \, 0} + \Upsilon_k^{\omega, \varepsilon} (t). 
\end{align*}
Defining the operator $L_{-,k,t} = - \partial_x^2 + \omega_k (t) - \phi_{\omega_k (t)}^2 (x-c_kt- \sigma_k (t))$, we have
\begin{align*}
& \, \, \text{Re} \int_{\R} \phi_{\omega_k} (x_k^t) e^{-i \Gamma_k} (i \partial_x^2 \varepsilon_u - i S^u \varepsilon_n - i S^n \varepsilon_u ) \\
=& \, \, \text{Im} \int_{\R} (L_{-,k,t} \varepsilon_u) \phi_{\omega_k} (x_k^t) e^{-i \Gamma_k} - \omega_k(t) \, \text{Im} \int_{\R} \phi_{\omega_k} (x_k^t) e^{-i \Gamma_k} \varepsilon_u + \Upsilon_k^{\omega, \varepsilon} (t).
\end{align*}
Since $L_{-,k,t}$ is self-adjoint, we have $\int_{\R} (L_{-,k,t} \varepsilon_u) \phi_{\omega_k} (x_k^t) e^{-i \Gamma_k} = \int_{\R} \varepsilon_u L_{-,k,t} \left ( \phi_{\omega_k} (x_k^t) e^{-i \Gamma_k} \right )$. Let us compute this integral. Note that
\[  L_{-,k,t} \left ( \phi_{\omega_k} (x_k^t) e^{-i \Gamma_k} \right ) = e^{-i \Gamma_k} L_{-,k,t} \left ( \phi_{\omega_k} (x_k^t) \right ) - 2 \partial_x \left ( \phi_{\omega_k} (x_k^t) \right )   \partial_x \left ( e^{-i \Gamma_k} \right ) - \phi_{\omega_k} (x_k^t) \partial_x^2 \left ( e^{-i \Gamma_k} \right ). \]
We compute $\partial_x \left ( e^{-i \Gamma_k} \right ) = - \frac{i c_k}{2} e^{-i \Gamma_k}$ and $\partial_x^2 \left ( e^{-i \Gamma_k} \right ) = - \frac{c_k^2}{4} e^{-i \Gamma_k}$. Recalling that $L_- Q = 0$ and $\phi_\omega (x) = \sqrt{\omega} \, Q ( \sqrt{\omega} \, x)$, we have $L_{-,k,t} \left ( \phi_{\omega_k} (x_k^t) \right ) = 0$. Gathering these relations, it follows that
\begin{equation}
\begin{split} 
& \text{Re} \int_{\R} \phi_{\omega_k} (x_k^t) e^{-i \Gamma_k} (i \partial_x^2 \varepsilon_u - i S^u \varepsilon_n - i S^n \varepsilon_u ) \\
= \, \,  & c_k \, \text{Re} \int_{\R}  \phi_{\omega_k} ' (x_k^t) e^{-i \Gamma_k} \varepsilon_u + \frac{c_k^2}{4} \Imm \int_{\R} \phi_{\omega_k} (x_k^t) e^{-i \Gamma_k} \varepsilon_u  - \omega_k \text{Im} \int_{\R} \phi_{\omega_k} (x_k^t) e^{-i \Gamma_k} \varepsilon_u + \Upsilon_k^{\omega, \varepsilon }. 
\end{split}
\label{modul3}
\end{equation}
\noindent Combining equations \eqref{eqeps}, \eqref{modul1}, \eqref{modul2} and \eqref{modul3}, we obtain
\begin{align*}
& \frac{\text{d}}{\text{d}t} \Ree \int_{\R} \phi_{\omega_k (t)} (x-c_kt - \sigma_k (t)) e^{-i \Gamma_k (t  , x)} \varepsilon_u (t \, , x) \, \text{d}x \\
=& - \dot{\omega}_k (t) \frac{\sqrt{1-c_k^2}}{\sqrt{ \omega_k}} \int_{\R} Q \Lambda Q + \Upsilon_k^{\omega,0} + \sum\limits_{\substack{j=1 \\ j \neq k}}^K \left ( \dot{\omega}_j (t) \Upsilon_{k,j}^{\omega,\omega} + \dot{\sigma}_j (t)   \Upsilon_{k,j}^{\omega,\sigma} + \left ( \dot{\gamma}_j (t) - \omega_j (t) + \omega_j^0 \right ) \Upsilon_{k,j}^{\omega,\gamma} \right ) \\
& \midspace + \boxed{c_k}   \text{Re} \int_{\R}  \phi_{\omega_k} ' (x_k^t) e^{-i \Gamma_k} \varepsilon_u + \boxed{\frac{c_k^2}{4}} \Imm \int_{\R} \phi_{\omega_k} (x_k^t) e^{-i \Gamma_k} \varepsilon_u - \omega_k \, \text{Im} \int_{\R} \phi_{\omega_k} (x_k^t) e^{-i \Gamma_k} \varepsilon_u \\
& \midspace  + \Upsilon_k^{\omega, \varepsilon } + \text{Im} \int_{\R} \phi_{\omega_k} (x_k^t) e^{-i \Gamma_k} \varepsilon_u \varepsilon_n + \dot{\omega}_k (t) \, \text{Re} \int_{\R} \Lambda_{\omega_k} (x_k^t) e^{-i \Gamma_k} \varepsilon_u \\
& \midspace  - ( \, \boxed{c_k} + \dot{\sigma}_k (t))   \text{Re} \int_{\R} \phi_{\omega_k} ' (x_k^t) e^{-i \Gamma_k} \varepsilon_u + \left ( \dot{\gamma}_k (t) + \omega_k^0 \,  \boxed{- \frac{c_k^2}{4}} \, \right ) \text{Im} \int_{\R} \phi_{\omega_k} (x_k^t) e^{-i \Gamma_k} \varepsilon_u. 
\end{align*}
The framed terms cancel one another. We find
\begin{equation}
\begin{split}
& \frac{\text{d}}{\text{d}t} \Ree \int_{\R} \phi_{\omega_k (t)} (x-c_kt - \sigma_k (t)) e^{-i \Gamma_k (t  , x)} \varepsilon_u (t \, , x) \, \text{d}x \\
=& - \frac{\sqrt{1-c_k^2}}{\sqrt{\omega_k}} \left [  \mathbf{d}^\omega \dot{\omega}_k (t) + \sum\limits_{j=1}^K \left ( \dot{\omega}_j (t) \Theta_{k,j}^{\omega,\omega} + \dot{\sigma}_j (t)   \Theta_{k,j}^{\omega,\sigma} + \left ( \dot{\gamma}_j (t) - \omega_j (t) + \omega_j^0 \right ) \Theta_{k,j}^{\omega,\gamma} \right ) \right. \\
& \midspace \midspace \midspace \midspace \midspace \left. - \left ( \Theta_k^{\omega,q} + \Theta_k^{\omega , \Upsilon} \right ) \right ]. 
\end{split}
\label{explomega}
\end{equation}
Using the first $K$ lines of the matricial system \eqref{modulEDO} , we get, for $k \in \{ 1 \, , ... \, , K \}$,
\[ \frac{\text{d}}{\text{d}t} \Ree \int_{\R} \phi_{\omega_k (t)} (x-c_kt - \sigma_k (t)) e^{-i \Gamma_k (t  , x)} \varepsilon_u (t \, , x) \, \text{d}x = 0. \]
Hence, for all $k$, the function $\text{Ort}_k^\omega \, : \, t \mapsto \text{Re} \int_{\R} \phi_{\omega_k (t)} (x_k^t) e^{-i \Gamma_k (t)} \varepsilon_u (t) \, \text{d}x$ is constant on $[T_* \, , T_p]$. Since $\varepsilon_u (T_p) = 0$, $\text{Ort}_k^\omega (T_p) = 0$ and then $\text{Ort}_k^\omega = 0$ on $[T_* \, , T_p ]$. \\
\\ Now, using the relation $L_{-,k,t} \left ( x_k^t \phi_{\omega_k} (x_k^t) \right ) = -2 \phi_{\omega_k} ' (x_k^t)$ and the lines $K+1$ to $2K$ in the matricial system \eqref{modulEDO}, we follow the same proof as above and find that:
\begin{equation}
\begin{split}
& \frac{\text{d}}{\text{d}t} \Ree \int_{\R} (x-c_kt- \sigma_k (t)) \phi_{\omega_k (t)} (x-c_kt - \sigma_k (t)) e^{-i \Gamma_k (t  , x)} \varepsilon_u (t \, , x) \, \text{d}x \\
=& - \sqrt{\omega_k} \sqrt{1-c_k^2} \left [ \mathbf{d}^\sigma \dot{\sigma}_{k} (t) + \sum\limits_{j=1}^K \left ( \dot{\omega}_j (t) \Theta_{k,j}^{\sigma,\omega} + \dot{\sigma}_{j} (t) \Theta_{k,j}^{\sigma ,\sigma} + \dot{\sigma}_{j} (t) \Theta_{k,j}^{\sigma ,\sigma} + \left ( \dot{\gamma}_j (t) - \omega_j(t) + \omega_j^0 \right ) \Theta_{k,j}^{\sigma ,\gamma} \right ) \right. \\
& \midspace \midspace \midspace \midspace \midspace \midspace \left. - \left ( \Theta_k^{\sigma  , q} + \Theta_k^{\sigma  , \Upsilon} \right ) \right ] \\
=&  \, \, 0. 
\end{split}
\label{explsigma}
\end{equation}
Hence, for all $k$, the function $\text{Ort}_k^\sigma \, : \, t \mapsto \text{Re} \int_{\R} x_k^t \phi_{\omega_k (t)} (x_k^t) e^{-i \Gamma_k (t)} \varepsilon_u (t) \, \text{d}x$ is constant on $[T_* \, , T_p]$. Since $\varepsilon_u (T_p) = 0$, $\text{Ort}_k^\sigma (T_p) = 0$ and then $\text{Ort}_k^\sigma = 0$ on $[T_* \, , T_p ]$. \\
\\ Lastly, using the relation $L_{-,k,t} \left ( \Lambda_{\omega_k} (x_k^t) \right ) = - \phi_{\omega_k} (x_k^t) + 2 \left ( \phi_{\omega_k}^2 \Lambda_{\omega_k} \right ) (x_k^t)$ and the last $K$ lines of the matricial system \eqref{modulEDO}, we follow the same proof as above and find that:
\begin{equation}
\begin{split}
& \frac{\text{d}}{\text{d}t} \Imm \int_{\R} \Lambda_{\omega_k (t)} (x-c_kt - \sigma_k (t)) e^{-i \Gamma_k (t  , x)} \varepsilon_u (t \, , x) \, \text{d}x \\
=& \, \frac{\sqrt{1-c_k^2}}{\sqrt{\omega_k}} \left [  \mathbf{d}^\gamma \left ( \dot{\gamma}_k (t) - \omega_k(t) + \omega_k^0 \right ) + \sum\limits_{j=1}^K \left ( \dot{\omega}_j (t) \Theta_{k,j}^{\gamma,\omega} + \dot{\sigma}_j (t)   \Theta_{k,j}^{\gamma,\sigma} + \left ( \dot{\gamma}_j (t) - \omega_j (t) + \omega_j^0 \right ) \Theta_{k,j}^{\gamma,\gamma} \right ) \right. \\
& \midspace \midspace \midspace \midspace \midspace \left. - \left ( \Theta_k^{\gamma , q} + \Theta_k^{\gamma , \Upsilon} \right ) \right ] \\
=& \, \, 0.
\end{split}
\label{explgamma}
\end{equation}
Hence, for all $k$, the function $\text{Ort}_k^\gamma \, : \, t \mapsto \text{Im} \int_{\R}  \Lambda_{\omega_k (t)} (x_k^t) e^{-i \Gamma_k (t)} \varepsilon_u (t) \, \text{d}x$ is constant on $[T_* \, , T_p]$. Since $\varepsilon_u (T_p) = 0$, $\text{Ort}_k^\gamma (T_p) = 0$ and then $\text{Ort}_k^\gamma = 0$ on $[T_* \, , T_p ]$. \\
\\ To control $| \dot{\omega}_k |$, $| \dot{\sigma}_k |$ and $| \dot{\gamma}_k - \omega_k + \omega_k^0 |$, we simply use \eqref{modulEDO} to write that $\dot{\Pi} (t) - \mathbf{N} ( \Pi - \Pi^0 ) = ( \mathbf{I}_{3K} + \mathbf{A} (t))^{-1} \mathbf{B} (t)$. From \eqref{modulAB} we know that $| ( \mathbf{I}_{3K} + \mathbf{A} (t))^{-1} \mathbf{B} (t) | \leqslant C || \varepsilon (t) ||_{\mathbf{H}} + C e^{-16 \theta_0 t}$ thus
\[ | \dot{\omega}_k (t) | +  | \dot{\sigma}_k (t) | +  | \dot{\gamma}_k (t) - ( \omega_k (t) - \omega_k^0 ) | \leqslant C || \varepsilon (t) ||_{\mathbf{H}} + C e^{-2 \theta_0 t} \]
for all $k \in \{ 1 \, , ... \, , K \}$ and all $t \in [T_* \, , T_p]$. \\
\\ \textit{(Proof that $T_*=t^*$).} Assume that $T_*>t^*$. Since \eqref{modulHYP} holds on $[T_* \, , T_p]$, the solution $\Pi (t)$ cannot blow up as $t \to T_*^+$ (take $T_0$ large enough), and continues to exist on $[T_* - \eta \, , T_p]$ for a certain $\eta >0$. Consequently, the hypothesis \eqref{modulHYP} must cease to hold for any $t$ close to $T_*^-$. Let us find a contradiction to this assumption. To facilitate the computations, we first prove the following lemma.

\begin{leftbar}
\noindent \textbf{Lemma 16.} Let $c \in (-1 \, , 1)$, $\omega , \omega^0 \in [\omega_- \, , \omega^+]$, $\sigma , \sigma^0 \in \R$ and $\gamma , \gamma^0 \in \R$. Denote $\pi = (\omega \, , \sigma \, , \gamma )$, $\pi^0 = (\omega^0 \, , \sigma^0 \, , \gamma^0 )$, $\Gamma^{\pi} (t \, , x) = \frac{c   x}{2} - \frac{c^2t}{4} + \omega^0 t + \gamma$ and $\Phi_\pi (t \, , x) = \phi_\omega (x-ct- \sigma ) e^{i \Gamma^\pi (t,x)}$. Then
\begin{align*}
& | \omega - \omega^0 | \leqslant C | \pi - \pi^0 |^2 + C \left | \text{Re} \int_{\R} e^{-i \Gamma^{\pi}} \phi_\omega (x-ct- \sigma ) ( \Phi_\pi - \Phi_{\pi^0} ) (t \, , x) \, \text{d}x \right | , \\
& | \sigma - \sigma^0 | \leqslant C | \pi - \pi^0 |^2 + C \left | \text{Re} \int_{\R} e^{-i \Gamma^{\pi}} (x-ct- \sigma ) \phi_\omega (x-ct- \sigma ) ( \Phi_\pi - \Phi_{\pi^0} ) (t \, , x) \, \text{d}x \right | \\
\text{and} \quad & | \gamma - \gamma^0 | \leqslant C | \pi - \pi^0 |^2 + C \left | \text{Im} \int_{\R} e^{-i \Gamma^{\pi}} \Lambda_\omega (x-ct- \sigma ) ( \Phi_\pi - \Phi_{\pi^0} ) (t \, , x) \, \text{d}x \right | .
\end{align*}
The constant $C$ depends on $\omega_-$ and $\omega^+$ but does not depend on $\omega$ itself. It also does not depend on $c$, $\sigma$ or~$\gamma$.
\end{leftbar}

\noindent \textit{Proof.} We begin with a Taylor expansion for the function $\pi = ( \omega \, , \sigma \, , \gamma ) \mapsto \Phi_\pi$:
\begin{equation} 
\left | \Phi_\pi - \Phi_{\pi^0} - ( \pi - \pi^0 ) \cdot \nabla_\pi \Phi ( \pi^0 ) \right | \leqslant C | \pi - \pi^0 |^2 
\label{Phi1}
\end{equation}
where the constant $C$ depends on the Hessian matrix of $\pi \mapsto \Phi_\pi$, thus depends on $\omega_-$ and $\omega^+$ but not on $\omega$ itself. We now express the gradient of $\pi \mapsto \Phi_\pi$. Denote $x^{t,\pi} = x-ct-\sigma$ (and $x^{t,\pi^0} = x-ct- \sigma^0$). We have
\[ \frac{\partial \Phi}{\partial \omega} ( \pi ) = \Lambda_\omega (x^{t,\pi}) e^{i \Gamma^\pi} , \quad \frac{\partial \Phi}{\partial \sigma} ( \pi ) = - \phi_\omega ' (x^{t,\pi}) e^{i \Gamma^\pi} \quad \text{and} \quad \frac{\partial \Phi}{\partial \gamma} ( \pi ) = i \phi_\omega (x^{t,\pi}) e^{i \Gamma^\pi} . \]
Previously, for the obtention of the orthogonality relations \eqref{explomega}, \eqref{explsigma} and \eqref{explgamma}, we only gave details for $\omega$; here let us conduct the proof for $\sigma$ (the others are similar). 
\begin{itemize}
	\item[$\bullet$] First, 
	\begin{align*}
	& \, \text{Re} \int_{\R} e^{-i \Gamma^\pi} x^{t,\pi} \phi_\omega (x^{t,\pi}) \frac{\partial \Phi}{\partial \omega} ( \pi^0 ) (t \, , x ) \, \text{d}x \\
	=& \, \text{Re} \int_{\R} e^{-i \Gamma^\pi} x^{t,\pi} \phi_\omega (x^{t,\pi}) \Lambda_{\omega^0} (x^{t,\pi^0}) e^{i \Gamma^{\pi^0}} \text{d}x \\
	=& \, \text{Re} \int_{\R} e^{-i \Gamma^{\pi^0}} x^{t, \pi^0} \phi_{\omega^0} (x^{t,\pi^0}) \Lambda_{\omega^0} (x^{t,\pi^0}) e^{i \Gamma^{\pi^0}} \text{d}x \\
	& \quad + \text{Re} \int_{\R} e^{-i \Gamma^\pi} (x^{t,\pi} \phi_\omega (x^{t,\pi}) - e^{-i \Gamma^{\pi^0}} x^{t, \pi^0} \phi_{\omega^0} (x^{t,\pi^0}) ) \Lambda_{\omega^0} (x^{t,\pi^0}) e^{i \Gamma^{\pi^0}} \text{d}x
	\end{align*}
	The first term vanishes by symmetry:
	\[ \text{Re} \int_{\R} e^{-i \Gamma^{\pi^0}} x^{t, \pi^0} \phi_{\omega^0} (x^{t,\pi^0}) \Lambda_{\omega^0} (x^{t,\pi^0}) e^{i \Gamma^{\pi^0}} \text{d}x = \int_{\R} x \phi_{\omega^0} \Lambda_{\omega^0} = 0. \]
	For the second term, we use a Taylor expansion for $\pi \mapsto e^{-i \Gamma^\pi} x^{t,\pi} \phi_\omega (x^{t,\pi})$, which gives
	\[ \left | e^{-i \Gamma^\pi} x^{t,\pi} \phi_\omega (x^{t,\pi}) - e^{-i \Gamma^{\pi^0}} x^{t, \pi^0} \phi_{\omega^0} (x^{t,\pi^0}) \right | \leqslant C | \pi - \pi^0 | . \]
	Again, the constant $C$ depends on the gradient of $\pi \mapsto e^{-i \Gamma^\pi} x^{t,\pi} \phi_\omega (x^{t,\pi})$, thus depends on $\omega_-$ and $\omega^+$ but not on $\omega$ itself. Gathering the estimates above, it follows that
	\[ \left | \text{Re} \int_{\R} e^{-i \Gamma^\pi} x^{t,\pi} \phi_\omega (x^{t,\pi}) \frac{\partial \Phi}{\partial \omega} ( \pi^0 ) (t \, , x ) \, \text{d}x \right | \lesssim C | \pi - \pi^0 | . \]
	\item[$\bullet$] Similarly, using the cancellation
	\[ \text{Re} \int_{\R} e^{-i \Gamma^{\pi^0}} x^{t,\pi^0} \phi_{\omega^0} (x^{t,\pi^0}) i \phi_{\omega^0} (x^{t,\pi^0}) e^{i \Gamma^{\pi^0}} \text{d}x = \text{Re} \int_{\R} i x \phi_{\omega^0}^2 = 0, \]
	we find
	\[ \left | \text{Re} \int_{\R} e^{-i \Gamma^\pi} x^{t,\pi} \phi_\omega (x^{t,\pi}) \frac{\partial \Phi}{\partial \gamma} ( \pi^0 ) (t \, , x ) \, \text{d}x \right | \lesssim C | \pi - \pi^0 | . \]
	\item[$\bullet$] For the derivative with regards to $\sigma$, we have
	\[ \text{Re} \int_{\R} e^{-i \Gamma^{\pi^0}} x^{t,\pi^0} \phi_{\omega^0} (x^{t,\pi^0}) (- \phi_{\omega^0}') (x^{t,\pi^0}) e^{i \Gamma^{\pi^0}} \text{d}x = - \int_{\R} x \phi_{\omega^0} \phi_{\omega^0} ' = \frac{1}{2} \phi_{\omega^0}^2 = 2 \sqrt{\omega^0} . \]
	It follows that
	\[ \left | \text{Re} \int_{\R} e^{-i \Gamma^\pi} x^{t,\pi} \phi_\omega (x^{t,\pi}) \frac{\partial \Phi}{\partial \sigma} ( \pi^0 ) (t \, , x ) \, \text{d}x - 2 \sqrt{\omega^0} \right | \leqslant C | \pi - \pi^0 | . \]
\end{itemize}
Combining the three estimates above, we get
\begin{equation}
\left | \text{Re} \int_{\R} e^{-i \Gamma^\pi} x^{t,\pi} \phi_\omega (x^{t,\pi}) (\pi - \pi^0) \cdot \nabla_\pi \Phi ( \pi^0 ) (t \, , x) \, \text{d}x - 2 \sqrt{\omega^0} ( \sigma - \sigma^0 ) \right | \leqslant C | \pi - \pi^0 |^2
\label{Phi2}
\end{equation}
where the constant $C$ depends on $\omega^-$ and $\omega^+$ but does not depend on $\omega$, $c$, $\sigma$ or $\gamma$. Combining \eqref{Phi1} and \eqref{Phi2} we immediately obtain
\begin{align*}
2 \sqrt{\omega^0} | \sigma - \sigma^0 | \leqslant & \, \left | \text{Re} \int_{\R} e^{-i \Gamma^\pi} x^{t,\pi} \phi_\omega (x^{t,\pi}) (\pi - \pi^0) \cdot \nabla_\pi \Phi ( \pi^0 ) (t \, , x) \, \text{d}x - 2 \sqrt{\omega^0} ( \sigma - \sigma^0 ) \right | \\
& \quad + \left | \text{Re} \int_{\R} e^{-i \Gamma^\pi} x^{t,\pi} \phi_\omega (x^{t,\pi}) \left ( \Phi_\pi - \Phi_{\pi^0} - ( \pi - \pi^0 ) \cdot \nabla_\pi \Phi ( \pi^0 ) \right ) (t \, , x) \, \text{d}x \right | \\
& \quad + \left | \text{Re} \int_{\R} e^{-i \Gamma^\pi} x^{t,\pi} \phi_\omega (x^{t,\pi}) ( \Phi_\pi - \Phi_{\pi^0} ) (t \, , x) \, \text{d}x \right | \\
\leqslant & \, C | \pi - \pi^0 |^2 + C | \pi - \pi^0 |^2 + \left | \text{Re} \int_{\R} e^{-i \Gamma^\pi} x^{t,\pi} \phi_\omega (x^{t,\pi}) ( \Phi_\pi - \Phi_{\pi^0} ) (t \, , x) \, \text{d}x \right | .
\end{align*}
The second estimate stated in the lemma follows. The other estimates of the lemma result from analogous proofs and the relations
\begin{align*}
& \text{Re} \int_{\R} e^{-i \Gamma^{\pi^0}} \phi_{\omega^0} (x^{t,\pi^0}) \Lambda_{\omega^0} (x^{t,\pi^0}) e^{i \Gamma^{\pi^0}} \text{d}x = \int_{\R} \phi_{\omega^0} \Lambda_{\omega^0} = \frac{1}{\sqrt{\omega^0}} , \\
& \text{Re} \int_{\R} e^{-i \Gamma^{\pi^0}} \phi_{\omega^0} (x^{t,\pi^0}) ( - \phi_{\omega^0} ' ) (x^{t,\pi^0}) e^{i \Gamma^{\pi^0}} \text{d}x = - \int_{\R} \phi_{\omega^0} \phi_{\omega^0} ' = 0 , \\
& \text{Re} \int_{\R} e^{-i \Gamma^{\pi^0}} \phi_{\omega^0} (x^{t,\pi^0}) i \phi_{\omega^0} (x^{t,\pi^0}) e^{i \Gamma^{\pi^0}} \text{d}x = \text{Re} \int_{\R} i \phi_{\omega^0}^2 = 0 , \\
& \text{Im} \int_{\R} e^{-i \Gamma^{\pi^0}} \Lambda_{\omega^0} (x^{t,\pi^0}) \Lambda_{\omega^0} (x^{t,\pi^0}) e^{i \Gamma^{\pi^0}} \text{d}x = \text{Im} \int_{\R} \Lambda_{\omega^0}^2 = 0 , \\
& \text{Im} \int_{\R} e^{-i \Gamma^{\pi^0}} \Lambda_{\omega^0} (x^{t,\pi^0}) ( - \phi_{\omega^0}') (x^{t,\pi^0}) e^{i \Gamma^{\pi^0}} \text{d}x = - \text{Im} \int_{\R} \Lambda_{\omega^0} \phi_{\omega^0} ' = 0 \\
\text{and} \quad & \text{Im} \int_{\R} e^{-i \Gamma^{\pi^0}} \Lambda_{\omega^0} (x^{t,\pi^0}) i \phi_{\omega^0} (x^{t,\pi^0}) e^{i \Gamma^{\pi^0}} \text{d}x = \int_{\R} \Lambda_{\omega^0} \phi_{\omega^0} = \frac{1}{\sqrt{\omega^0}} .
\end{align*}
Similar proofs to the one conducted for $\sigma$ enable to conclude. \hfill \qedsymbol

\noindent \\ Now take $t \in [T_* \, , T_p]$ and $k \in \{ 1 \, , ... \, , K \}$. We apply Lemma 16 to $\pi = \pi_k = ( \omega_k (t) \, , \sigma_k (t) \, , \gamma_k (t))$ and $\pi^0 = \pi_k^0 = ( \omega_k^0 \, , \sigma_k^0 \, , \gamma_k^0)$. We give details for the estimates on $\sigma$; the others are similar. By definition,
\[ \Phi_{\pi_k} = \frac{1}{\sqrt{1-c_k^2}} S_k^u \quad \text{and} \quad \Phi_{\pi_k^0} = \frac{1}{\sqrt{1-c_k^2}} R_k^u \]
thus
\begin{align}
| \sigma_k (t) - \sigma_k^0 | \leqslant & \, C | \pi_k (t) - \pi_k^0 |^2 + C \left | \text{Re} \int_{\R} e^{-i \Gamma_k (t)} x_k^t \phi_{\omega_k (t)} (x_k^t) ( \Phi_{\pi_k (t)} - \Phi_{\pi_k^0} ) (t \, , x) \, \text{d}x \right | \nonumber \\
\leqslant & \, C | \pi_k (t) - \pi_k^0 |^2 + C \left | \text{Re} \int_{\R} e^{-i \Gamma_k (t)} x_k^t \phi_{\omega_k (t)} (x_k^t) (S_k^u - R_k^u ) (t \, , x) \, \text{d}x \right | \label{x1}
\end{align}
where the constant $C$ depends on $c_1,...,c_K,\omega_-$ and $\omega^+$ (but not on any time-dependent parameter). Now we remark that, by definition, 
\[ S^u - R^u = (u - \varepsilon_u ) - (u - U_p) = U_p - \varepsilon_u . \]
It follows that
\begin{align*} 
& \, \left | \text{Re} \int_{\R} e^{-i \Gamma_k (t)} x_k^t \phi_{\omega_k (t)} (x_k^t) (S_k^u - R_k^u ) (t \, , x) \, \text{d}x \right | \\
\leqslant & \, \left | \text{Re} \int_{\R} e^{-i \Gamma_k (t)} x_k^t \phi_{\omega_k (t)} (x_k^t) (S^u - R^u) (t \, , x) \, \text{d}x \right | + \sum\limits_{\substack{j=1 \\ j \neq k}}^K \left | \text{Re} \int_{\R} e^{-i \Gamma_k (t)} x_k^t \phi_{\omega_k (t)} (x_k^t) (S_j^u - R_j^u)(t \, , x) \, \text{d}x \right | \\
\leqslant & \, \left | \text{Re} \int_{\R} e^{-i \Gamma_k (t)} x_k^t \phi_{\omega_k (t)} (x_k^t) U_p (t \, , x) \, \text{d}x \right |  + \left | \text{Re} \int_{\R} e^{-i \Gamma_k (t)} x_k^t \phi_{\omega_k (t)} (x_k^t) \varepsilon_u (t \, , x) \, \text{d}x \right | \\
& \quad + \sum\limits_{\substack{j=1 \\ j \neq k}}^K \left | \int_{\R} e^{-i \Gamma_k (t)} x_k^t \phi_{\omega_k (t)} (x_k^t) (S_j^u - R_j^u)(t \, , x) \, \text{d}x \right | . 
\end{align*}
We estimate these three terms.
\begin{itemize}
	\item[$\bullet$] First, it follows from \eqref{bootstrap} that
	\[ \left | \text{Re} \int_{\R} e^{-i \Gamma_k (t)} x_k^t \phi_{\omega_k (t)} (x_k^t) U_p (t \, , x) \, \text{d}x \right | \leqslant \| x \phi_{\omega_k (t)} \|_{L^2} \| U_p \|_{L^2} \leqslant \omega_-^{-1/4} \| U_p \|_{L^2} \leqslant C A_0 e^{- \theta_0 t} . \]
	\item[$\bullet$] Next, from the orthogonality relation $\text{Ort}_k^\sigma =0$ (see \eqref{explsigma}),
	\[ \text{Re} \int_{\R} e^{-i \Gamma_k (t)} x_k^t \phi_{\omega_k (t)} (x_k^t) \varepsilon_u (t \, , x) \, \text{d}x = 0. \]
	\item[$\bullet$] Last, taking $T_0$ large enough, $\sigma_k (t)$ is far from $\sigma_j (t)$ and $\sigma_j^0$ (for $j \neq k$), hence
	\[ \sum\limits_{\substack{j=1 \\ j \neq k}}^K \left | \int_{\R} e^{-i \Gamma_k (t)} x_k^t \phi_{\omega_k (t)} (x_k^t) (S_j^u - R_j^u)(t \, , x) \, \text{d}x \right | \leqslant C e^{-16 \theta_0 t} . \]
\end{itemize}
Gathering these three estimates, it follows that
\begin{equation} 
\left | \text{Re} \int_{\R} e^{-i \Gamma_k (t)} x_k^t \phi_{\omega_k (t)} (x_k^t) (S_k^u - R_k^u ) (t \, , x) \, \text{d}x \right | \leqslant C A_0 e^{- \theta_0 t} . 
\label{x2}
\end{equation}
Combining \eqref{x1} and \eqref{x2}, we obtain
\[ | \sigma_k (t) - \sigma_k^0 | \leqslant C | \pi_k (t) - \pi_k^0 |^2 + C A_0 e^{- \theta_0 t} . \]
Proceeding analogously for $\omega_k (t)$ and $\gamma_k (t)$, we get
\begin{equation}
| \pi_k (t) - \pi_k^0 | \leqslant C_2 | \pi_k (t) - \pi_k^0 |^2 + C_3 A_0 e^{- \theta_0 t} 
\label{x3}
\end{equation}
where $C_2$ and $C_3$ do not depend on $t$, $A_0$ or $C_0$. Recall the constant $C_1$ from \eqref{RS}; we now finally define 
\[ C_0 = 2(1+2K(1+C_1)C_2) \]
(the constant in \eqref{modulHYP}). Take $T_0$ large enough (depending on $A_0$ and the constants $C_i$) such that $C_0 A_0 e^{- \theta_0 T_0} \leqslant \frac{1}{2C_3}$. From \eqref{modulHYP} it follows that
\begin{equation}
| \pi_k (t) - \pi_k^0 | \leqslant C_0 A_0 e^{- \theta_0 T_0} \leqslant \frac{1}{2C_3} .
\label{x4}
\end{equation}
Gathering \eqref{x3} and \eqref{x4} leads to
\begin{align*} 
& \frac{1}{2} | \pi_k (t) - \pi_k^0 | \leqslant C_2 A_0 e^{- \theta_0 t} \\
\text{therefore} \quad & | \Pi (t) - \Pi^0 | \leqslant 2KC_2 A_0 e^{- \theta_0 t}.
\end{align*}
It follows from \eqref{bootstrap} and \eqref{RS} that
\begin{align*}
\| \varepsilon (t) \|_{\mathbf{H}} \leqslant & \, \| (U_p \, , N_p \, , V_p) (t) - (S^u \, ,S^n \, , S^v) (t) + (R^u \, , R^n \, , R^v)(t) \|_{\mathbf{H}} \\
\leqslant & \, A_0 e^{- \theta_0 t} + C_1 | \Pi (t) - \Pi^0 | \\
\leqslant & \, (1+2KC_1C_2)A_0 e^{- \theta_0 t}.
\end{align*}
Eventually, for all $t \in [T_* \, , T_p]$,
\[ \| \varepsilon (t) \|_{\mathbf{H}} + | \Pi (t) - \Pi^0 | \leqslant 2KC_2A_0 e^{- \theta_0 t} + (1+2KC_1C_2) A_0 e^{- \theta_0 t} \leqslant \frac{1}{2} C_0 A_0 e^{- \theta_0 t} \]
by definition of $C_0$. Hence, for $\eta >0$ small enough, \eqref{modulHYP} still holds on $[T_* - \eta \, , T_p]$. This offers the contradiction we were looking for. As a consequence, $T_* = t^*$. \hfill \qedsymbol

\end{document}